%% file: Main_ROM_Coarse.tex
\documentclass[final]{siamltex}
\usepackage{todonotes}

\usepackage{verbatim}
\usepackage[linesnumbered,lined,boxed,commentsnumbered]{algorithm2e}
\usepackage{latexsym, graphicx, epsfig}
\usepackage{pgfplots}
\usepackage[notcite,notref]{showkeys}
\usepackage{showkeys}
\usepackage{url}
\usepackage{color}
\usepackage{array}
\usepackage{amsfonts,amsmath,amssymb}
\usepackage{verbatim}
\usepackage{setspace}
\usepackage{multirow}
\usepackage{slashbox}
\usepackage{algpseudocode}

\newtheorem{remark}{Remark}[section]

\newcommand{\beq}{\begin{equation}}
\newcommand{\eeq}{\end{equation}}
\newcommand{\beqq}{\begin{equation*}}
\newcommand{\eeqq}{\end{equation*}}
\newcommand{\beqas}{\begin{eqnarray*}}
\newcommand{\eeqas}{\end{eqnarray*}}
\newcommand{\bsp}{\begin{split}}
\newcommand{\eesp}{\end{split}}

\renewcommand{\div}{\mathop{\rm div}\nolimits}

\usepackage{amscd}
\usepackage{colordvi}
\usepackage{epstopdf}
\usepackage{hyperref}
\usepackage{subcaption}
\usepackage{cleveref}
\def\norm#1{\|#1\|}

\title{ A non-intrusive bi-fidelity reduced basis method for time-independent problems 
}
\author{ Jun Sur Richard Park\thanks{Department of Mathematical Sciences, Korea Advanced Institute of Science \& Technology. Korea. Email: 
pjss1223@kaist.ac.kr.}
           \and  Xueyu Zhu\thanks{Department of Mathematics, University of Iowa, Iowa City, IA 52246. USA. Email: xueyu-zhu@uiowa.edu.}}

\begin{document}

\maketitle 

\begin{abstract}
Scientific and engineering problems often involve parametric partial differential equations (PDEs), such as uncertainty quantification, optimizations, and inverse problems. However, solving these PDEs repeatedly can be prohibitively expensive, especially for large-scale complex applications.
To address this issue, reduced order modeling (ROM) has emerged as an effective method to reduce computational costs. However, ROM often requires significant modifications to the existing code, which can be time-consuming and complex, particularly for large-scale legacy codes.
Non-intrusive methods have gained attention as an alternative approach. However, most existing non-intrusive approaches are purely data-driven and may not respect the underlying physics laws during the online stage, resulting in less accurate approximations of the reduced solution.
In this study, we propose a new non-intrusive bi-fidelity reduced basis method for time-independent parametric PDEs. Our algorithm utilizes the discrete operator, solutions, and right-hand sides obtained from the high-fidelity legacy solver. By leveraging a low-fidelity model, we efficiently construct the reduced operator and right-hand side for new parameter values during the online stage. Unlike other non-intrusive ROM methods, we enforce the reduced equation during the online stage. In addition, 
the non-intrusive nature of our algorithm makes it straightforward and applicable to general nonlinear time-independent problems. We demonstrate its performance through several benchmark examples, including nonlinear and multiscale PDEs.
\end{abstract}



\pagestyle{myheadings}
\thispagestyle{plain}


\input Intro_v2

\input Setup

\input Examples


\input Summary

\section*{Acknowledgments}
XZ was supported by Simons Foundation.
\clearpage

\bibliographystyle{plain}
\bibliography{referencesROM_Coarse}


\end{document}

%% file: intro_v2.tex
\section{Introduction}
Partial differential equations (PDEs) are fundamental to many scientific and engineering problems, and they are widely used to model physical phenomena. However, solving PDEs numerically with high accuracy can be computationally expensive, particularly when they need to be repeatedly computed at a large number of parameters. Examples of such scenarios include uncertainty quantification \cite{chen2016model, zhu2017multi, zhu2014computational} and optimization problems \cite{forrester2007multi}. To tackle this challenge, reduced-order modeling (ROM) has emerged as a popular and successful technique \cite{quarteroni2015reduced,hesthaven2016certified,rozza2008reduced}.

ROM typically involves an offline and an online stage. The offline stage is performed once to construct a set of reduced basis, and the online stage is used to construct the reduced problem by projecting the full-order model onto the reduced approximation space. However, traditional projection-based methods often have a major limitation: they are intrusive and often require extensive modifications to the existing sophisticated legacy solver of the full-order model. Moreover, these methods can be computationally inefficient for nonlinear problems, as the online computational complexity of the nonlinear term remains dependent on the degrees of freedom (DOFs) of the full-order solution. Additional techniques are usually needed to improve the computational efficiency, such as the empirical interpolation method (EIM) and its variant \cite{barrault2004empirical,chaturantabut2010nonlinear,negri2015efficient,bonomi2017matrix}. 

To address these issues, non-intrusive reduced order modeling techniques have received lots of attention.
Many efforts have focused on efficiently constructing the surrogate of the high-fidelity reduced coefficients without the projection of the full-order model, including interpolation \cite{xiao2015non,xiao2017parameterized}, neural network \cite{hesthaven2018non,wang2019non}, Gaussian process regression \cite{guo2018reduced}, to name a few. Nonetheless, 
the methods are developed based on high-fidelity samples, which can  be computationally expensive if many high-fidelity samples are needed to construct an accurate surrogate.   To mitigate this issue, multi-fidelity approaches have gained popularity recently. These approaches leverage both high and low-fidelity models to further reduce the number of high-fidelity samples needed in the context of reduced order modeling.   For example, a non-intrusive two-grid finite element method 
is proposed to construct a linear correction map for the reduced coefficients obtained by the coarse mesh model to the reduced coefficients based on the fine mesh model \cite{chakir2009two}. In \cite{kast2020non}, multi-fidelity Gaussian process regression (GPR) is employed to approximate the reduced coefficients of high-fidelity model.  A nonintrusive bi-fidelity neural network is proposed to improve predictive capability by augmenting additional  input features extracted from low-fidelity solutions with a limited number of high-fidelity samples \cite{lu2021bifidelity}.

It is worth noting that the  aforementioned non-intrusive approaches are purely data-driven and do not necessarily respect the underlying physics laws at the reduced approximation space during the online stage, potentially resulting in less accurate approximations of the reduced solution.
Recent efforts have been made to address this issue by enforcing physics laws in the reduced approximation space. For example, in [6], a physics-reinforced neural network was developed to enforce the reduced model as a soft constraint in a non-intrusive reduced order modeling context  \cite{chen2021physics}. 
Particularly relevant for our context, a non-intrusive bi-fidelity technique has been proposed to efficiently construct reduced coefficients while enforcing the underlying reduced equation \cite{etter2023coarse}. This method builds reduced operators using a linear combination of high-fidelity operators  at parameters selected using low-fidelity operators.
However, to approximate the expansion coefficients of the reduced operators, they employ a random subsampling technique that still requires access to partial information from the high-fidelity operator matrices. Moreover, the approximation may fail due to the randomized algorithm employed.

Motivated by recent developments in non-intrusive reduced-order modeling techniques and multi-fidelity methods, we propose an efficient ROM algorithm to address the challenges for time-independent problems. Our method capitalizes on the information available from the legacy solvers of full-order models over fine and coarse meshes, including the solution, discrete operator, and the right-hand side. Building on this foundation, we make two key assumptions. First, we assume that the high-fidelity solution/operator at any parameter can be well approximated by a linear combination of the high-fidelity solution/operators, respectively, at a finite number of parameter points. The similar assumptions can be found in \cite{etter2023coarse,negri2015efficient}. Secondly, we assume that the low-fidelity solutions and operators can capture the essential behaviors and general trends of the high-fidelity solution and operator, respectively, in the parameter space. 
Based on these assumptions, we first leverage the low-fidelity model to inform the selection of parameter points, which are then used to construct reduced basis sets for the high-fidelity models. After that, we approximate the expansion coefficients of the high-fidelity reduced operators  by leveraging information learned from the low-fidelity operators. Once these expansion coefficients are computed, the reduced operator of the high-fidelity model can be efficiently constructed during the online stage. Importantly, the computational complexity of the online stage is independent of the degrees of freedom of the full-order model but instead depends on the degrees of freedom of the coarse model.
Notably, our algorithm does not require significant modifications to the existing legacy solver, rendering it straightforward to implement and applicable to both linear and nonlinear problems.


The rest of the paper is organized as follows:  we first discuss the setup of the  problem. The build blocks of our proposed reduced basis method are presented in Section 3. Then the complexity analysis is provided in Section 4. In Section 5, we demonstrate the  performance of our approach through several benchmark examples.

%% file: Setup.tex
\section{Problem setup}
In this section, we shall discuss the basics of reduced order models  \cite{quarteroni2015reduced,hesthaven2016certified,rozza2008reduced} for general linear and nonlinear parametric partial differential equations. 
\subsection{Linear parametric partial differential equations}
\label{sec:linear}
 We first consider the following time-independent linear parametrized partial differential equations: 
\beq
\label{eq:main}
{\mathcal L}(\mu) u(x,\mu) = f(x,\mu), \ \  x \in  \Omega \subset \mathbb{R}^n,
\eeq
where ${\mathcal L} (\mu) $ is the differential operator and some proper boundary conditions are given. 
The physical/uncertain parameter $\mu=(\mu_1, \mu_2, \dots, \mu_d)$ is from the parametric domain $D \subset \mathbb{R}^d$.
The solutions to (\ref{eq:main}) are often numerically approximated by solving the corresponding discrete system:
\beq
\label{eq:main_discrete}
L_h(\mu) u_h(\mu) = f_h(\mu),
\eeq
where $L_h(\mu)\in \mathbb{R}^{N_h\times N_h}$ represent the discrete operator and $u_h(\mu)\in\mathbb{R}^{N_h}$ is the approximated  solution. The number $N_h$ is typically very large to provide a sufficiently accurate approximation to $u(x,\mu)$ in (\ref{eq:main}), therefore solving such a problem with a large $N_h$ is typically computationally expensive, particularly when a large number of high-fidelity solver runs are needed.  This is particularly evident in applications involving uncertainty quantification, design and optimization, as well as inverse problems. 

Our goal is to efficiently approximate the solutions of the partial differential equations at a given parameter $\mu$.  The traditional projection-based reduced models make the assumption that  the truth solution $u_h(\mu)$ lives in a low-dimensional  manifold, which can be approximated by a reduced approximation space with a dimension $N_{rb}$ much smaller than $N_h$. 
 Then the solution $u_h(\mu)$ to (\ref{eq:main_discrete}) is approximated by the reduced solution as follows:
 \beq
 \label{eq:reconstructed}
 u_r(\mu) = Q u_{rb}(\mu),
 \eeq
where  $u_{rb}(\mu)\in \mathbb{R}^{N_{rb}}$ is the reduced coefficient vector and $Q \in \mathbb{R}^{N_h\times N_{rb}}$ is the reduced basis that span the reduced approximation space.
By Galerkin projection on the problem (\ref{eq:main}) on the reduced approximation space, we can get the following reduced system:
\beq
\label{eq:reduced}
L_{rb}(\mu) u_{rb}(\mu) = f_{rb}(\mu),
\eeq
where  the reduced operator $L_{rb}(\mu) \in \mathbb{R}^{N_{rb}\times N_{rb}}$ and the reduced right-hand side $L_{rb}(\mu) \in \mathbb{R}^{N_{rb}\times N_{rb}}$ are given as follows:
\beq
\label{eq:reduced2-operator}
L_{rb}(\mu)  = Q^T L_h(\mu)Q, \quad f_{rb}(\mu) =Q^T f_h(\mu),
\eeq
For linear problems of affine operators ${\mathcal L}(\mu)$, the reduced system (\ref{eq:reduced}) can be assembled and solved efficiently thanks to the reduced dimension $N_{rb} \ll N_h$ for $u_{rb}(\mu)$. 
More specifically, if the operator $L_h(\mu)$ is affine with respect to functions of $\mu$, it can be expressed as linear combination of parameter-independent operators:
\beq
\label{eq:affine}
L_{h}(\mu) = \displaystyle\sum_{q=1}^{M} L_q c_q(\mu) ,
\eeq
where $c_q: D \to \mathbb{R}$ is a function of parameter and $L_q \in \mathbb{R}^{N_h \times N_h}$ are the operator matrices. 
Then the corresponding reduced operator is
\beq
\label{eq:affine_reduced}
L_{rb}(\mu) = Q^T L_h(\mu) Q = \displaystyle\sum_{q=1}^{M} [Q^T L_q Q] \, c_q(\mu). 
\eeq
Since the matrices $Q^T L_q Q$ are independent of the parameter $\mu$, it can be precomputed in the offline stage so that for every new parameter $\mu$, we can evaluate the reduced operator efficiently, i.e., the computation complexity is independent of $N_h$ in the online stage. 

However, with a non-affine operator, the computation complexity of the assembly of the reduced matrix $L_{rb}(\mu)$ at the online stage can depend on the dimension of the high-fidelity problem, $N_h$. 
To further improve the efficiency, EIM and its variants are usually needed \cite{barrault2004empirical,chaturantabut2010nonlinear}, where one often needs major modifications to the source code, which can be time-consuming and sophisticated, particularly for large-scale complex applications.

To address this challenge, we develop a non-intrusive reduced-order model, particularly to efficiently construct the reduced operator $L_{rb}(\mu)$ during the online stage in this work. 

 We assume we can access the solutions, assembly matrices, and the right-hand side of our model from the existing legacy solver. We further assume that for any given parameter $\mu \in {\mathcal D}$,  the high-fidelity model operator $L_h(\mu) \in \mathbb{R}^{N_h\times N_h}$, and right-hand side $f_h(\mu) \in \mathbb{R}^{N_h}$ can be well approximated by a linear combination of the high-fidelity operators/right-hand side, respectively
\beq
\label{eq:LF_expansion}
L_h(\mu) \approx \sum\limits_{j=1}^{n_L} L_h(\mu^L_j) \, a_h^{(j)}(\mu), \ 
f_h(\mu) \approx \sum\limits_{j=1}^{n_f} f_h(\mu^f_j) \, b_h^{(j)}(\mu),
\eeq
by properly chosen parameter points  $\{\mu^L_j \in {\mathcal D}, j=1, \hdots,n_L\}$ and  $\{\mu^f_j \in {\mathcal D}, j=1, \hdots,n_f\}$. The similar assumptions are found in \cite{negri2015efficient, etter2023coarse}.
The expansion coefficients $a_h^{(j)}(\mu)$ and $b_h^{(j)}(\mu)$ in (\ref{eq:LF_expansion}) need to be approximated.

In this case,  reduced operators $L_{rb}(\mu) $ and reduced right-hand side $f_{rb}(\mu) $ can be approximated as follows:
\beq
\begin{split}
\label{eq:QLF_expansion}
L_{rb}(\mu) &= Q^T L_h(\mu) Q \approx \sum\limits_{j=1}^{n_L} [Q^T L_h(\mu^L_j) Q ]\, a_h^{(j)}(\mu) := \sum\limits_{j=1}^{n_L} L_{rb}(\mu^L_j) \, a_h^{(j)}(\mu) , \\ 
f_{rb}(\mu) &=Q^T f_h(\mu) \approx \sum\limits_{j=1}^{n_f} [Q^T f_h(\mu^f_j)] \, b_h^{(j)}(\mu) := \sum\limits_{j=1}^{n_f} f_{rb}(\mu^f_j) \, b_h^{(j)}(\mu),
\end{split}
\eeq
where the matrices $L_{rb}(\mu^L_j)$ and vectors $f_{rb}(\mu^f_j)$ are independent of the parameter $\mu$, thus can be precomputed during the offline stage. Then in the online stage, the expansion coefficients $a_h^{(j)}(\mu)$ and $b_h^{(j)}(\mu)$ are computed to assemble the reduced model. 
In the next section, we introduce the corresponding problem formulation for nonlinear partial differential equations. 

\subsection{Nonlinear parametric partial differential equation}
Next, we consider the following nonlinear time-independent parametrized partial differential equations: 
\beq
\label{eq:main_nl}
{\mathcal R}(u(x,\mu);\mu) = 0, \ \  x \in  \Omega \subset \mathbb{R}^n,
\eeq
where ${\mathcal R} (\cdot;\mu) $ is the differential operator and proper boundary conditions are given. The high-fidelity solutions are often obtained numerically by solving the following nonlinear discrete system:
\beq
\label{eq:main_nl_discrete}
{\mathcal R}_h(u_h(\mu);\mu) = 0,
\eeq
where $u_h(\mu) \in \mathbb{R}^{N_h}$.

For a given parameter $\mu$
, the nonlinear problems (\ref{eq:main_nl_discrete}) are typically solved using iterative methods given an initial guess $u_h^0(\mu)$ for the solution. In general, the $k$-th iteration of such methods can be written as 
the following linear system.
\beq
\label{eq:main_discrete_nl_iter}
A_h(u_h^{k-1}(\mu);\mu) u_h^{k}(\mu) = g_h(u_h^{k-1}(\mu);\mu),
\eeq
where $u_h^{k}(\mu)\in\mathbb{R}^{N_h}$ represent the solution at the $k$-th iteration, and $A_h(\cdot;\mu)\in \mathbb{R}^{N_h\times N_h}$ and $g_h(\cdot;\mu)$ depend on the iteration method used. To illustrate this, we present  two commonly used iterative methods, 
\begin{itemize}
    \item  \textbf{Picard's iteration}: At $k^{th}$ iteration step of Picard's iteration, the solution $u_h^{k}$ is found by solving a linearized equation of the form (\ref{eq:main_discrete_nl_iter}) that depends on the solution $u_h^{k-1}$ from the previous iteration \cite{langtangen2016solving,larson2013finite,efendiev2014generalized}.
    \item \textbf{Newton's method}: At $k^{th}$ iteration step, we find $u_h^k(\mu)$ such that
\beq
\bsp
\label{eq:main_discrete_nl_newton}
{\mathcal J}_h(u_h^{k-1}(\mu);\mu) \delta u_h(\mu) &= - {\mathcal R}_h(u_h^{k-1}(\mu);\mu),\\
u_h^{k}(\mu) &= u_h^{k-1}(\mu) + \delta u_h(\mu),
\end{split}
\eeq
where ${\mathcal J}_h(u_h^{k-1}(\mu);\mu)$ is so-called Jacobian matrix of the residual ${\mathcal R}_h(u_h^{k-1}(\mu);\mu)$ \cite{langtangen2016solving,larson2013finite,bonomi2017matrix}. 
Combining the two equations in (\ref{eq:main_discrete_nl_newton}), we obtain
\beq
\bsp
\label{eq:main_discrete_nl_newton2}
{\mathcal J}_h(u_h^{k-1}(\mu);\mu) u^{k}_h(\mu) &= - {\mathcal R}_h(u_h^{k-1}(\mu);\mu)+{\mathcal J}_h(u_h^{k-1}(\mu);\mu) u_h^{k-1}(\mu),\\
\end{split}
\eeq
which can be written in the form (\ref{eq:main_discrete_nl_iter}) with
\end{itemize}
\beq
\bsp
A_h(u_h^{k-1}(\mu);\mu) &= {\mathcal J}_h(u_h^{k-1}(\mu);\mu), \\
g_h(u_h^{k-1}(\mu);\mu) &= - {\mathcal R}_h(u_h^{k-1}(\mu);\mu)+{\mathcal J}_h(u_h^{k-1}(\mu);\mu) u_h^{k-1}(\mu).
\end{split}
\eeq

We note that, given $u_h^{k-1}(\mu)$ obtained from the previous step,  \eqref{eq:main_discrete_nl_iter} is linear in terms of the solution $u_h^{k}(\mu)$. The iterative procedure ends at the $K$th iteration when it satisfies a specific stopping criterion. The required iteration number often depends on the parameter $\mu$, and we denote the last iteration step as $K(\mu)$. Then the high-fidelity numerical solution of (\ref{eq:main_nl}) is given by $u_h(\mu) = u_h^{K(\mu)}(\mu)$, the solution of 
\beq
\label{eq:main_discrete_nl_last}
A_h(u_h^{K(\mu)-1}(\mu);\mu) u_h(\mu) = g_h(u_h^{K(\mu)-1}(\mu);\mu).
\eeq
In this work, we consider on the above linearized equation at the last iteration step. 
For nonlinear problems, we have the following assumption that is similar to (\ref{eq:LF_expansion}) for any given parameter $\mu$ in the parameter domain.
\beq
\bsp
\label{eq:LF_expansion_nl}
A_h(u_h^{K(\mu)-1}(\mu);\mu) &\approx \sum\limits_{j=1}^{n_L} A_h(u_h^{K(\mu^L_j)-1}(\mu^L_j);\mu^L_j) \, a_h^{(j)}(\mu), \\ 
g_h(u_h^{K(\mu)-1}(\mu);\mu) &\approx \sum\limits_{j=1}^{n_f} g_h(u_h^{K(\mu^f_j)-1}(\mu^f_j);\mu^f_j) \, b_h^{(j)}(\mu),
\end{split}
\eeq
For simplicity, we denote the linearized operators and the right-hand sides as
\beq
\label{eq:simplify_hf}
L_h(\mu) := A_h(u_h^{K(\mu)-1}(\mu);\mu), \ f_h(\mu) := g_h(u_h^{K(\mu)-1}(\mu);\mu),
\eeq
which is consistent with the notation for the linear operators and the right-hand sides in section \eqref{sec:linear}. 
Similar to \eqref{eq:main_discrete_nl_last} as a linear system, the corresponding reduced model is written as follows:
\beq
\label{eq:reduced_nl_last}
L_{rb}(\mu) u_{rb}(\mu) = f_{rb}(\mu),
\eeq
where $L_{rb}(\mu) = Q^T L_h(\mu)Q $ and $f_{rb}(\mu) = Q^T f_h(\mu)$ for a reduced basis matrix $Q \in \mathbb{R}^{N_h\times N_{rb}}$.

By adopting the unified notations, the construction of the reduced model of both linear and nonlinear problems can be formulated in the same form of \eqref{eq:QLF_expansion}. Therefore, both linear and nonlinear problems can be handled by our proposed method (discussed in the next section) in the same manner.

In (\ref{eq:QLF_expansion}), everything other than the expansion coefficients $a_h^{(j)}(\mu)$ and $b_h^{(j)}(\mu)$ are independent of $\mu$ and can be pre-computed in offline stage.
Thus, for efficient online computation, it is crucial to develop an efficient algorithm to recover the expansion coefficients. 
In \cite{etter2023coarse}, the authors suggest 
to recover the expansion coefficients $a_h^{(j)}(\mu)$ of the operators by solving least squares problems via (\ref{eq:LF_expansion}) but only with randomly sampled rows of the high-fidelity operators.  However, it requires access to partial information of the high-fidelity operator for the new given parameter $\mu$ and the approximation may fail due to the randomness. Alternatively,
an interpolation approach (Matrix discrete empirical interpolation method) is proposed in \cite{negri2015efficient} to deterministically sub-sample the entries of the high-fidelity operator matrices to efficiently compute the expansion coefficients. {However, the implementation 
is less black-box, depending on
the underlying high-fidelity legacy solver.
Motivated by these previous approaches, we propose a new non-intrusive bi-fidelity framework to approximate high-fidelity coefficients $a_h^{(j)}(\mu)$ and $b_h^{(j)}(\mu)$ in order to construct and solve the reduced order model efficiently and accurately without modification of the legacy solver.  In the next section,
We shall introduce the detailed algorithm of our proposed method.

\section{Bi-fidelity reduced-based methods}

To set up the stage for the discussion later on, we assume a coarse mesh model corresponding to the high-fidelity (fine-mesh) model (\ref{eq:main}) is available: 
\beq
\label{eq:coarse_main}
L_l(\mu) u_l(\mu) = f_l(\mu),
\eeq
where $L_l(\mu) \in \mathbb{R}^{N_l \times N_l}$, $u_l(\mu)$, $f_l(\mu) \in \mathbb{R}^{N_l}$ for $N_l \ll N_h$. 
The coarse mesh model can be solved much faster than the high-fidelity model due to less number of degrees of freedom.  Although the coarse model is less accurate, 
the low-fidelity functions $u_l(\mu)$, $L_l(\mu)$ and $f_l(\mu)$ can capture some important features of corresponding high-fidelity functions $u_h(\mu)$, $L_h(\mu)$ and $f_h(\mu)$ respectively.
We remark that for nonlinear problems, (\ref{eq:coarse_main}) represents the low-fidelity counterpart of the high-fidelity linearized model (\ref{eq:main_discrete_nl_last}), i.e.,
\beq
\label{eq:simplify_lf}
u_l(\mu) = u^{K(\mu)}_l(\mu), \ L_l(\mu) = A_l(u_l^{K(\mu)-1}(\mu);\mu), \ f_l(\mu) = g_l(u_l^{K(\mu)-1}(\mu);\mu).
\eeq


For a parameter set $\gamma = \{\mu_1, \mu_2, \dots, \mu_{n}\}$, we denote snapshot matrices of low-fidelity solutions, operators, and the right-hand sides evaluated at the points in $\gamma$ as 
\beq
\label{low-snapshots}
\begin{split}
&{\mathcal U}_{l}(\gamma) = [u_l(\mu_1)\ \ u_l(\mu_2)\ \ \dots \ \ u_l(\mu_{n})],\\
&{\mathcal L}_{l}(\gamma) = [\vec{L}_l(\mu_1)\ \ \vec{L}_l(\mu_2)\ \ \dots \ \ \vec{L}_l(\mu_{n})], \\
&{\mathcal F}_{l}(\gamma) = [f_l(\mu_1)\ \ f_l(\mu_2) \ \ \dots \ \ f_l(\mu_{n})],
\end{split}
\eeq
where $\vec{L}$ means the vectorized operators. For example, $\vec{L}_l(\mu_i)$ is a vector of size $N_l^2 \times 1$ containing all entries of the matrix $L_l(\mu_i)\in \mathbb{R}^{N_l \times N_l}$. But in this work, we only consider the nonzero entries of the matrices $L_l(\mu_i)$ to generate $\vec{L}_l(\mu_i)$ for both computational and memory efficiency as the matrix is typically sparse. As a result, the dimension of the vector $\vec{L}_l(\mu_i)$ is ${\mathcal O}(N_l)$, corresponding to the number of nonzero entries of $L_l(\mu_i)$.

Similarly, the snapshot matrices for the high-fidelity solutions can be denoted accordingly.
\beq
\label{high-snapshots}
\begin{split}
& {\mathcal U}_{h}(\gamma) = [u_h(\mu_1)\ \ u_h(\mu_2)\ \ \dots\ \ u_h(\mu_{n})],\\
& {\mathcal L}_{h}(\gamma) = [\vec{L}_h(\mu_1)\ \ \vec{L}_h(\mu_2)\ \ \dots\ \ \vec{L}_h(\mu_{n})], \\
&{\mathcal F}_{h}(\gamma) = [f_h(\mu_1)\ \ f_h(\mu_2)\ \ \dots \ \ f_h(\mu_{n})],
\end{split}
\eeq
 The dimension of the vector $\vec{L}_h(\mu_i)$ is ${\mathcal O}(N_h)$, when the operator is sparse.
In this work, we shall leverage both low-fidelity and high-fidelity models to construct the reduced system efficiently.
Consider a set of parameters $\Gamma = \{\mu_1, \mu_2, \dots,\mu_{n_p}\}$ dense enough in the parameter domain $D$, our proposed algorithm consists the following two stages: \\
{\bf Offline stage}:
\begin{enumerate}
\item Low-fidelity runs
\begin{itemize}
\item[$-$] Perform low-fidelity simulations at all parameters in $\Gamma$ to obtain ${\mathcal U}_{l}(\Gamma)$, ${\mathcal L}_{l}(\Gamma)$, and ${\mathcal F}_{l}(\Gamma)$ in \eqref{low-snapshots}. 
\item[$-$] Select important parameters from $\Gamma$ to construct $\gamma_u = \{\mu^u_1,\mu^u_2, \dots, \mu^u_{N_{rb}}\}$,  $\gamma_L = \{\mu^L_1,\mu^L_2, \dots, \mu^L_{n_L}\}$, $\gamma_f = \{\mu^f_1,\mu^f_2, \dots, \mu^f_{n_f}\}$ in \eqref{eq:LF_expansion} 
\end{itemize}
\item High-fidelity runs
\begin{itemize}
\item[$-$] Perform high-fidelity simulations at the selected parameter set in $\gamma_u$, $\gamma_L$, and $\gamma_f$ to construct ${\mathcal U}_{h}(\gamma_u)$, ${\mathcal L}_{h}(\gamma_L)$, and ${\mathcal F}_{h}(\gamma_f)$  in \eqref{high-snapshots}. 
\end{itemize}
\item Reduced basis construction
\begin{itemize}
\item[$-$]
 Construct the reduced basis $Q$ using high-fidelity snapshots matrix ${\mathcal U}_{h}(\gamma_u)$.
\item[$-$] Compute the basis for reduced  operator and reduced 
 right-hand side: $L_{rb}(\mu^L)$ and $f_{rb}(\mu^f)$ in \eqref{eq:QLF_expansion} for $\mu^L \in \gamma_L$, $\mu^f \in \gamma_f$
\end{itemize}
\end{enumerate}
{\bf Online stage}:
\begin{enumerate}
\setcounter{enumi}{3}
\item For a given parameter $\mu$, run the corresponding low-fidelity solver to find the expansion coefficients $a_l(\mu)$ and $b_l(\mu)$.
\item Construct and solve reduced problems
\begin{itemize}
\item[$-$] Assemble the reduced operators and reduced right-hand sides by \eqref{eq:QLF_expansion}
\item[$-$] Solve the reduced system \eqref{eq:reduced}.
\end{itemize}
\end{enumerate}

The schematic description of our algorithm is presented in Figure \ref{fig:flowchart}.  We shall discuss the details of each building block of the proposed method.
\begin{figure}
    \centering
    \includegraphics[width=\textwidth]{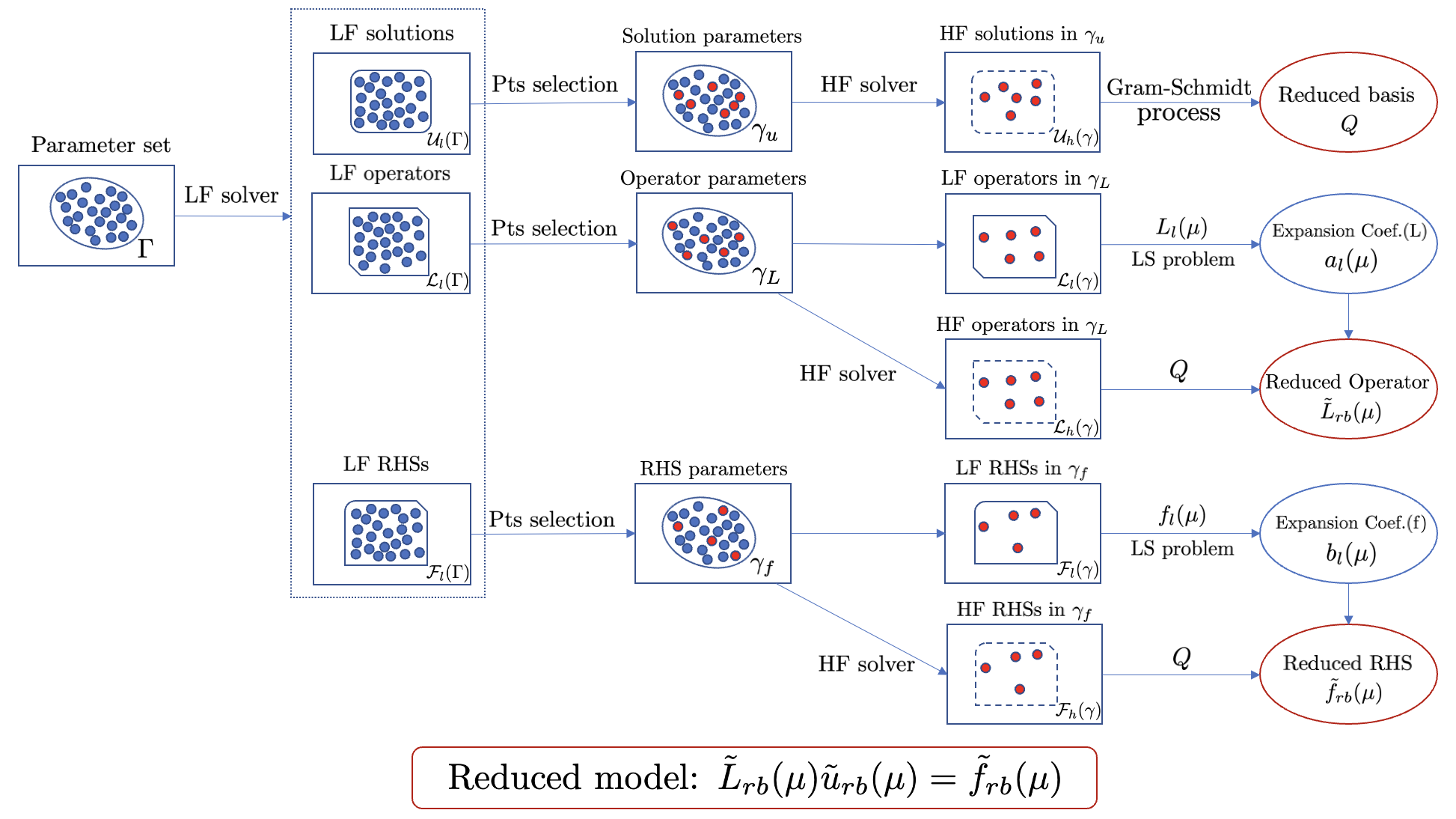}
    \caption{Flowchart of the proposed reduced method.}
    \label{fig:flowchart}
\end{figure}

\begin{remark}
According to the above algorithm, it can be seen that we only need the operator matrices, solution vectors, and the right-hand side vectors from the high- and low-fidelity legacy solvers.  
\end{remark}

\subsection{Point Selection and Reduced basis construction (Offline)}
We first choose a candidate set of parameters $\Gamma = \{\mu_1, \mu_2, \dots,\mu_{n_p}\} \subset D$ and assume  $\Gamma$ is large enough to represent  the parameter space. 
We consider the following matrix in $\mathbb{R}^{N_h\times n_p}$ of high-fidelity solutions at parameters in $\Gamma$,
\beq
\label{eq:S}
{\mathcal U}_{h}(\Gamma) = [u_h(\mu_1) \ \ u_h(\mu_2) \ \ \dots \ \ u_h(\mu_{n_p})]
\eeq
To construct the projection-based reduced basis  $Q \in \mathbb{R}^{N_h\times N_{rb}}$, one typically performs $n_p$ high-fidelity simulations to collect all solution snapshots and perform Gram-Schmidt process to extract the reduced basis set. However, the computation of the high-fidelity solutions $u_h(\mu_i)$ at many different parameters $\mu_i$, $i = 1,\dots, n_p$ is very expensive. To reduce the number of high-fidelity samples, we assume the coarse-model solutions can capture the main characteristics of the high-fidelity model in the parameter space. 
This idea is widely used in the multi-fidelity literature \cite{zhu2014computational,zhu2017multi,etter2023coarse}. We refer the interested readers for more details in \cite{etter2023coarse}.

Specially, we    
perform the low-fidelity simulations for (\ref{eq:coarse_main}) on 
 the candidate set $\Gamma$  to collect the following $N_l\times n_p$ snapshot matrix of the coarse solutions:
\beq
\label{eq:S_c}
{\mathcal U}_{l}(\Gamma) = [u_l(\mu_1) \ \ u_l(\mu_2) \ \ \dots \ \ u_l(\mu_{n_p})] 
\eeq
To select important parameter points, we perform pivoted Cholesky factorization on ${\mathcal U}_{l}(\Gamma)$.
This procedure provides a permutation of the matrix ${\mathcal U}_{l}(\Gamma)$. We select the columns corresponding to $N_{rb}$ largest diagonal entries of the matrix $R$, and choose the $N_{rb}$ parameters associated with the selected important columns, where $N_{rb}$ represents the dimension of reduced basis. We refer to the selected parameters as $\gamma_u = \{\mu^u_1,\mu^u_2, \dots, \mu^u_{N_{rb}}\} \subset \Gamma$. 
To construct the reduced basis for the fine mesh solutions,
we run the fine-mesh solver to obtain the high-fidelity solutions (\ref{eq:main}) {\it only} for each point in $\gamma_u$. We  denote the snapshot matrix as follows:
\beq
\label{eq:S_hat}
{\mathcal U}_{h}(\gamma) = [u_h(\mu^u_1) \ \ u_h(\mu^u_2) \ \ \dots \ \ u_h(\mu^u_{N_{rb}})],
\eeq
where ${\mathcal U}_{h}(\gamma)$ is   $N_h \times N_{rb}$ and the columns represent important solutions in ${\mathcal U}_{h}(\Gamma)$. 
To form the orthonormal reduced bases, we then perform Gram-Schmidt orthogonalization on ${\mathcal U}_{h}(\gamma)$ to obtain the reduced basis set $Q$. The detailed algorithm of this section is presented in Algorithm $1$.

\subsection{Constructing the reduced operator and right-hand side (Offline)}

{\bf Construct the reduced operator}. To select the important points for the reduced operator $L_h(\mu)$ in \eqref{eq:LF_expansion}, 
we consider the following snapshot matrix for the  fine mesh operator $L_h$ whose columns are the vectorized high-fidelity operators $L_h(\mu)$ in for candidate set $\Gamma$ in the parameter space.
\beq
\label{eq:L}
{\mathcal L}_{h}(\Gamma) = [\vec{L}_h(\mu_1) \ \ \vec{L}_h(\mu_2) \ \ \dots \ \ \vec{L}_h(\mu_{n_p})]
\eeq
if ${\mathcal L}_{h}(\Gamma)$ is available, we can perform pivoted Cholesky decomposition on it to select the important columns.  However, this requires a large number of high-fidelity runs to collect the snapshot matrix. By the same argument used to select important points for $u_h(\mu)$ in the previous section, we assume the coarse-model operator can capture the main characteristics of the high-fidelity operator in the parameter space. Based on this assumption, we construct the corresponding snapshot matrices of operators on a coarse mesh for the candidate set $\Gamma$:
\beq
\label{eq:L_c}
{\mathcal L}_{l}(\Gamma) = [\vec{L}_l(\mu_1) \ \ \vec{L}_l(\mu_2) \ \ \dots \ \ \vec{L}_l(\mu_{n_p})]
\eeq
Note that these matrices ${\mathcal L}_{l}(\Gamma)$ have already been assembled when we build the snapshot matrix ${\mathcal U}_{l}(\Gamma)$ in the previous section. 
We then apply pivoted Cholesky decomposition to ${\mathcal L}_{l}(\Gamma)$. The first $n_L$ elements of the permutation vector give the most important parameters for the operator $L$. We denote the selected parameters as 
\beq
\gamma_L = \{\mu^L_1, \mu^L_2, \dots , \mu^L_{n_L} \} \subset \Gamma.
\eeq
We then perform high-fidelity simulations for each point in $\gamma_L$ to compute the corresponding high-fidelity operators $\{ L_h(\mu^L_i) \ | \ \mu^L_i \in \gamma_L\}$. Then we project them onto the reduced approximation space spanned by the solution reduced basis set $Q$ and construct the reduced basis set $\{ Q^T L_h(\mu^L_i) Q \ | \ \mu^L_i \in \gamma_L \}$ for the reduced operator in (\ref{eq:QLF_expansion}).
During the online stage, the above reduced matrices will be used in the expansion (\ref{eq:QLF_expansion}) of the reduced operator $L_{rb}(\mu)$ at the given parameter $\mu$.

\begin{remark}
In many applications, the operator matrices are extremely sparse, and the number of nonzero entries of $L_h(\mu)$, $L_l(\mu)$ are ${\mathcal O}(N_h)$, ${\mathcal O}(N_l)$ respectively. We remark that to generate the vectorized operators $\vec{L}_h(\mu)$ and $\vec{L}_l(\mu)$, we only store the nonzero entries of the operator matrices $L_h(\mu)$ and $L_l(\mu)$. Thus, the dimensions of $\vec{L}_h(\mu)$ and $\vec{L}_l(\mu)$ are ${\mathcal O}(N_h)$, ${\mathcal O}(N_l)$ respectively.
\end{remark}

{\bf Construct reduced right-hand side}.
 To select the important points for the reduced right-hand side  in \eqref{eq:LF_expansion}, we denote  the following right-hand side snapshot matrix in $\mathbb{R}^{N_h\times n_p}$ as follows:
\beq
\label{eq:F}
{\mathcal F}_{h}(\Gamma) = [f_h(\mu_1) \ \ f_h(\mu_2) \ \ \dots \ \ f_h(\mu_{n_p})]
\eeq
Similarly, we denote the following matrix in $\mathbb{R}^{N_l\times n_p}$ of coarse right-hand sides instead of (\ref{eq:F}). 
\beq
\label{eq:Fc}
{\mathcal F}_{l}(\Gamma) = [f_l(\mu_1) \ \ f_l(\mu_2) \ \ \dots \ \ f_l(\mu_{n_p})]
\eeq

By the same argument for point selection for constructing the reduced operator, we shall select important points informed coarse mesh snapshot matrix ${\mathcal F}_{l}(\Gamma)$, which are already available as they were generated when we obtain ${\mathcal U}_{l}(\Gamma)$ using the coarse solver. We then use pivoted Cholesky decomposition on ${\mathcal F}_{l}(\Gamma)$ to find the corresponding $n_f$ dominant columns. Then we select the parameters corresponding to those columns and denote them by $\gamma_f =  \{\mu^f_1, \mu^f_2, \dots, \mu^f_{n_f} \} \subset \Gamma$.  
We compute the reduced basis $\{Q^T f_h(\mu^f) \ | \ \mu^f \in \gamma_f \}$ for the reduced right-hand side in (\ref{eq:QLF_expansion}). 
The detailed algorithm of this section can be found in Algorithm $2$.

\subsection{Least squares problems (Offline$\backslash$Online)}
In the online stage, the reduced operators and the right-hand sides for a new parameter $\mu$ will be approximated by the expansions (\ref{eq:QLF_expansion}). Thus, the expansion coefficients in (\ref{eq:QLF_expansion}) need to be computed. This is done by solving least-squares problems using (\ref{eq:LF_expansion}). In this section, we introduce the least-squares problem and its offline-online decomposition.  

To efficiently approximate the expansion coefficients in (\ref{eq:LF_expansion}):
\beq
\bsp
a_h(\mu) &= [a_h^{(1)}(\mu) \ \ \cdots \ \ a_h^{(n_L)}(\mu) ]^T\in \mathbb{R}^{n_L\times 1}, \\
b_h(\mu) &= [b_h^{(1)}(\mu) \ \ \cdots \ \ b_h^{(n_f)}(\mu) ]^T \in \mathbb{R}^{n_f\times 1},
\end{split}
\eeq
we consider the following matrix form corresponding to the equation (\ref{eq:LF_expansion}).
\beq
\label{eq:matrixlsq}
\bsp
\vec L_h(\mu) \approx {\mathcal L}_{h}(\gamma) a_h(\mu), \ F_h(\mu)\approx {\mathcal F}_{h}(\gamma)  b_h(\mu),
\end{split}
\eeq
where $\vec L_h(\mu) \in \mathbb{R}^{{\mathcal O}(N_{h})\times 1}$, 
and
\beq
\label{eq:HF_LF_gamma}
\begin{split}
{\mathcal L}_{h}(\gamma) &= [\vec{L}_h(\mu^L_1) \ \ \vec{L}_h(\mu^L_2) \ \ \cdots \ \ \vec{L}_h(\mu^L_{n_L})]\in\mathbb{R}^{{\mathcal O}(N_{h})\times n_L} \\
{\mathcal F}_{h}(\gamma) &= [f_h(\mu^f_1) \ \ f_h(\mu^f_2) \ \ \dots \ \ f_h(\mu^f_{n_f})]\in \mathbb{R}^{N_h\times n_L}.
\end{split}
\eeq
Given a parameter $\mu$, we aim to compute $a_h(\mu)$ and $b_h(\mu)$ using the least squares approach. However, the least squares problem on using (\ref{eq:matrixlsq}) is very expensive as we need to access the high-fidelity operator and right-hand side: $L_h(\mu)$, $f_h(\mu)$. 
 To mitigate this issue, we solve the following low-fidelity equations to compute the coefficients $a_l(\mu)$, $b_l(\mu)$ as the surrogates of $a_h(\mu)$, $b_h(\mu)$ respectively in \eqref{eq:matrixlsq}.
\beq
\label{eq:matrixlsq_c}
\bsp
\vec{L}_l(\mu)\approx {\mathcal L}_{l}(\gamma) a_l(\mu), \ f_l(\mu)\approx {\mathcal F}_{l}(\gamma) b_l(\mu),
\end{split}
\eeq
where the low-fidelity expansion coefficients are
\beq
\bsp
a_l(\mu) &= [a_l^{(1)}(\mu) \ \ \cdots \ \ a_l^{(n_L)}(\mu) ]^T\in \mathbb{R}^{n_L\times 1},\\
b_l(\mu) &= [b_l^{(1)}(\mu) \ \ \cdots \ \ b_l^{(n_f)}(\mu) ]^T\in \mathbb{R}^{n_f\times 1},
\end{split}
\eeq
and
\beq
\label{eq:LF_LF}
\begin{split}
{\mathcal L}_{l}(\gamma) &= [\vec{L}_l(\mu^L_1)) \ \ \vec{L}_l(\mu^L_2)) \ \ \cdots \ \ \vec{L}_l(\mu^L_{n_L}))] \in \mathbb{R}^{{\mathcal O}(N_{l})\times n_L},\\
{\mathcal F}_{l}(\gamma) &= [f_l(\mu^f_1) \ \ f_l(\mu^f_2) \ \ \dots \ \ f_l(\mu^f_{n_f})]  \in \mathbb{R}^{N_l\times n_f}.
\end{split}
\eeq
We note again that ${\mathcal L}_{l}(\gamma)$ and ${\mathcal F}_{l}(\gamma)$ were already obtained at offline stage when we assembled ${\mathcal U}_{l}(\gamma)$ (\ref{eq:S}).
We solve (\ref{eq:matrixlsq_c}) for $a_l(\mu)$, $b_l(\mu)$ using least squared approaches. Specifically, we solve the following problems:

\beq
\label{eq:normal_c}
\begin{split}
G^L_{l,\gamma} a_l(\mu) &= g^L_{l,\gamma}(\mu), \\ 
G^F_{l,\gamma}  b_l(\mu) &= g^F_{l,\gamma}(\mu),
\end{split}
\eeq
where 
\beq
\label{eq:normal_d}
\begin{split}
G^L_{l,\gamma} &= {\mathcal L}_{l}(\gamma)^T {\mathcal L}_{l}(\gamma), \quad g^L_{l,\gamma}(\mu) = {\mathcal L}_{l}(\gamma)^T \vec L_l(\mu), \\ 
G^F_{l,\gamma} &= {\mathcal F}_{l}(\gamma)^T {\mathcal F}_{l}(\gamma), \quad 
 g^F_{l,\gamma}(\mu) = {\mathcal F}_{l}(\gamma)^T f_l(\mu)
\end{split}
\eeq

For the least squares problem solve, we have the following offline-online decomposition.\\
{\bf Offline} Calculate the matrices $G^L_{l,\gamma} \in \mathbb{R}^{n_L\times n_L}$ and $G^F_{l,\gamma} \in \mathbb{R}^{n_f\times n_f}$; \\
{\bf Online} Compute $L_l(\mu)$ and $f_l(\mu)$ for a new parameter $\mu$ and solve the least squares problem (\ref{eq:normal_c})

The algorithm for the least-squares problem can be found in Algorithm $2$ and $3$.
The expansion coefficient vectors $a_l(\mu)$, $b_l(\mu)$ computed in this step will be used for the constructions of the reduced operators and right-hand sides in the next step. We note that the least squares problem is independent of our high-fidelity solver. Nonetheless, it still depends on  coarse solver, where we assume the operators and the right-hand sides based on coarse solver is much cheaper to assembly.


\begin{remark}
Note that the coarse solution $u_l(\mu)$ can be easily computed using the coarse solver. Although the coarse model is usually not an accurate approximation to  the underlying  problem, it is still useful to guide the construction of reduced bases, operators, and right-hand sides.
\end{remark}

\subsection{Assembly and solving reduced problem (Online)}
\label{subsec:rb_prob}
Once the expansion coefficients are found, we can approximate the reduced operator $L_{rb}(\mu)$ and reduce right-hand side $f_{rb}(\mu)$ for a given parameter $\mu$ using the right-hand sides of (\ref{eq:QLF_expansion}),
\beq
\label{eq:QLF_expansion_red}
\bsp
\tilde{L}_{rb}(\mu) = \sum\limits_{j=1}^{n_L} L_{rb}(\mu^L_j) \, a_l^{(j)}(\mu), \  
\tilde{f}_{rb}(\mu) = \sum\limits_{j=1}^{n_f} f_{rb}(\mu^f_j) \, b_l^{(j)}(\mu) 
\end{split}
\eeq
We note that all matrices in the right-hand side of (\ref{eq:QLF_expansion_red}) have been pre-computed in the offline stage.
This step does not depend on the size of our high-fidelity model and its complexity is independent of $N_h$.

Then we solve the following reduced problem to obtain the reduced coefficient $\tilde{u}_{rb}(\mu)$:
\beq
\label{eq:reduced_short}
\tilde{L}_{rb}(\mu) \tilde{u}_{rb}(\mu) = \tilde{f}_{rb}(\mu).
\eeq
Then the reconstructed reduced solution $u_r(\mu)$ defined in (\ref{eq:reconstructed}) can be approximated by $\tilde{u}_r(\mu) = Q \tilde{u}_{rb}(\mu) \in \mathbb{R}^{N_h}$.

\RestyleAlgo{ruled}
\SetKwComment{Comment}{/* }{ */}

\begin{algorithm}
\caption{Constructing reduced basis (Offline)}\label{alg:offlineQ}
\SetKwInOut{Input}{input}\SetKwInOut{Output}{output}
\Input{$\Gamma = \{\mu_1, \mu_2, \dots,\mu_{n_p}\}$, $N_{rb}$}
\Output{Reduced basis matrix $Q$}
\setstretch{1.35}
\tcc{Construct coarse-grid solutions ${\mathcal U}_{l}(\Gamma)$}
${\mathcal U}_{l}(\Gamma)$, ${\mathcal L}_{l}(\Gamma)$, ${\mathcal F}_{l}(\Gamma)$ $\longleftarrow$ LFSolver($\Gamma$)  \;
\tcc{Perform Cholesky factorization on ${\mathcal U}_{l}(\Gamma)$}
Permutation vector $P_U$ $\longleftarrow$ Cholesky$({\mathcal U}_{l}(\Gamma))$\;
\tcc{Select $N_{rb}$ indices ${\mathcal I}_u$ of important columns of ${\mathcal U}_{l}(\Gamma)$}
${\mathcal I}_u = P_{U}(1:N_{rb})$, $\gamma_u = \{\mu_i \in \Gamma \ | \ i \in {\mathcal I}_u\}$\;
\tcc{Construct ${\mathcal U}_{h}(\gamma)$}
${\mathcal U}_{h}(\gamma)$ $\longleftarrow$ HFSolver($\gamma_u$)\;
\tcc{Construct $Q$ by Gram-Schmidt orthogonalization on ${\mathcal U}_{h}(\gamma)$}
$Q$ $\longleftarrow$ G-S$({\mathcal U}_{h}(\gamma))$\;
\label{alg:alg_off1}
\end{algorithm}

\begin{algorithm}
\label{alg:alg_off2}
\caption{Constructing coarse and reduced operators, sources (Offline)}\label{alg:offlineLF}
\SetKwInOut{Input}{input}\SetKwInOut{Output}{output}
\Input{$\Gamma = \{\mu_1, \mu_2, \dots,\mu_{n_p}\}$, $Q$, ${\mathcal L}_{l}(\Gamma)$, ${\mathcal F}_{l}(\Gamma)$, $n_L$, $n_f$}
\Output{$G^L_{l,\gamma}$, $G^F_{l,\gamma}$, ${\mathcal L}_{l}(\gamma)$, ${\mathcal F}_{l}(\gamma)$, Reduced L-basis, Reduced F-basis} 
\setstretch{1.35}
\tcc{Construct a basis of reduced operator}
Permutation vector $P_L$ $\longleftarrow$ Cholesky$({\mathcal L}_{l}(\Gamma))$\;
${\mathcal I}_L = P_L(:,1:n_L)$, $\gamma_L = \{\mu_i \in \Gamma \ | \ i\in {\mathcal I}_L\}$ \;
${\mathcal L}_{l}(\gamma)$ $\longleftarrow$ ${\mathcal L}_{l}(\Gamma)$, $\gamma_L$\;
Reduced Operator basis (Reduced L-basis): $\{Q^T L_h(\mu_{i}) Q \ | \ \mu_{i}\in \gamma_L \}$ $\longleftarrow$ HFSolver($\gamma_L$)\;
\tcc{Construct a basis of reduced right-hand side}
Permutation vector $P_F$ $\longleftarrow$ Cholesky$({\mathcal F}_{l}(\Gamma))$\;
${\mathcal I}_f = P_F(:,1:n_f)$, $\gamma_f = \{\mu_i \in \Gamma \ | \ i\in {\mathcal I}_f\}$\;
${\mathcal F}_{l}(\gamma)$ $\longleftarrow$ ${\mathcal F}_{l}(\Gamma)$, $\gamma_f$\;
Reduced right-hand side basis (Reduced F-basis): $\{Q^T f_h(\mu_{i}) \ | \ \mu_i \in \gamma_f \}$$\longleftarrow$ HFSolver($\gamma_f$)\;
\tcc{Compute Gramian matrices for the least squares problems}
$G^L_{l,\gamma}$, $G^F_{l,\gamma}$ $\longleftarrow$ ${\mathcal L}_{l}(\gamma)^T {\mathcal L}_{l}(\gamma)$, ${\mathcal F}_{l}(\gamma)^T {\mathcal F}_{l}(\gamma)$\; 
\end{algorithm}

\begin{algorithm}
\label{alg:alg_on}
\caption{Construct and solve reduced problem(Online)}\label{alg:online}
\SetKwInOut{Input}{input}\SetKwInOut{Output}{output}
\Input{A parameter $\mu \in D$, $G^L_{l,\gamma}$, $G^F_{l,\gamma}$, ${\mathcal L}_{l}(\gamma)$, ${\mathcal F}_{l}(\gamma)$, Reduced L-basis, Reduced F-basis}
\Output{Reduced solution $u_r(\mu)$}
\setstretch{1.35}
$L_l(\mu)$, $f_l(\mu)$ $\longleftarrow$ LFSolver$(\mu)$\;
\tcc{Find expansion coefficients $a_l$, $b_l$}
$a_l$, $b_l$ $\longleftarrow$  Solve LS problem (\ref{eq:normal_c})\;
\tcc{Assembling reduced operator and right-hand side}
$\tilde{L}_{rb}(\mu)$ $\longleftarrow$ Reduced L-basis, $a_l$\;
$\tilde{f}_{rb}(\mu)$ $\longleftarrow$ Reduced F-basis, $b_l$\;
\tcc{Solving the reduced problem}
$\tilde{u}_{rb}(\mu)$ $\longleftarrow$ $\tilde{L}_{rb}(\mu) \backslash \tilde{f}_{rb}(\mu)$\;
\tcc{Reconstructing high-fidelity reduced solution} 
$\tilde{u}_r(\mu)$ $\longleftarrow$ $Q \tilde{u}_{rb}(\mu)$
\end{algorithm}

\section{Computational Cost}
In this section, we summarize the cost associated with the proposed algorithm. 

Offline cost:
\begin{itemize}
    \item $n_p$ low-fidelity runs to collect the snapshot matrices.
    \item Identifying important parameters ($\gamma_u$, $\gamma_L$, $\gamma_f$) by pivoting Cholesky decomposition: ${\mathcal O}(n_p N_{rb}^2) + {\mathcal O}(n_p n_L^2)+{\mathcal O}(n_p n_f^2)$ 

    \item $N_{rb}+n_L+n_f$ high-fidelity runs. 

    \item Performing Gram-Schmidt orthogalization on high-fidelity solution snapshot matrix ${\mathcal U}_{h}(\gamma)$ (\ref{eq:S_hat}) to construct the reduced basis $Q$: ${\mathcal O}(N_h N_{rb}^2)$ 
    \item Assembling reduced basis for operator, right-hand sides, $L_{rb}(\mu^L_j)$, $f_{rb}(\mu^f_j)$ in (\ref{eq:QLF_expansion}):${\mathcal O}(n_L(N_{rb}N_h^2+N_{rb}^2 N_h))+{\mathcal O}(n_f N_{rb} N_h)$ \\

\end{itemize}

Online cost:
\begin{itemize}
    \item For a given  new parameter $\mu$, run a low-fidelity solver to obtain $g^L_{l,\gamma}(\mu)$, $g^F_{l,\gamma}(\mu)$ (\ref{eq:normal_d}) for the least squares problem : 
    one low-fidelity run + ${\mathcal O}((n_L+n_f) N_{l})$
    \item Solving the least squares problems to the low-fidelity coefficients for  reduced operator and right-hand side in \eqref{eq:normal_c}: ${\mathcal O}(n_L^3)$
    \item Assembling reduced operators ${\mathcal O}(n_L N_{rb}^2)$ and right-hand sides ${\mathcal O}(n_f N_{rb})$ in \eqref{eq:QLF_expansion_red}.
    \item Solving the reduced system  in \eqref{eq:reduced_short}: ${\mathcal O}(N_{rb}^3)$.
\end{itemize}
Note that the offline stage is performed only once, and the major cost of the offline stage arises from the
$N_{rb}+n_L+n_f$ number of high fidelity simulations. The online cost is dominated by the low-fidelity solver.

%% file: Examples.tex
\section{Numerical examples}
In this section, we present numerical results to validate the robustness of our proposed approach for some benchmark examples,  including nonlinear problems and the problems featuring multiscale high-contrast media. 

To quantify the accuracy of our method, we compute the following errors: 
\begin{equation}
e_{u}(\mu) = \frac{\norm{\tilde{u}_r(\mu)  -u_{h}(\mu)}_{l^2}}{\norm{ u_{h}(\mu)}_{l^2}}, \ \ 
e_{u}^l(\mu) = \frac{\norm{u_{l}(\mu)  -u_{h}(\mu)}_{l^2}}{\norm{ u_{h}(\mu)}_{l^2}}, \ \ 
\end{equation}
where $\tilde{u}_{r}$ is our reduced solution and
 the reference solution $u_h$ represents the high-fidelity solution, and $u_l$ is the solution of the low-fidelity models used for each example. We report the errors averaged over an independent test data set of parameters that has the same size as the dense parameter set $\Gamma$.
 In all examples, the parameters in the dense parameter set $\Gamma$ and the test data set are randomly sampled 
 in the parameter domain.

Throughout the example sections, we compare our method with an existing bi-fidelity method where the reduced equation in section \ref{subsec:rb_prob} is not utilized.  The reference bi-fidelity method, described in \cite{narayan2014stochastic, zhu2014computational}, employs pivoted Cholesky decomposition to select $N_{rb}$ significant parameters from the low-fidelity solution data to form the reduced basis set for the high-fidelity model. For a new parameter $\mu$, the method determines the expansion coefficients $c_l(\mu)$ by solving least squares problems utilizing the low-fidelity solution data. These low-fidelity coefficients $c_l(\mu)$ are then used to approximate the high-fidelity solution at $\mu$ through the expression:
\begin{equation}
u_h(\mu) \approx \sum_{j=1}^{N_{rb}} c_l(\mu) u_h(\mu_j),
\end{equation}
where $\mu_j$ represents the selected parameters. For more comprehensive details, we refer interested readers to \cite{narayan2014stochastic, zhu2014computational}.
In contrast to this method, one of the key features of our proposed approach is the incorporation of equation information (\ref{eq:reduced_short}) during the online stage, which is not encompassed in the aforementioned reference method. Thus, we refer to this method as the {\it reference bi-fidelity method} and our proposed method as the {\it proposed bi-fidelity method}.

We used MATLAB implementations for both reference and proposed solutions in all examples. All implementations are carried out on a MacBook pro using Apple M1 pro chip with a 10-core CPU.


\subsection{2d nonlinear elliptic equation}
\label{sec:2DNL}
We first consider the following nonlinear elliptic equations defined in $\Omega = [-\frac{\pi}{2},\frac{\pi}{2}]^2$ inspired by the first example in \cite{hesthaven2018non}.
\beq
\label{eq:2dnonlinear}
\begin{split}
-\nabla \cdot [(2+\sin(2\pi \mu_2 u(x,\mu)+\mu_1)) \nabla u(x,\mu)] = f(x,\mu), \ \  &\textrm{in} \ \ \Omega,
\end{split}
\eeq
where $x = (x_1,x_2)$ and the source is given by $f(x,\mu) = \frac{1}{1+\mu_3^2}\sin(4x_1) + \mu_2 x_2$. The Dirichlet boundary condition is imposed and the parameters $\mu = 
(\mu_1,\mu_2,\mu_3)$ are chosen from $[0,1]^3$.



For high-fidelity model, we use a FEM solver over equal-distant meshes of $128\times 128$ points over the spatial domain. The low-fidelity model is based on the same FEM solver on a coarse mesh with a mesh of $8\times 8$. The FEM solver utilizes the Picard's iteration where the $k+1^{th}$ iteration step is as follows:
\beq
-\nabla \cdot [\kappa(u^{k},\mu)\nabla u^{k+1}] = f(x,\mu),
\eeq
where $\kappa(u,\mu) = 2+\sin(2\pi \mu_2 u+\mu_1)$.
The dense point set $\Gamma$
is chosen to be $512$ points randomly sampled in the parameter space. 
$N_{rb}$ important points are selected based on low-fidelity solutions over the set $\Gamma$, where high-fidelity solutions are collected to form the high-fidelity reduced bases. The same strategy is applied to select $n_L = 30$ and $n_f = 2$ important parameter points to form the basis set for high-fidelity operators and right-hand sides, in \eqref{eq:QLF_expansion}. In summary, total $N_T= N_{rb}+n_L+n_f$ high-fidelity runs are required at most during the offline stage. We remark that this number is reduced avoiding duplicated high-fidelity runs when some parameter points are selected multiple times for reduced basis of the solution, operator or the right-hand side.

\begin{figure}[!htb]
	\centering
 \begin{subfigure}{0.45\textwidth}
\includegraphics[width=\textwidth]{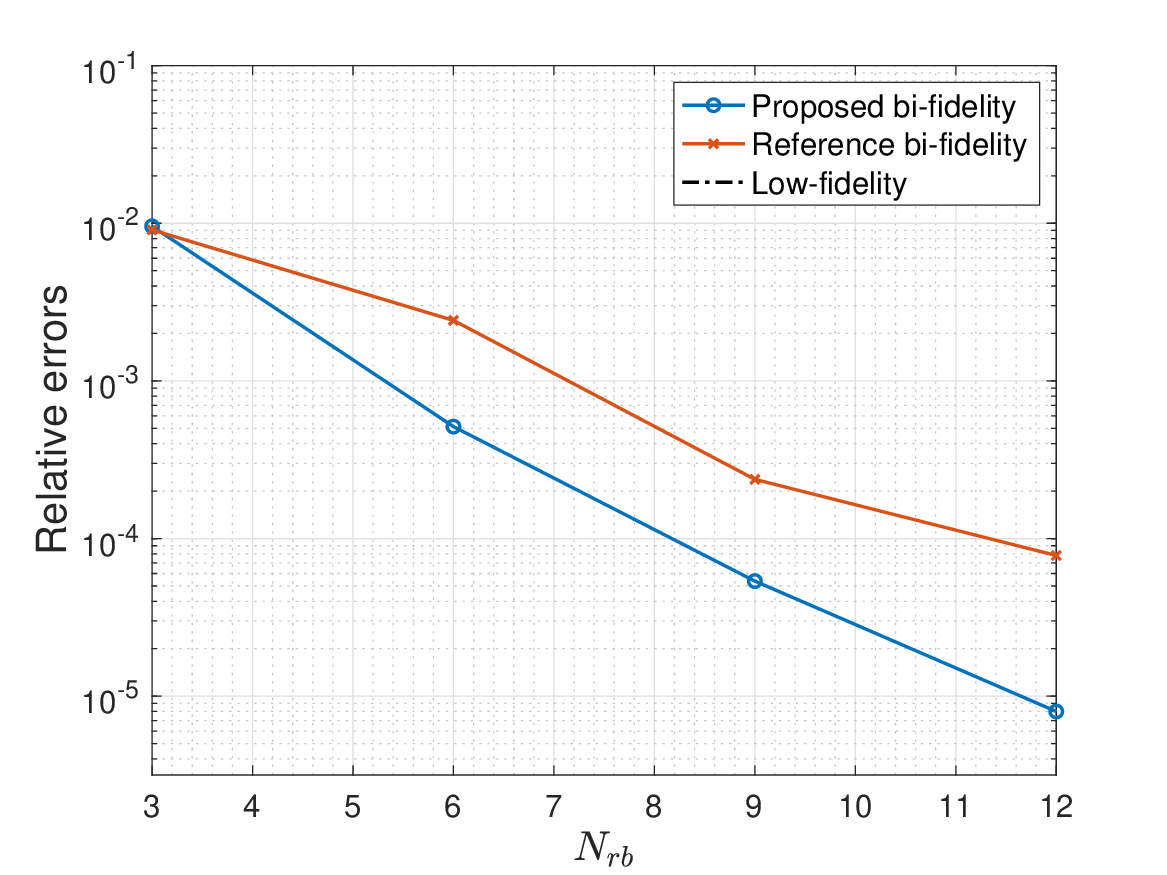}
\caption{}
    \label{fig:Errors_2DNL}
    \end{subfigure}
    \hspace{.1in}
 \begin{subfigure}{0.45\textwidth}
  \includegraphics[width=\textwidth]{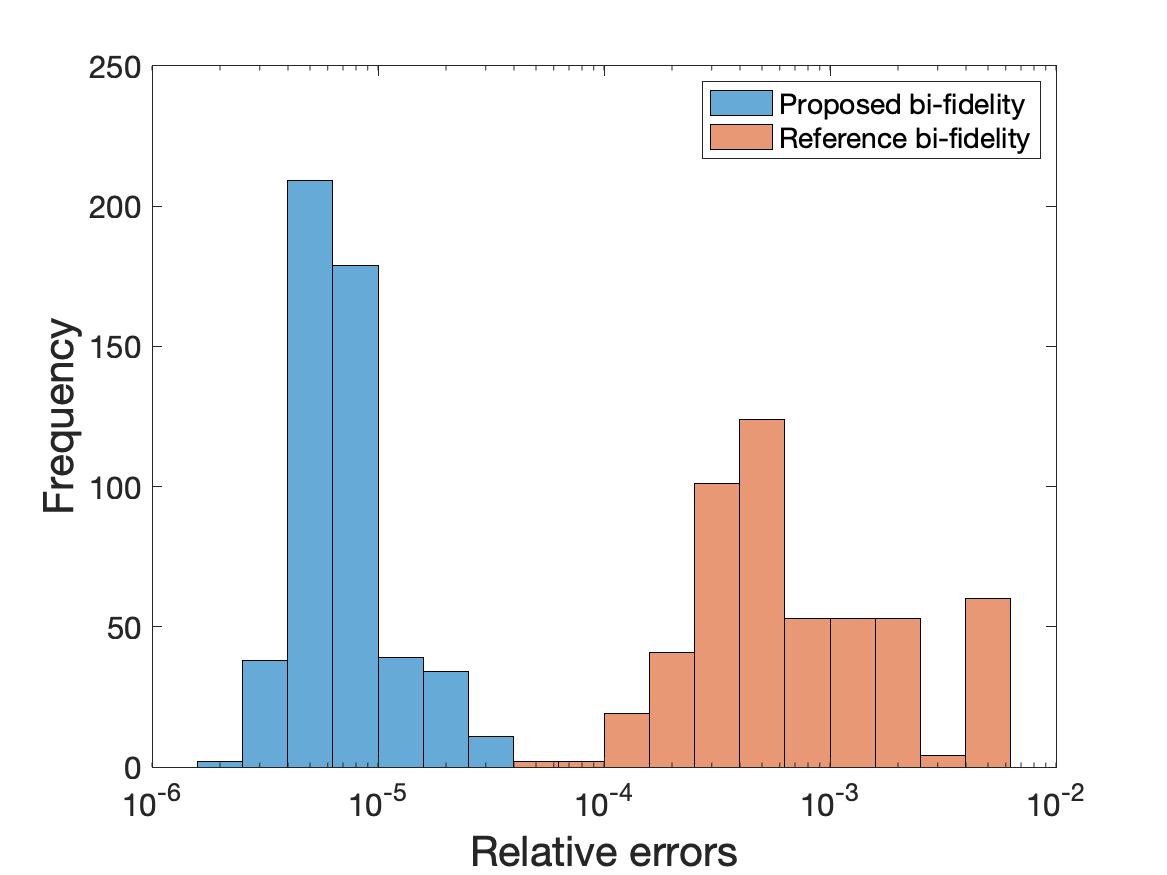}
  \caption{}
  \label{fig:Histo_2DNL}
\end{subfigure}
\vspace{-.2in}
  \caption{Example \ref{sec:2DNL}: a. The relative errors of the reduced solution by our new bi-fidelity method and the reference bi-fidelity method with different number of reduced basis; b. Histogram of the relative errors for both new and the reference bi-fidelity methods when $N_{rb} = 12$.}
  \label{figure:Errors_2D_NL_new2}
 \end{figure}

Figure \ref{figure:Errors_2D_NL_new2} represents the relative errors of our bi-fidelity model and the reference bi-fidelity model, as well as the low-fidelity model. 
The relative errors are computed and averaged over a test set of $512$ randomly sampled parameter points that is independent of $\Gamma$.
  We conducted tests for different numbers of reduced bases, specifically $N_{rb}=3, 6, 9, 12$, corresponding to a total of $N_{T}=33, 33, 35, 37$ high-fidelity simulations  conducted during the offline stage.
It is evident from the figure that the error converges as the number of reduced bases increases. Notably, the relative $l^2$ errors can reach a level of ${\mathcal O}(10^{-6})$ when employing $12$ reduced bases. This demonstrates the effectiveness of our reduced basis method in accurately approximating the high-fidelity solutions. In contrast, the low-fidelity model exhibits relatively large errors.
It is worth noting that when the size of the reduced basis is sufficiently large, our method outperforms the reference bi-fidelity method in terms of accuracy which can be also observed in a histogram of the relative errors for both bi-fidelity methods plotted in the right figure. This is attributed to leveraging equation information through reduced equation  (\ref{eq:reduced_short})  during the online stage.

\begin{figure}[!htb]
	\centering
		\begin{subfigure}{0.35\textwidth}
  \includegraphics[width=\textwidth]{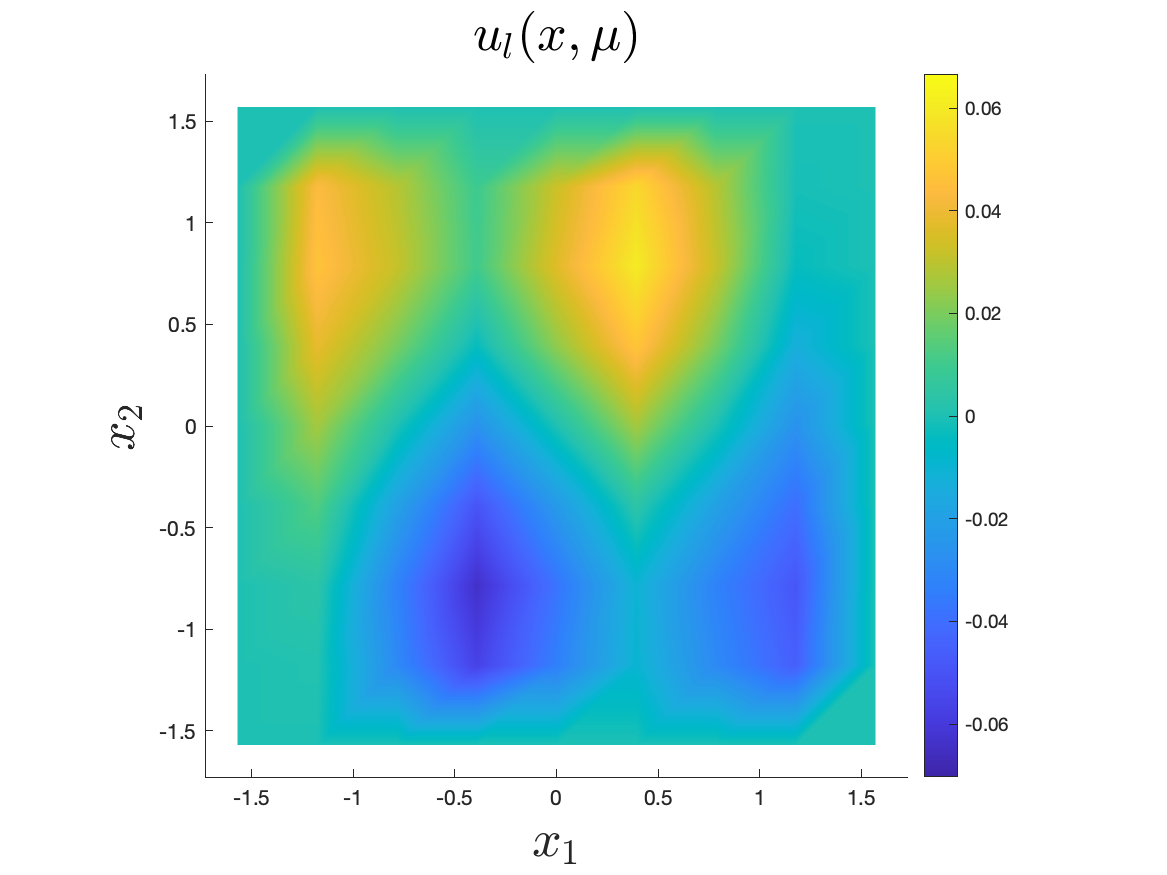}
  \caption{$u_l$}
  \label{fig:Sol_plots_ul_nrb12_2DNL}
\end{subfigure}
\hspace{-.3in}
	\begin{subfigure}{0.35\textwidth}
  \includegraphics[width=\textwidth]{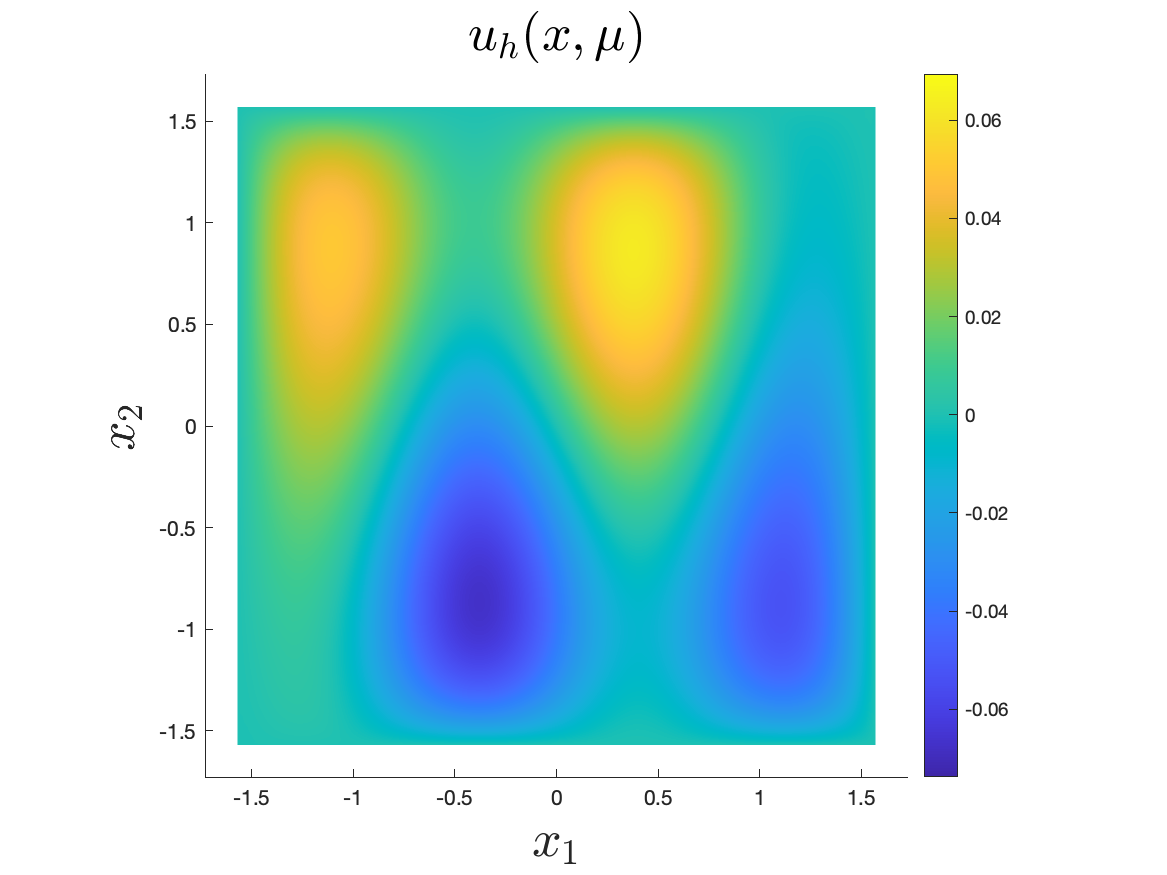}
  \caption{$u_h$}
  \label{fig:Sol_plots_uh_nrb12_2DNL}
\end{subfigure}
 \hspace{-.3in}
 	  \begin{subfigure}{0.35\textwidth}
  \includegraphics[width=\textwidth]{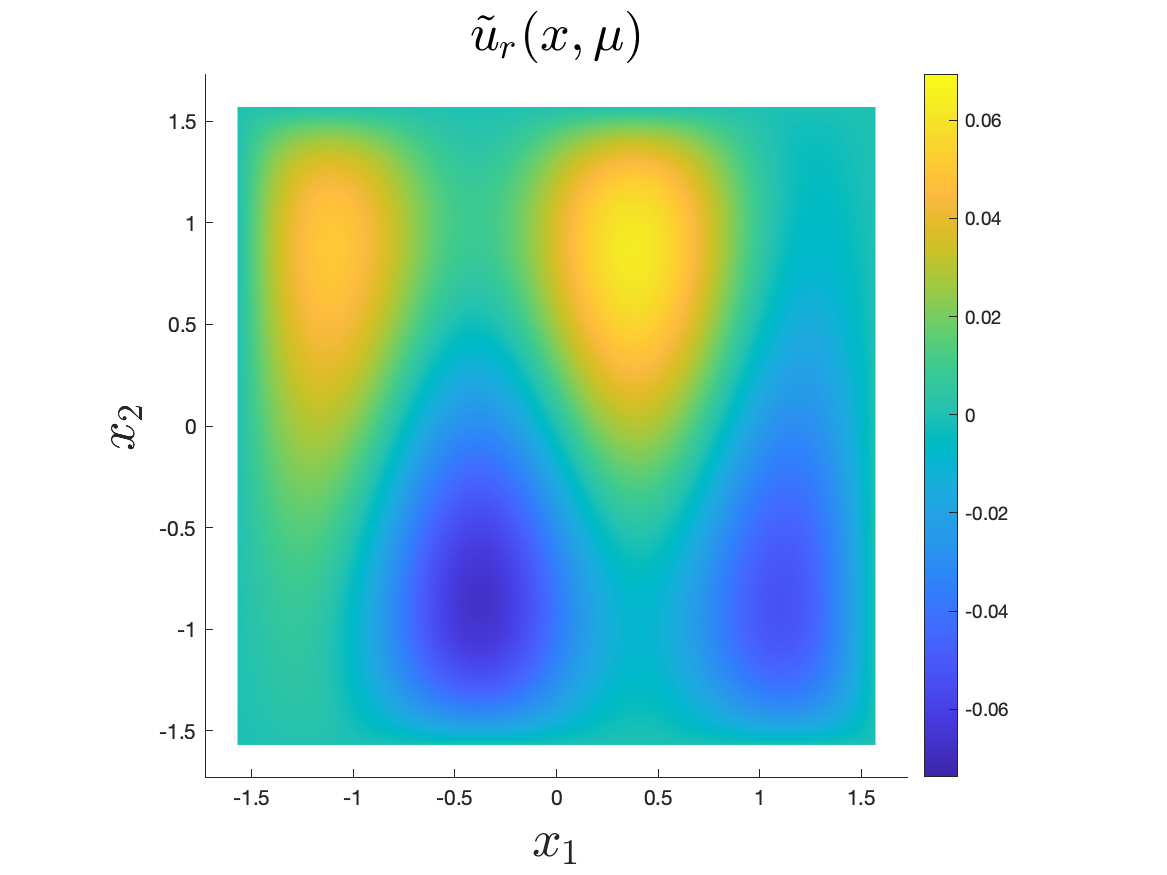}
  \caption{$\tilde{u}_{r}$}
  \label{fig:Sol_plots_ur_nrb12_2DNL}
 \end{subfigure}
  \vspace*{-4mm}
 \caption{Example \ref{sec:2DNL} : Solutions obtained by low-fidelity model $u_l(x,\mu)$, high-fidelity model $u_h(x,\mu)$ and reduced model $\tilde{u}_r(x,\mu)$ at $(\mu_1,\mu_2, \mu_3) \approx (0.083,0.346, 0.010)$ with the number of reduced basis is $N_{rb}=12$,  $n_L=30$ and $n_f=2$.}
\label{figure:solutions_2DNL}
  \vspace*{-3mm}
\end{figure}

Figure \ref{figure:solutions_2DNL} depicts the solution plots for the low- and high-fidelity models and our reduced-order model utilizing $12$ bases, computed at a parameter value of $(\mu_1,\mu_2, \mu_3) \approx (0.083,0.346, 0.010)$. Notably, the solutions of high-fidelity and the proposed models have remarkable similarities across the entire domain. The solution from the low-fidelity model is not  accurate but it captures some major behavior of the high-fidelity solution which is crucial for our proposed bi-fidelity method. 

In addition, Table \ref{tab:2dnonlinear2} presents the computation times for the offline and online stages of our method with varying numbers of reduced bases, along with the corresponding computation savings. Our approach achieved a nearly $\mathcal{O}(200)$ times speed-up, demonstrating its computational efficiency compared to the high-fidelity model.

\begin{table}[!htb]
\centering
\begin{tabular}{|c||c|c|c|c|c|}
\hline
$N_{rb}$&$T_{rb}^{(off)}$ & $T_{rb}^{(on)}$ &$T_{l}$ & $T_{h}$ & $T_{h}/T_{rb}^{(on)}$  \\
\hline
3 &14.0393s& 1.6554e-3s  & 8.6353e-4s &3.2721e-1s& 198\\  
6 &16.3735s & 1.6945e-3s &8.6353e-4s &3.2721e-1s&193 \\    
9 & 16.6554s & 1.7378e-3s &8.6353e-4s &3.2721e-1s& 188\\  
12 & 17.2338s  & 1.7554e-3s &8.6353e-4s &3.2721e-1s&186 \\ 
\hline
\end{tabular}
\caption{Example \ref{sec:2DNL}: Computation time for offline stage $T_{rb}^{(off)}$, online stage $T_{rb}^{(on)}$, low-fidelity run $T_{l}$ and high-fidelity run $T_{h}$ with respect to the number of reduced bases; Computation times are averaged over $512$ parameters.}
\label{tab:2dnonlinear2}
\end{table}

 \subsection{Cubic reaction-diffusion equation}
 \label{sec:2DCubic_Newton2}
Next, we consider the following reaction-diffusion equation considered in \cite{chen2019l1}.
\beq
 \label{eq:cubic_Newton2}
 \bsp
-\mu_2\Delta u + u(u-\mu_1)^2 &= f(x) \ \ \textrm{in} \ \ \Omega,\\
u &= 0  \ \ \textrm{on} \ \ \partial \Omega.
 \end{split}
 \eeq
 The parameters $(\mu_1,\mu_2)$ are chosen from $[0.4,5]\times[0.4,2]$.
 The source term is given by $f(x) = 100 \sin(2 \pi x) \cos( 2 \pi x)$.
For high- and low-fidelity models, we use the same FEM solver using Newton's method on $200\times 200$ and $10\times 10$ meshes, respectively. 
The following shows the $k+1^{th}$ Newton's iteration of the FEM solver. 
\beq
\bsp
-\mu_2\Delta (\delta u^{k+1}) + h'(u^k,\mu_1) (\delta u^{k+1}) &= \mu_2\Delta u^{k} - h(u^k,\mu_1)+f(x),\\
u^{k+1} &= u^k + \delta u^{k+1},
\end{split}
\eeq
where $h(u^k,\mu_1) = u(u-\mu_1)^2$.
The dense point set $\Gamma$ is chosen to be $400$ random points sampled in the parameter space.
In this case, we need $n_L=33$ and $n_f=40$ high-fidelity operators and the right-hand sides to construct the reduced model.

 

\begin{figure}[!htb]
	\centering
 \begin{subfigure}{0.45\textwidth}
\includegraphics[width=\textwidth]{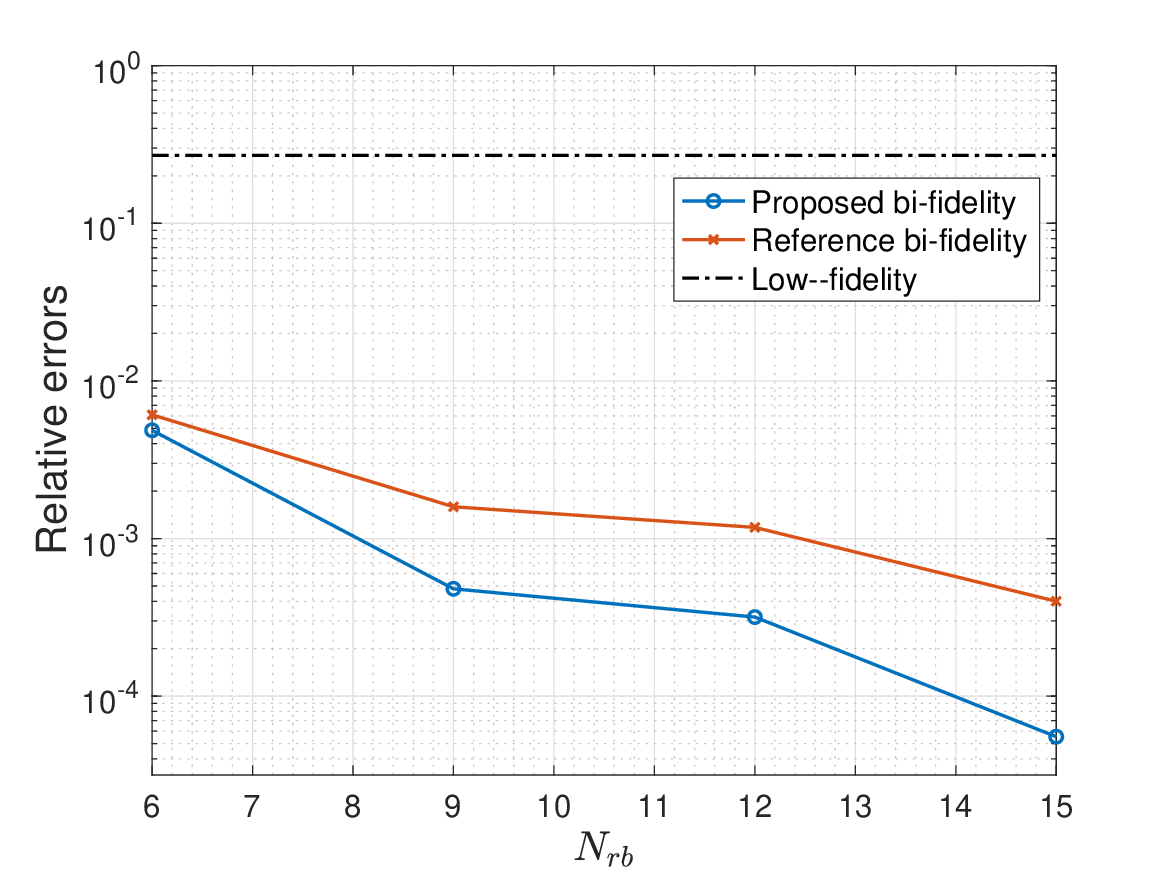}
\caption{}
    \label{fig:Errors_2DCubic_Newton2}
    \end{subfigure}
    \hspace{.1in}
 \begin{subfigure}{0.45\textwidth}
  \includegraphics[width=\textwidth]{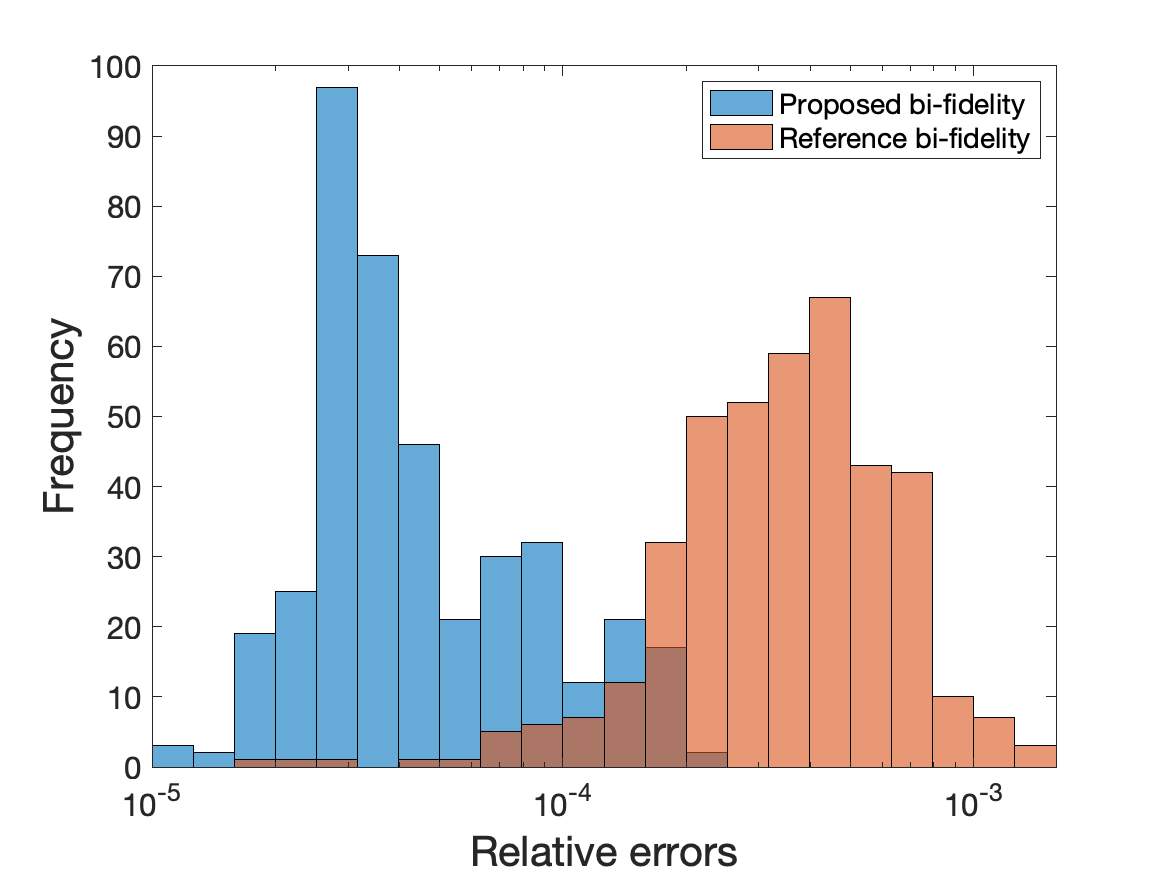}
  \caption{}
  \label{fig:Histo_2DCubic_Newton2}
\end{subfigure}
\vspace{-.2in}
  \caption{Example \ref{sec:2DCubic_Newton2}: a. The relative errors of the reduced solution by our new bi-fidelity method and the reference bi-fidelity method with different number of reduced basis; b. Histogram of the relative errors for both new and the reference bi-fidelity methods when $N_{rb} = 15$.}
  \label{figure:Errors_2D_cubic_Newton2}
 \end{figure}

Figure \ref{figure:Errors_2D_cubic_Newton2}  represents the relative errors of our bi-fidelity model and the reference bi-fidelity model, as well as the low-fidelity model. 
The relative errors are averaged over $400$ randomly sampled test parameter points not included in $\Gamma$.
We conducted tests for different numbers of reduced bases, specifically $N_{rb} = 6, 9, 12, 15$, corresponding to a total of $N_{T}=62, 63, 64, 65$ high-fidelity simulations performed during the offline stage. We also plot the error distributions of both the reference and our bi-fidelity methods over the $400$ test parameters. As the number of reduced bases increases, the averaged relative errors tend to decrease. Both bi-fidelity methods demonstrate significant improvements over the associated low-fidelity model in terms of  accuracy. 
Additionally, 
our proposed bi-fidelity algorithm can achieve a better accuracy compared to the reference bi-fidelity method, thanks to the  reduced order model (\ref{eq:reduced_short}) we enforced during the online stage.
  
\begin{figure}[!htb]
	\centering
 \begin{subfigure}{0.35\textwidth}
  \includegraphics[width=\textwidth]{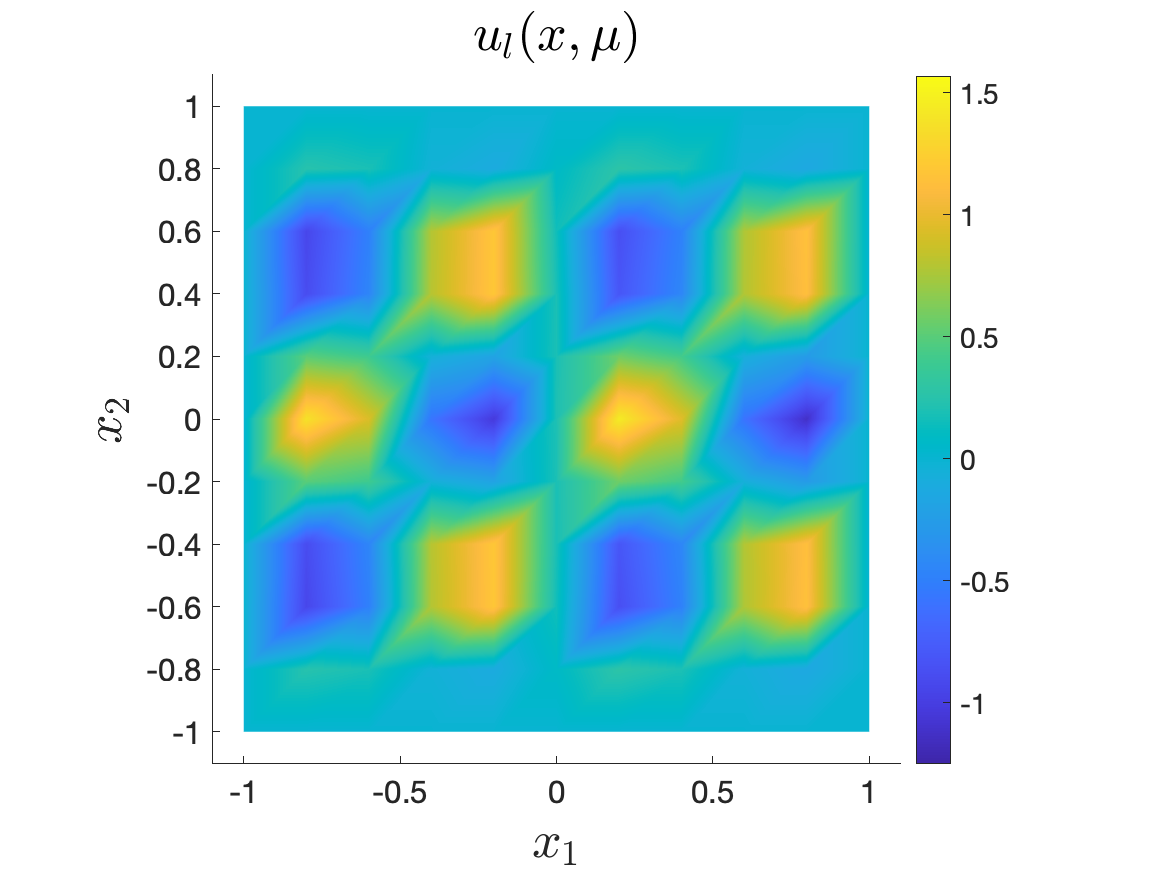}
  \caption{$u_l$}
  \label{fig:Sol_plots_ul_nrb6_2DCubic_Newton2}
 \end{subfigure}
 \hspace{-.3in}
	\begin{subfigure}{0.35\textwidth}
  \includegraphics[width=\textwidth]{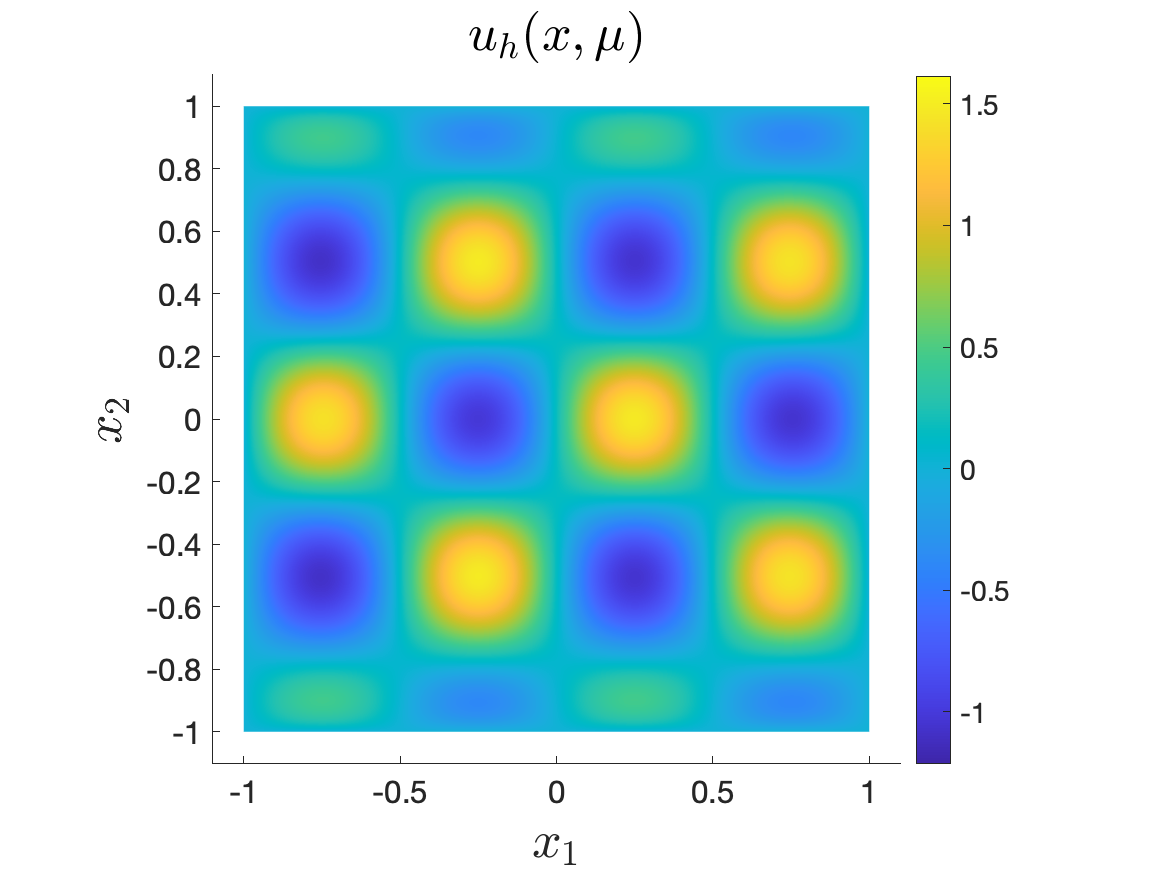}
  \caption{$u_h$}
  \label{fig:Sol_plots_uh_nrb6_2DCubic_Newton2}
\end{subfigure}
\hspace{-.3in}
	\begin{subfigure}{0.35\textwidth}
  \includegraphics[width=\textwidth]{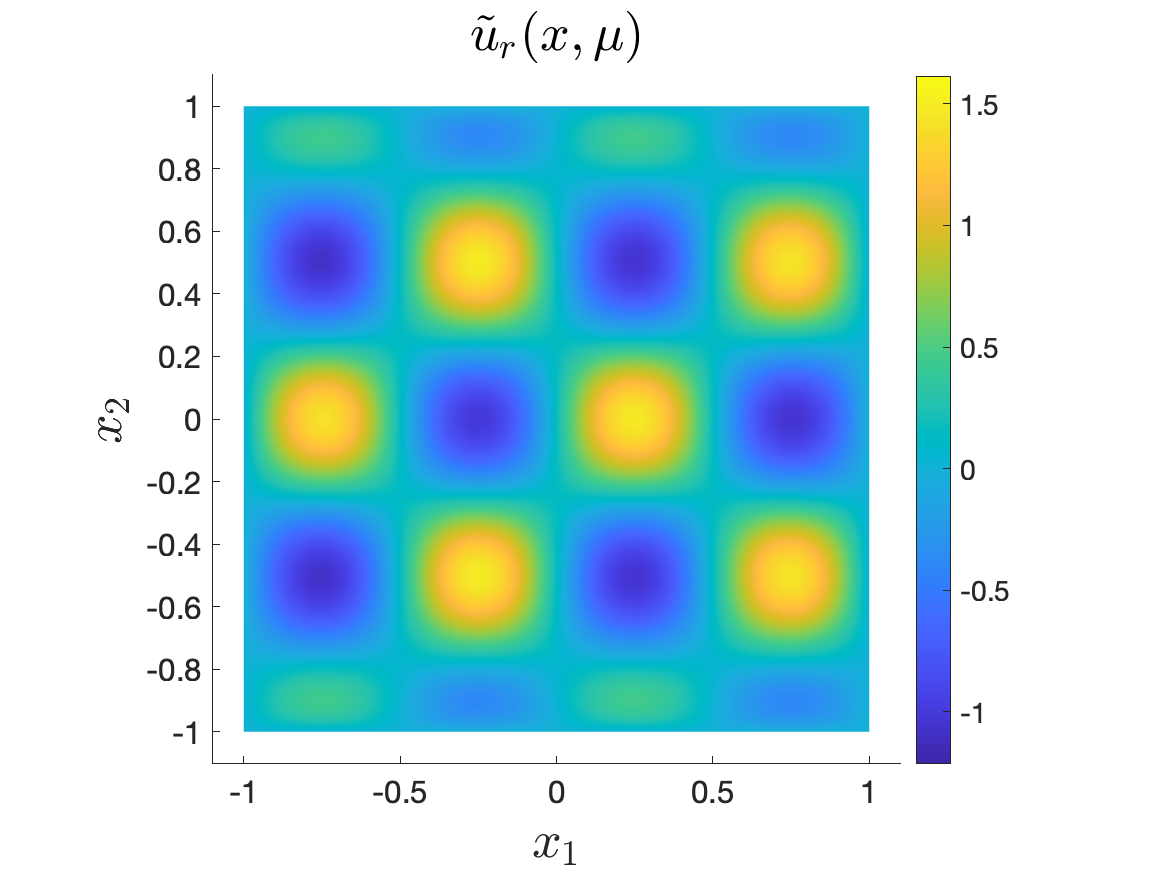}
  \caption{$\tilde{u}_{r}$}
  \label{fig:Sol_plots_ur_nrb6_2DCubic_Newton2}
 \end{subfigure}
  \vspace*{-4mm}
 \caption{Example \ref{sec:2DCubic_Newton2}: Solutions obtained by coarse model $u_l(x,\mu)$, full model $u_h(x,\mu)$ and reduced model $\tilde{u}_r(x,\mu)$ at $(\mu_1, \mu_2) \approx (4.202, 0.822)$ with the number of reduced basis  $N_{rb}=15$, $n_L=33$, $n_f=40$.} 
\label{figure:solutions_cubic_Newton2}
  \vspace*{-3mm}
\end{figure}

Figure \ref{figure:solutions_cubic_Newton2} showcases the solutions of the low-fidelity model, high-fidelity model, and our proposed reduced model at a parameter point $(\mu_1, \mu_2) \approx (4.202, 0.822)$. It is evident that our algorithm, utilizing $N_{rb} = 15$ reduced bases, provides a highly accurate approximation of the solution. While the low-fidelity solution depicted in the figure may not be as precise, it still manages to capture some of the key characteristics present in the high-fidelity solution, which is informative for our bi-fidelity model.

Table \ref{tab:cubic2_Newton2} lists the computation time for the offline/online  stage of our bi-fidelity model, low-fidelity and high-fidelity model. It is observed that the low-fidelity run is dominant in the online stage in terms of computation time. Furthermore, the table highlights that the online computation is more than $200$ times faster than the high-fidelity run. This substantial reduction in computation time demonstrates the significant efficiency gained by employing our bi-fidelity model.


 \begin{table}[ht!]
\centering
\begin{tabular}{|c||c|c|c|c|c|}
\hline
$N_{rb}$&$T_{rb}^{(off)}$ & $T_{rb}^{(on)}$ &$T_{l}$ & $T_{h}$ & $T_{h}/T_{rb}^{(on)}$   \\
\hline
6 &124.1668s & 5.4035e-3s  &4.0595e-3s &1.4967s& 277  \\ 
9 &132.1089s   &5.8559e-3s  &4.0595e-3s &1.4967s& 256\\ 
12 &134.8519s & 6.0180e-3s & 4.0595e-3s&1.4967s&249 \\
15 & 141.2292s &6.1942e-3s & 4.0595e-3s &1.4967s& 242 \\
\hline
\end{tabular}
\caption{Example \ref{sec:2DCubic_Newton2}: Computation time for offline stage $T_{rb}^{(off)}$, online stage $T_{rb}^{(on)}$, low-fidelity run $T_{l}$ and high-fidelity run $T_{h}$ with respect to the number of reduced bases; Computation times  are averaged over $400$ parameters.}
\label{tab:cubic2_Newton2}
\end{table}

\subsection{2D elliptic problem with high-contrast coefficient}
\label{sec:2DHC}

\begin{figure}[!htb]
	\centering
	\begin{subfigure}{0.33\textwidth}
  \includegraphics[width=\textwidth]{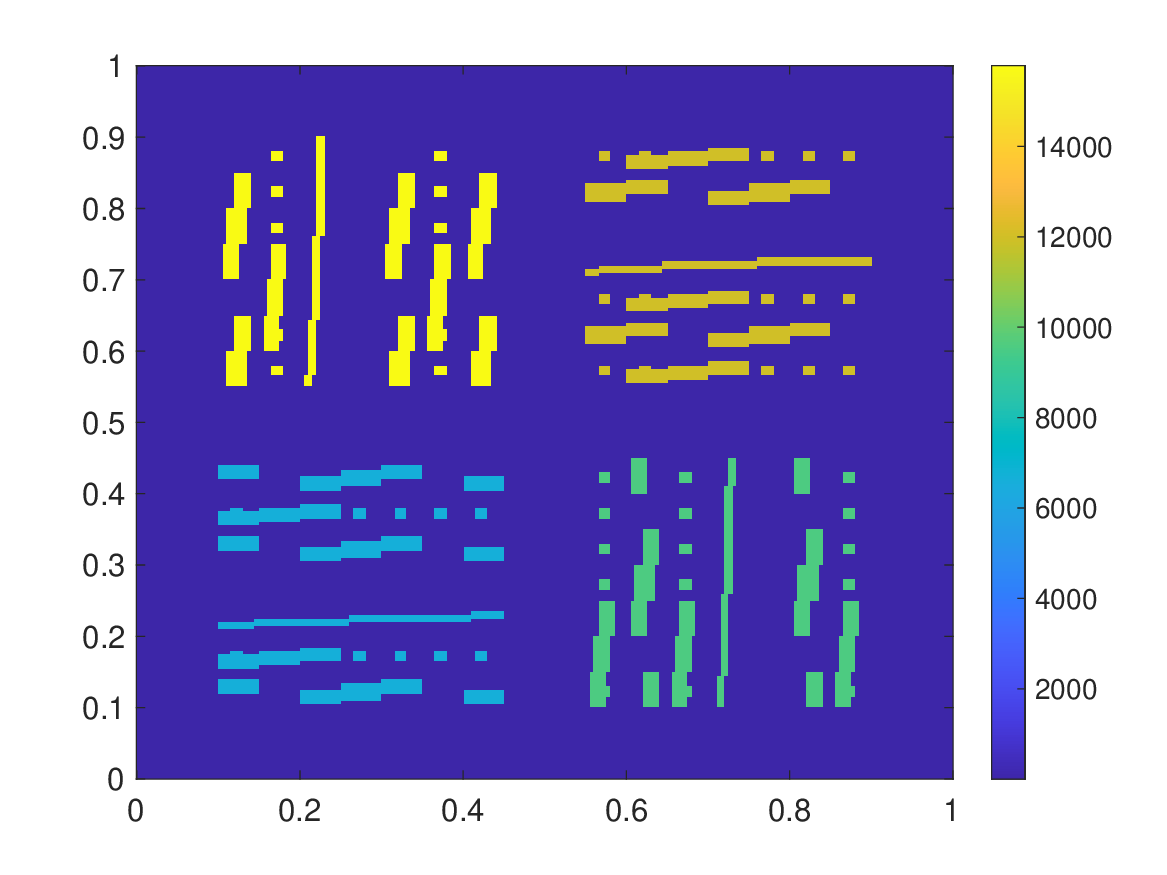}
  \caption{$\mu=(1,1,\frac{1}{2})$}
  \label{fig:kappa1_2DHCMS}
\end{subfigure}
\hspace{-.1in}
  \begin{subfigure}{0.33\textwidth}
  \includegraphics[width=\textwidth]{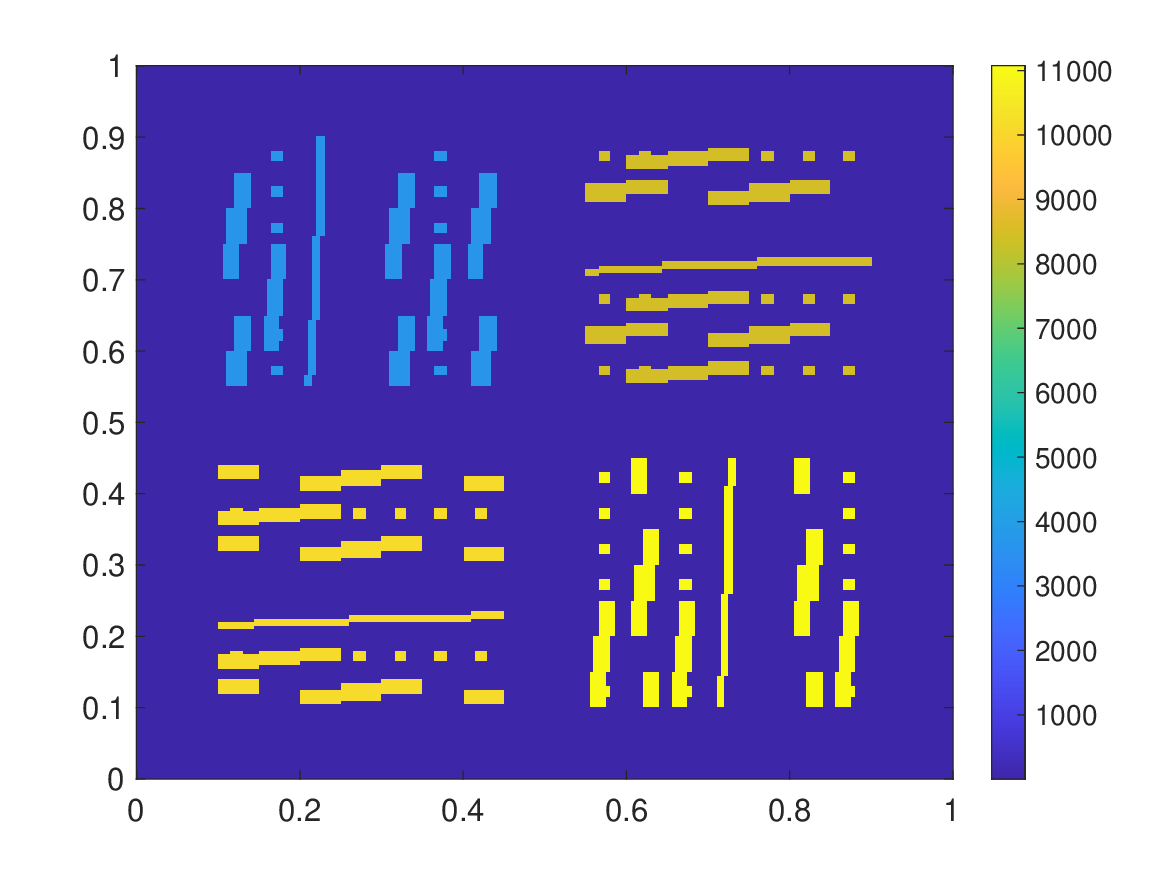}
  \caption{$\mu=(\frac{1}{2},-\frac{2}{3},-1)$}
  \label{fig:kappa2_2DHCMS}
 \end{subfigure}
 \hspace{-.1in}
  \begin{subfigure}{0.33\textwidth}
  \includegraphics[width=\textwidth]{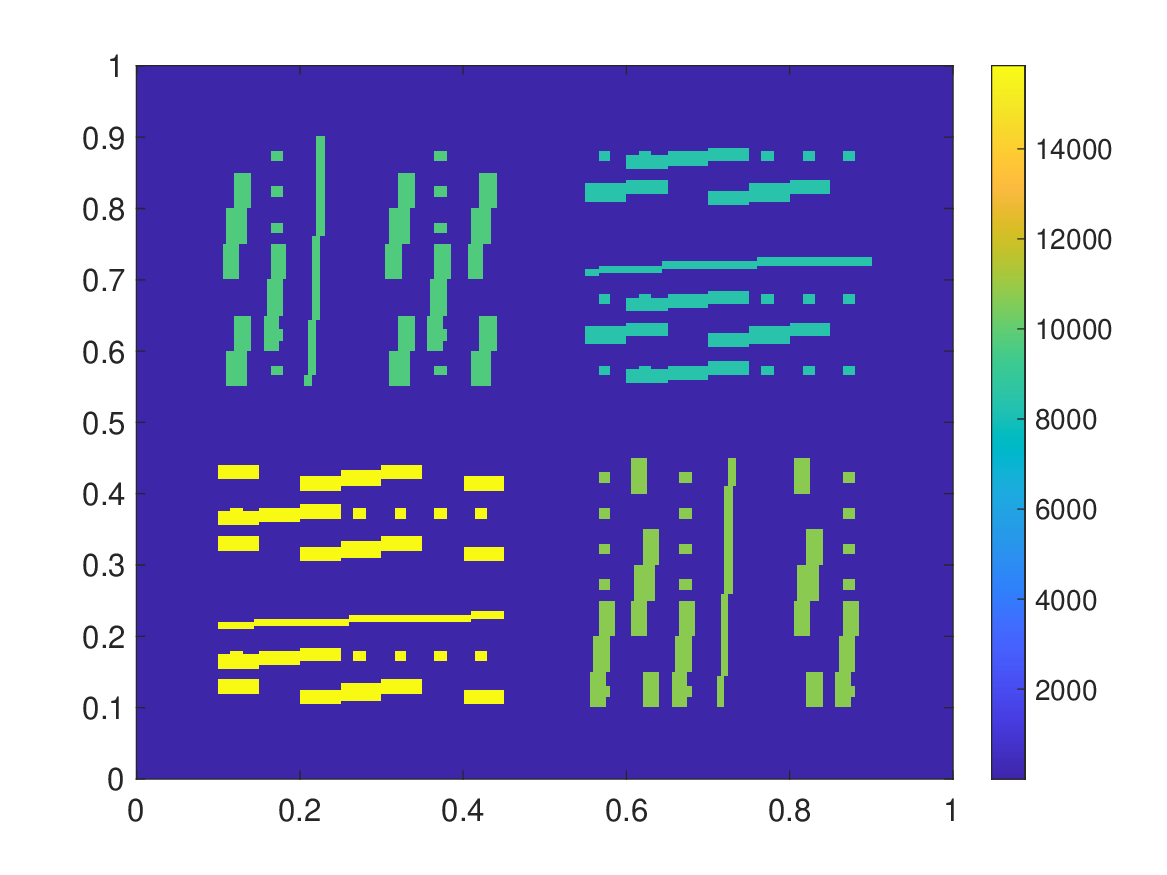}
  \caption{$\mu=(-\frac{1}{2},\frac{2}{3},-\frac{1}{3})$}
  \label{fig:kappa3_2DHCMS}
 \end{subfigure}
  \vspace*{-4mm}
 \caption{Example \ref{sec:2DHC}: Permeability coefficient $\kappa(x,\mu)$ with three different parameters. The value in the background (dark blue region) is $1$ for all cases.}
\label{figure:highcontrast_coeffs}
  \vspace*{-3mm}
\end{figure}

Motivated by an example in \cite{ou2020low}, we consider the following elliptic equation with high-contrast permeability coefficient:
\beq
\label{eq:2dhighcontrast}
\begin{split}
-\div(\kappa(x,\mu) \nabla u(x,\mu)) = f(x), \ \  &\textrm{in} \ \ \Omega,\\
u(x,\mu) = 0 \ \  &\textrm{on} \ \ \partial\Omega,
\end{split}
\end{equation}
where the permeability coefficient $\kappa(x,\mu)$ is given by
\beq
\label{eq:kappa_high}
\kappa(x,\mu) = \sum\limits_{i=1}^5  \alpha_i(\mu)\kappa_i(x).
\eeq
The value of $\kappa_i(x)$ for $i = 1,\dots,4$ is $10^4$ inside the channels and $0$ in the matrix (background), where the $\kappa_i(x)$ are the spatially distributed part of the coefficient. The value of $\kappa_5(x)$ is $1$ in the matrix and $0$ in channels.
The function $\alpha_i(\mu)$ is defined as follows:
\beq
\label{eq:alpha_mu}
\bsp
\alpha_1(\mu) &= \frac{0.8+1.6 \mu_1^4}{1+\mu_1^4},\\
\alpha_2(\mu) &= 1.1+0.8\sin(\mu_1+\mu_2+\mu_3),\\
\alpha_3(\mu) &= 1.1+0.7\cos(\mu_1^2+\mu_2^2+\mu_3^2),\\
\alpha_4(\mu) &= 1.2 - \frac{0.3\mu_3^2}{1+\mu_2^2 \mu_3^2},\\
\alpha_5(\mu) &= 1,
\end{split}
\eeq
where the parameters $\mu = (\mu_1,\mu_2, \mu_3)$ are taken from the parameter domain $[-1,1]^3$. 
We note that the contrast ratio (the ratio between the maximum and minimum value) of $\kappa(x,\mu)$ is  ${\mathcal O}(10^4)$.
Figure \ref{figure:highcontrast_coeffs} illustrates three different realizations of  the coefficients $\kappa(x,\mu)$. 
For the source function, we take $f(x) = \sin(\pi x_1) \sin(\pi x_2)$. 
Even though this problem is linear,  the computation cost can be high to achieve reasonable accuracy for each parameter due to the high contrast nature of the problem. The existing reduced order model techniques for this class of problems are intrusive \cite{efendiev2013generalized,efendiev2012systematic}.
For high- and low-fidelity models,  we use the same FEM solver on  $128\times 128$ and $4\times 4$ FEM grids respectively. 
The dense point set $\Gamma$ consists of $512$ random sample points in the parameter domain.
We use $n_L=5$ high-fidelity operators to construct the reduced operators and the right-hand side is independent of parameters.

\begin{figure}[!htb]
	\centering
 \begin{subfigure}{0.45\textwidth}
\includegraphics[width=\textwidth]{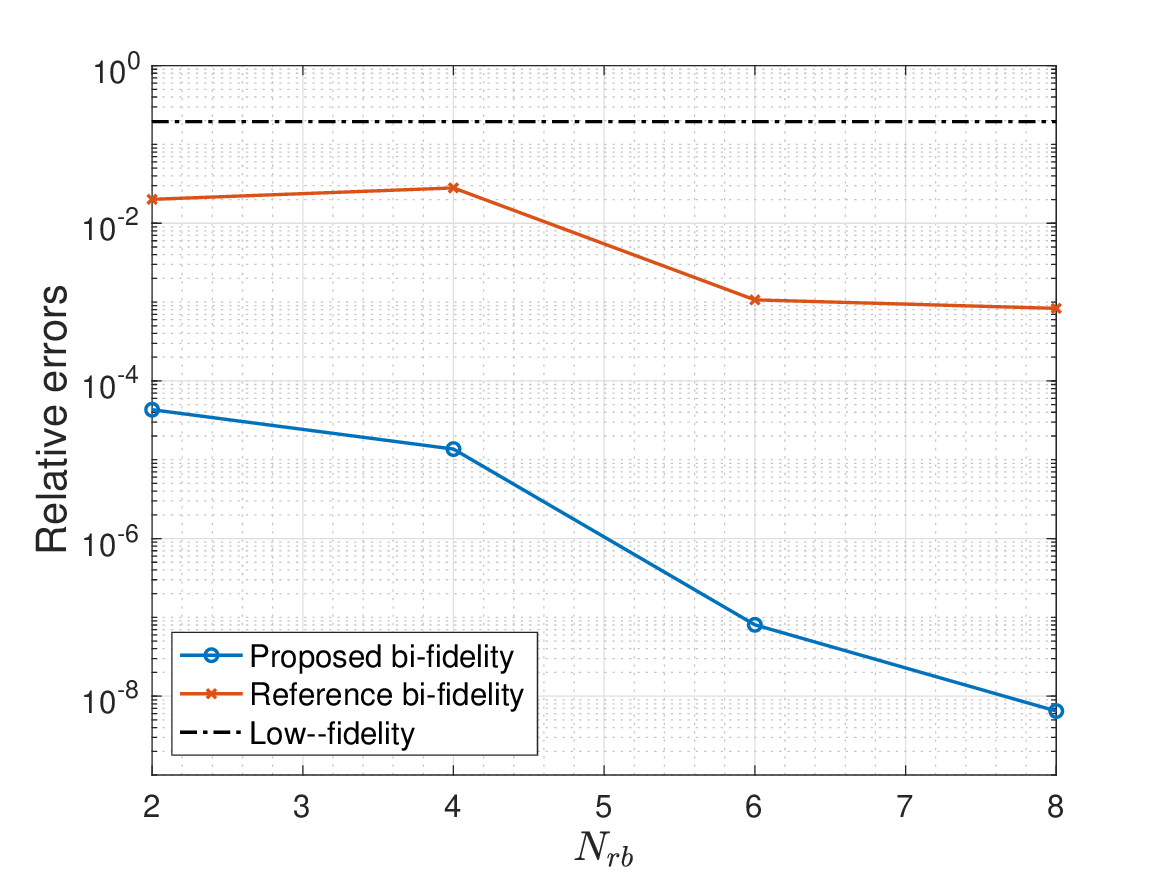}
\caption{}
    \label{fig:Errors_2DHC}
    \end{subfigure}
    \hspace{.1in}
 \begin{subfigure}{0.45\textwidth}
  \includegraphics[width=\textwidth]{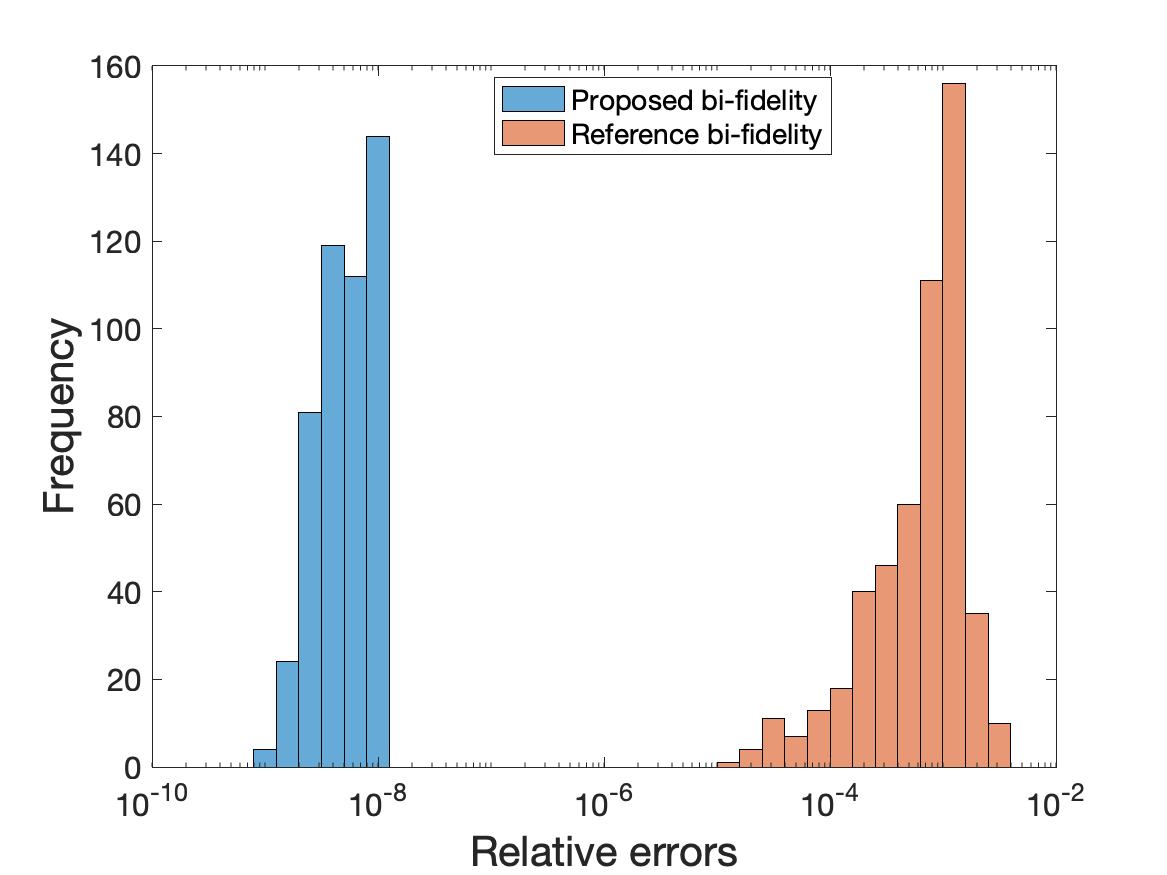}
  \caption{}
  \label{fig:Histo_2DHC}
\end{subfigure}
\vspace{-.2in}
  \caption{Example \ref{sec:2DHC}: a. The relative errors of the reduced solution by our new bi-fidelity method and the reference bi-fidelity method with different number of reduced basis; b. Histogram of the relative errors for both new and the reference bi-fidelity methods when $N_{rb} = 8$.}
  \label{figure:2dhighcontrast}
 \end{figure}

Figure \ref{figure:2dhighcontrast} plots the averaged relative errors of the reduced solutions obtained by the proposed, and reference bi-fidelity methods.
The relative errors computed and averaged over the independent test set of $512$ randomly sampled parameters.
We test the algorithm with $N_{rb} = 2, 4, 6, 8$ reduced bases which entail $N_{T} = 5, 5, 6, 8$ respective high-fidelity simulations in the offline stage. 
Both bi-fidelity algorithms can produce significantly improved accuracy compared to the associated low-fidelity model. 
The results clearly demonstrate the fast convergence and superior accuracy of our proposed bi-fidelity method compared to the reference method. The relative error of level ${\mathcal O}(10^{-9})$ is achieved when we utilize $N_{rb} = 8$ reduced bases. This can be further confirmed by the error distributions of both bi-fidelity frameworks, suggesting the effectiveness of leveraging the reduced  equation (\ref{eq:reduced_short})  in the online stage. 

\begin{figure}[!htb]
	\centering
		\begin{subfigure}{0.35\textwidth}
  \includegraphics[width=\textwidth]{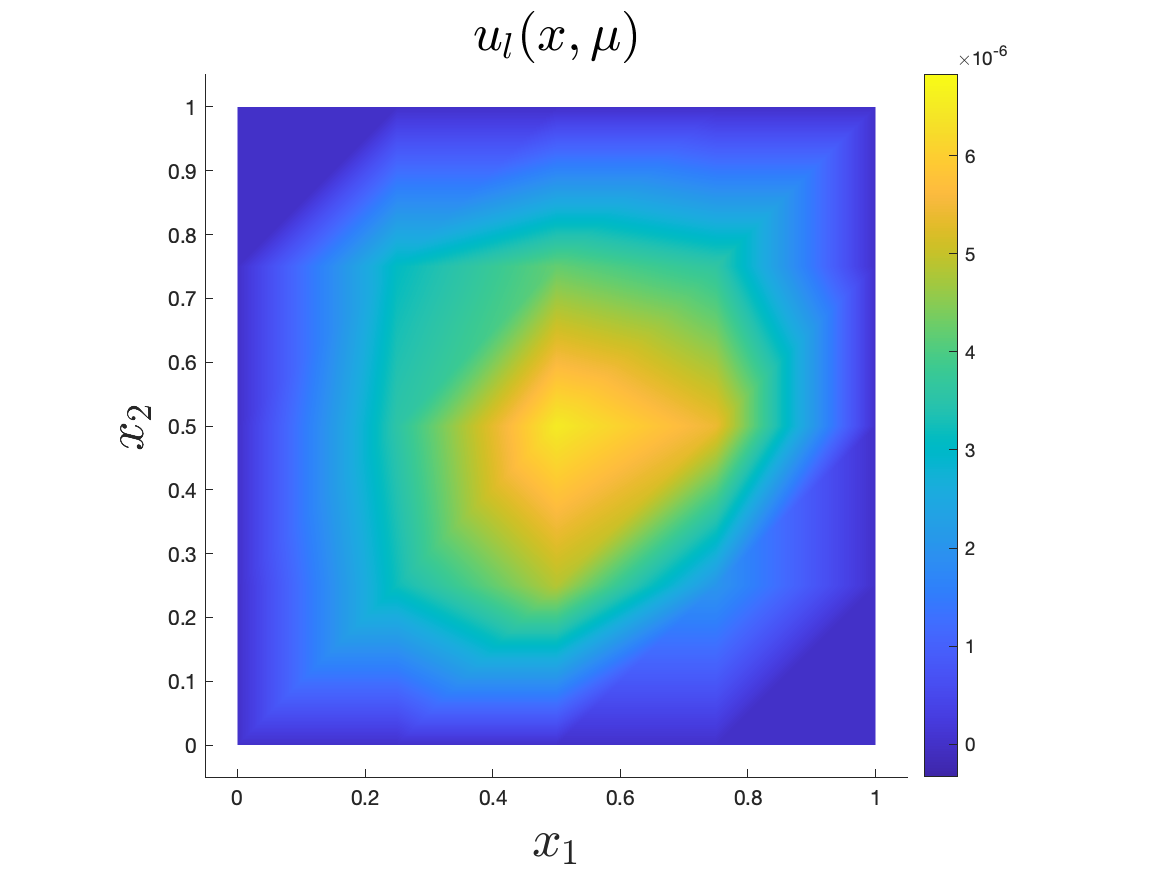}
  \caption{$u_l$}
  \label{fig:Sol_plots_ul_nrb8_2DHC}
\end{subfigure}
 \hspace{-.3in}
	\begin{subfigure}{0.35\textwidth}
  \includegraphics[width=\textwidth]{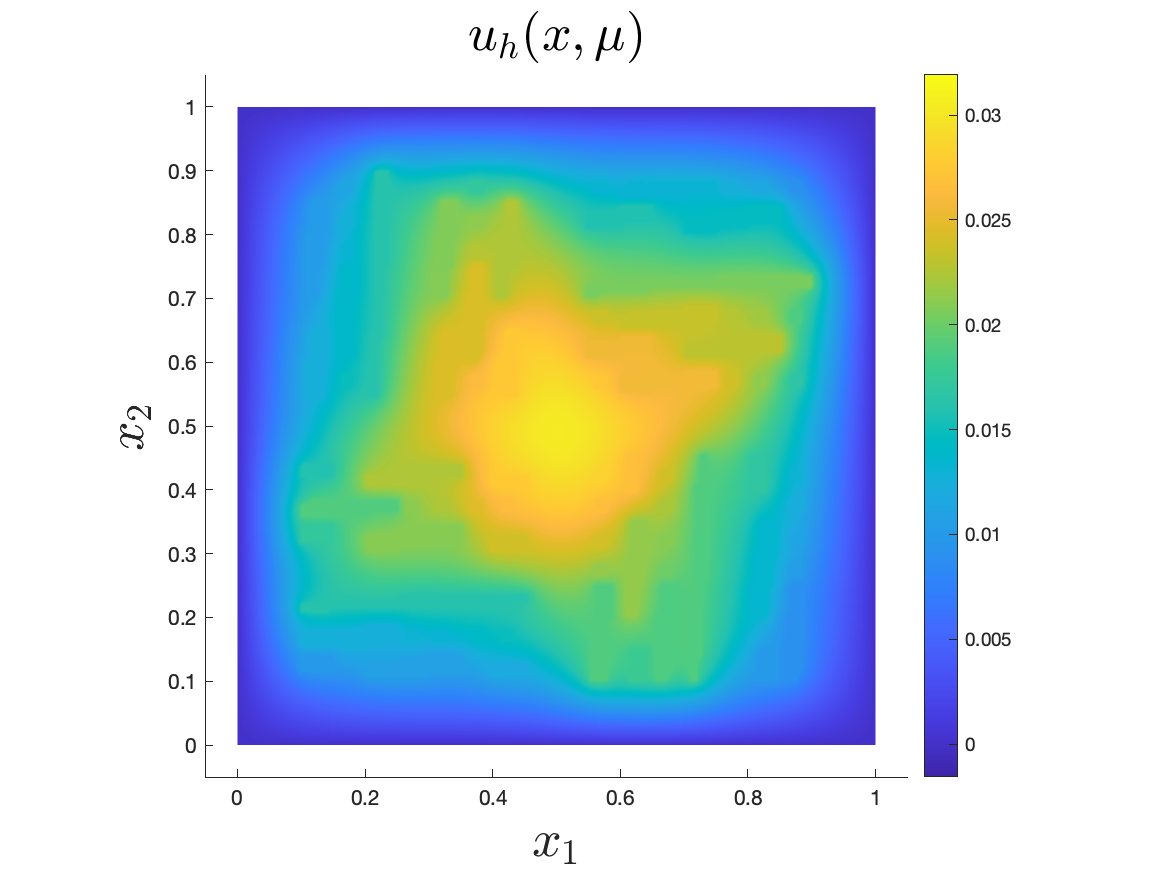}
  \caption{$u_h$}
  \label{fig:Sol_plots_uh_nrb8_2DHC}
\end{subfigure}
\hspace{-.3in}
	  \begin{subfigure}{0.35\textwidth}
  \includegraphics[width=\textwidth]{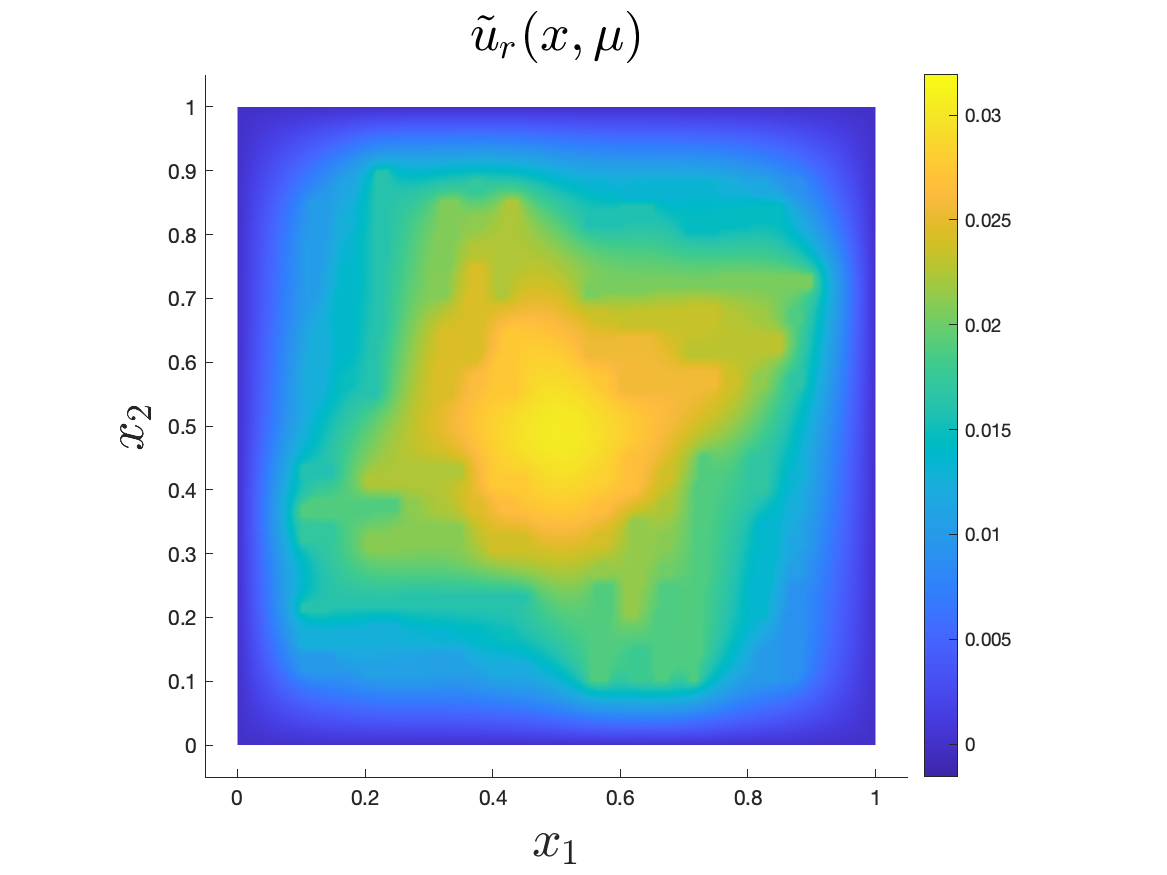}
  \caption{$\tilde{u}_{r}$}
  \label{fig:Sol_plots_ur_nrb8_2DHC}
 \end{subfigure}
  \vspace*{-4mm}
 \caption{Example \ref{sec:2DHC}: Solutions obtained by coarse model $u_l(x,\mu)$, full model $u_h(x,\mu)$ and reduced model $\tilde{u}_r(x,\mu)$ at $(\mu_1,\mu_2, \mu_3) \approx (0.712,-0.402,0,878)$ with the number of reduced basis  $N_{rb}=8$, $n_L=5, n_f=1$.}
\label{figure:solutions_2DHigh}
  \vspace*{-3mm}
\end{figure}

Figure \ref{figure:solutions_2DHigh} depicts the solutions obtained by low-fidelity, high-fidelity and our reduced models at parameter $(\mu_1,\mu_2, \mu_3) \approx (0.712,-0.402,0,878)$. Compare to the low-fidelity solution, our bi-fidelity 
solution shows better agreement with the high-fidelity solution.
The low-fidelity solution may not be highly accurate, but it does capture the main characteristics and fundamental parametric dependence of the high-fidelity solution. This rationale strengthens the validity of the proposed bi-fidelity method.

The computation time for the offline, online stage of our method and the associated low-, high-fidelity models are presented in Table \ref{tab:2dhighcontrast}. 
Given that the dominant cost in the online stage is attributed to the low-fidelity simulation, which is considerably less expensive than the high-fidelity model, our proposed bi-fidelity model enables a remarkable speedup of approximately 90 times compared to the corresponding high-fidelity runs. 

\begin{table}[!htb]
\centering
\begin{tabular}{|c||c|c|c|c|c|}
\hline
$N_{rb}$&$T_{rb}^{(off)}$ & $T_{rb}^{(on)}$ &$T_{l}$ & $T_{h}$ & $T_{h}/T_{rb}^{(on)}$   \\
\hline
2 & 13.7733s & 3.8753e-4s & 2.0885e-4s &3.7323e-2s& 96\\ 
4 & 13.7745s & 3.9741e-4s &2.0885e-4s& 3.7323e-2s& 94\\  
6 & 13.7765s & 4.0239e-4s &2.0885e-4s& 3.7323e-2s& 93\\  
8 & 13.7771s & 4.1041e-4s &2.0885e-4s& 3.7323e-2s& 91\\  
\hline
\end{tabular}
\caption{Example \ref{sec:2DHC}:Computation time for offline stage $T_{rb}^{(off)}$, online stage $T_{rb}^{(on)}$, low-fidelity run $T_{l}$ and high-fidelity run $T_{h}$ with respect to the number of reduced bases; Computation times $T_{rb}^{(on)}$, $T_{l}$ and $T_{h}$ are averaged over $512$ parameters.}
\label{tab:2dhighcontrast}
\end{table}


\subsection{2d nonlinear scalar multiscale equation}
\label{sec:2DNLHC}
We next consider the following nonlinear elliptic multiscale equation in $\Omega = [0,1]^2$ motivated by an example in \cite{mai2023constraint}.
\beq
\label{eq:2dnlhigh}
- \div \left(\kappa(x,\mu) e^{p(x,\mu)} \nabla p(x,\mu) \right) = 2+\sin(x_1)\cos(x_2)  \, \ \textrm{in} \, \ \Omega ,
\eeq
where the coefficient $\kappa(x,\mu)$ is given by (\ref{eq:kappa_high}) and (\ref{eq:alpha_mu}).
the parameters $\mu = (\mu_1,\mu_2, \mu_3)$ are taken from the parameter domain $[0,1]^3$.
This problem is challenging  due to the high-contrast feature and nonlinearity of the equation. Existing reduced basis techniques for efficient simulation of nonlinear multiscale high-contrast problems are often intrusive 

For our bi-fidelity framework, we use the same FEM solver using Picard's iteration with $128\times 128$ and $4\times 4$ grids for high-fidelity and low-fidelity models associated with our algorithm.
We randomly sample $512$ points in the parameter domain to construct the dense point set $\Gamma$.
We use $n_L=12$ selected high-fidelity operators to assemble the reduced operators and the right-hand side does not depend on parameters, i.e., $n_f=1$.

\begin{figure}[!htb]
	\centering
 \begin{subfigure}{0.45\textwidth}
\includegraphics[width=\textwidth]{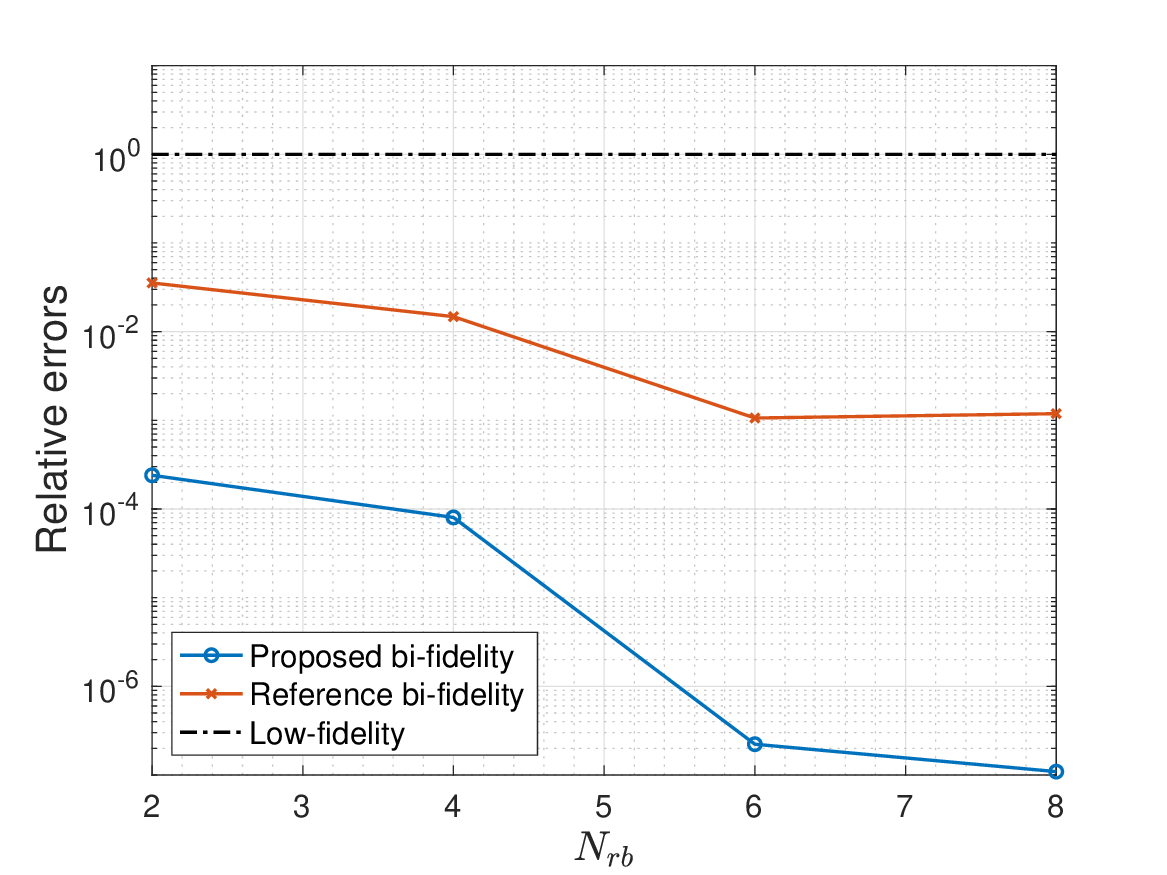}
\caption{}
    \label{fig:Errors_2DNLHC}
    \end{subfigure}
    \hspace{.1in}
 \begin{subfigure}{0.45\textwidth}
  \includegraphics[width=\textwidth]{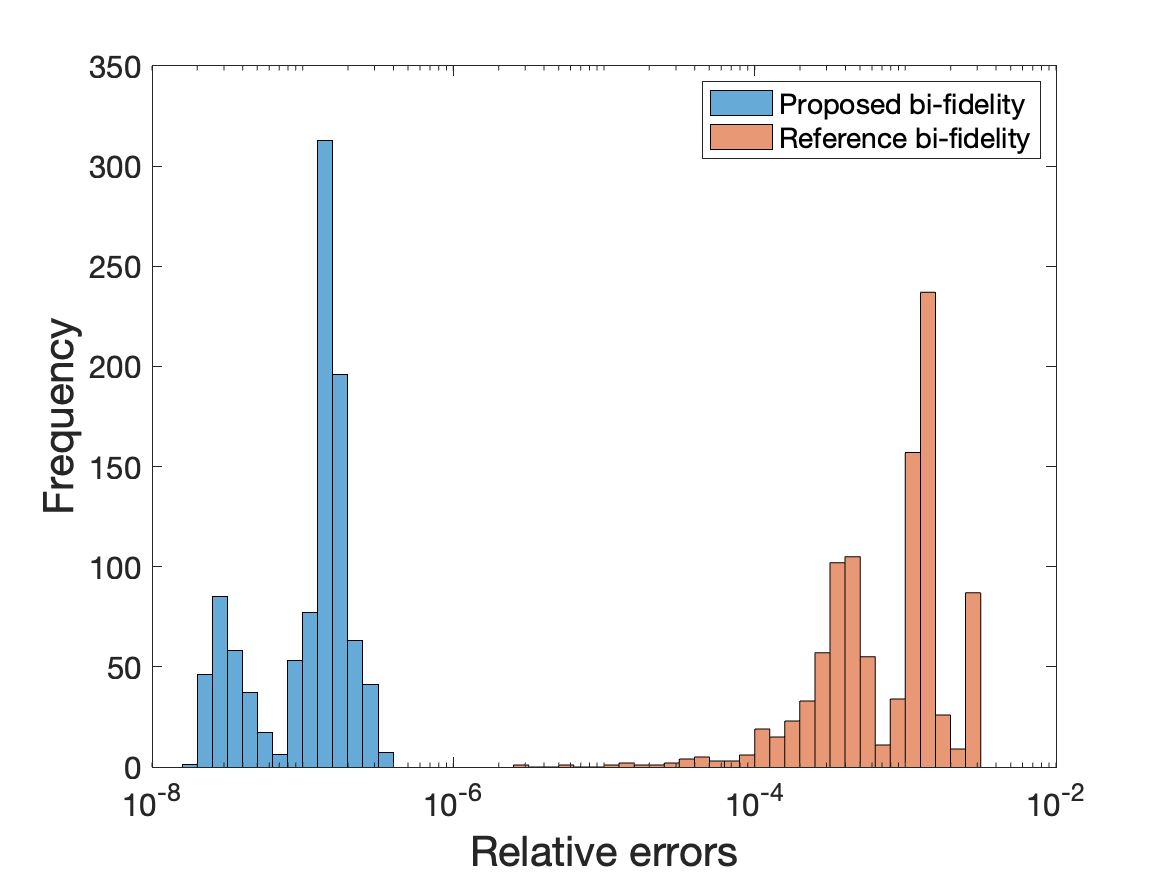}
  \caption{}
  \label{fig:Histo_2DNLHC}
\end{subfigure}
\vspace{-.2in}
  \caption{Example \ref{sec:2DNLHC} a. The relative errors of the reduced solution by our new bi-fidelity method and the reference bi-fidelity method with different number of reduced basis; b. Histogram of the relative errors for both new and the reference bi-fidelity methods when $N_{rb} = 6$.}
  \label{figure:Errors_2D_NL_HC}
 \end{figure}

Figure \ref{figure:Errors_2D_NL_HC} represents the relative errors of  both proposed and reference bi-fidelity approaches with the  number of reduced bases.
The relative errors are computed and averaged over a test set of $512$ randomly selected parameters independent of $\Gamma$
We tested with $N_{rb} = 2, 4, 6, 8$ reduced bases and the corresponding numbers of required high-fidelity simulations in the offline stage are $13, 14, 16, 18$. 
The fast convergence of the relative errors of the proposed method is clearly observed and the error level of ${\mathcal O}(10^{-7})$ is achieved with only $6$ reduced bases.  
We also present the error distributions of both the proposed and reference bi-fidelity methods for $N_{rb} = 6$ over a test set of size $512$.
The superiority of the proposed bi-fidelity method over the reference algorithm is apparent in terms of higher accuracy. 

\begin{figure}[!htb]
	\centering
		\begin{subfigure}{0.35\textwidth}
  \includegraphics[width=\textwidth]{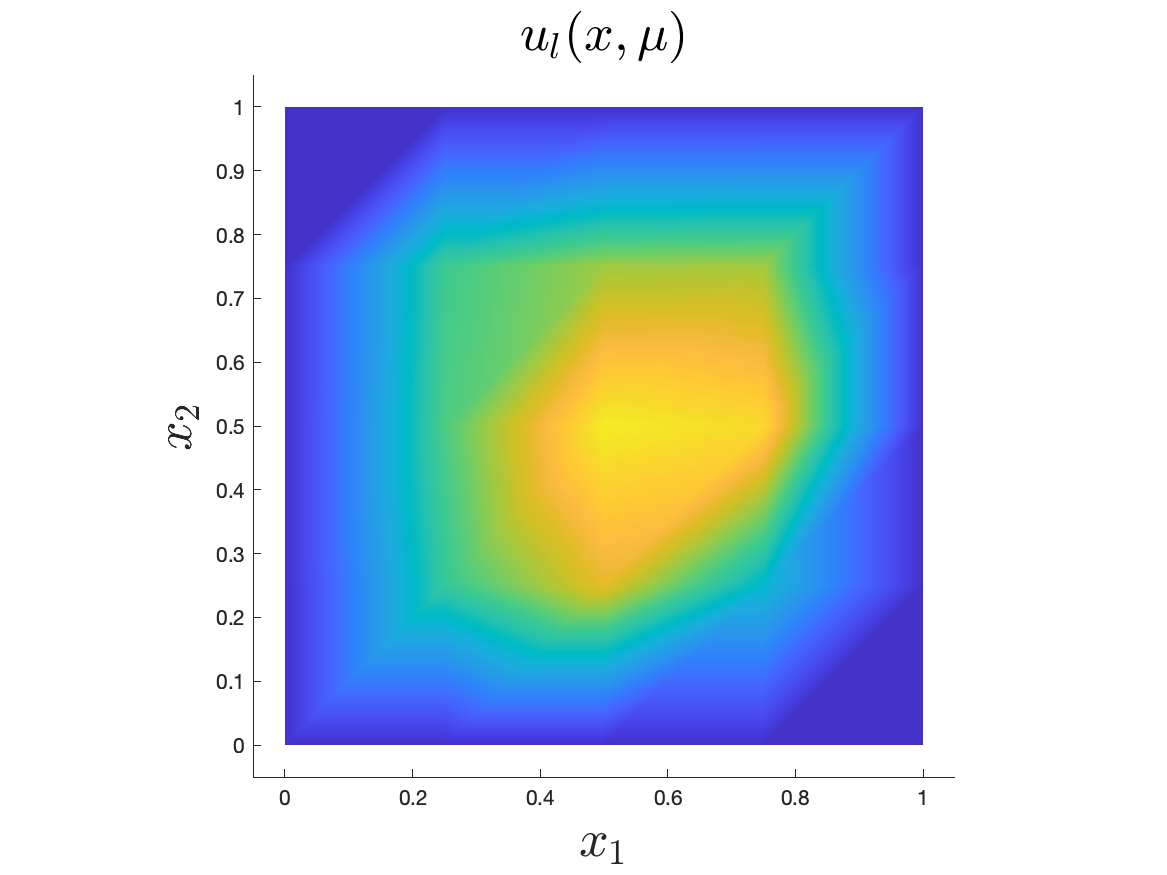}
  \caption{$u_l$}
  \label{fig:Sol_plots_ul_nrb6_2DNLHC}
\end{subfigure}
 \hspace{-.3in}
	\begin{subfigure}{0.35\textwidth}
  \includegraphics[width=\textwidth]{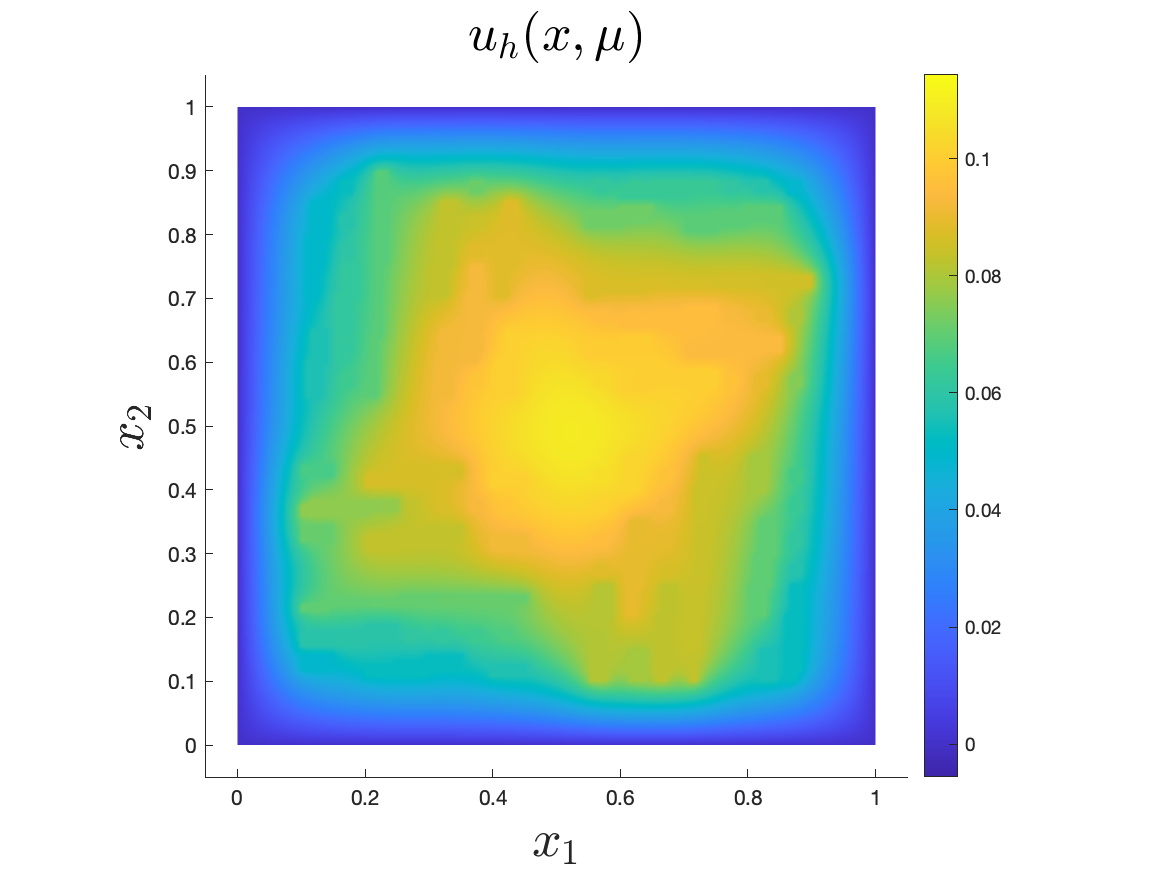}
  \caption{$u_h$}
  \label{fig:Sol_plots_uh_nrb6_2DNLHC}
\end{subfigure}
\hspace{-.3in}
	  \begin{subfigure}{0.35\textwidth}
  \includegraphics[width=\textwidth]{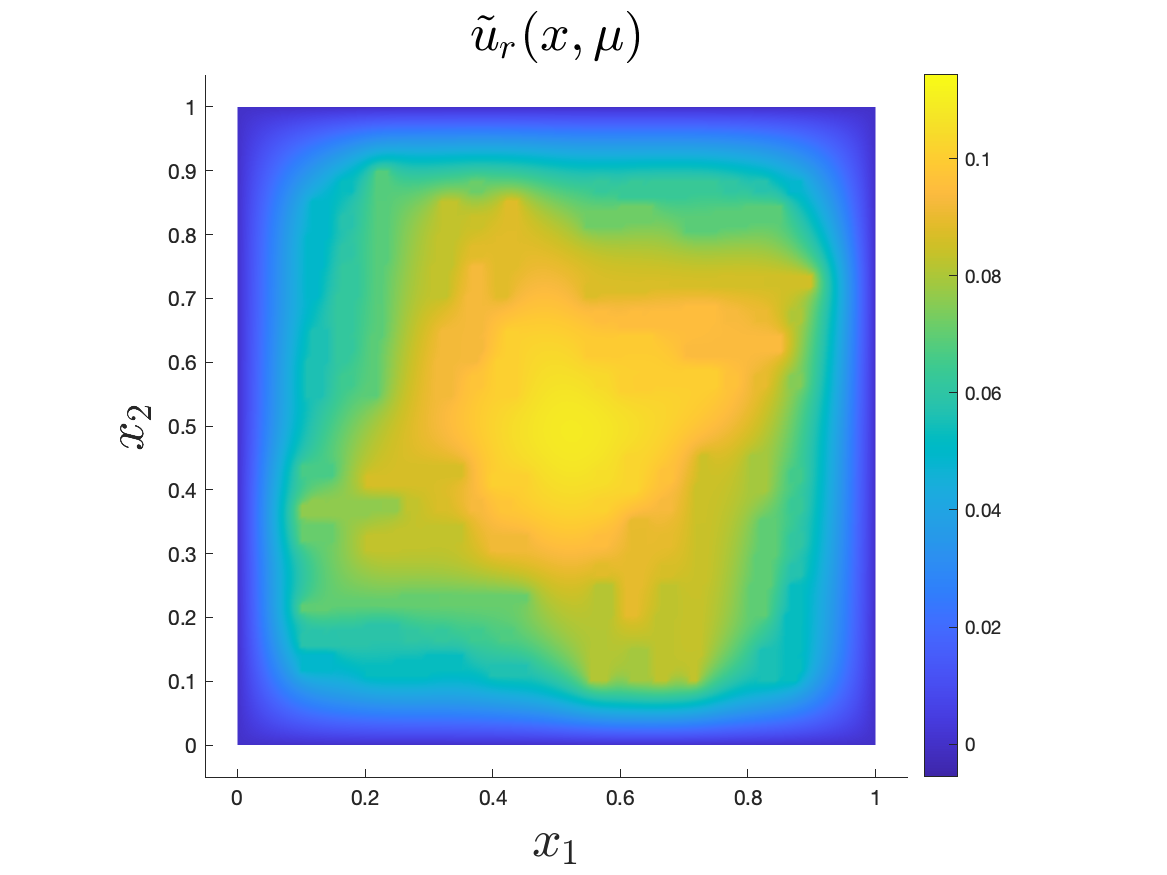}
  \caption{$\tilde{u}_{r}$}
  \label{fig:Sol_plots_ur_nrb6_2DNLHC}
 \end{subfigure}
  \vspace*{-4mm}
 \caption{Example \ref{sec:2DNLHC}: Solutions obtained by coarse model $u_l(x,\mu)$, full model $u_h(x,\mu)$ and reduced model $\tilde{u}_r(x,\mu)$ at $(\mu_1,\mu_2, \mu_3) \approx (0.644,0.435,0.311)$ with the number of reduced basis  $N_{rb}=6$, $n_L=12, n_f=1$. }
\label{figure:solutions_2DNLHigh}
  \vspace*{-3mm}
\end{figure}

Figure \ref{figure:solutions_2DNLHigh} plots the solutions obtained by the low-fidelity, high-fidelity, and our reduced models at $(\mu_1,\mu_2, \mu_3) \approx (0.644,0.435,0.311)$. 
A noticeable improvement in approximating the high-fidelity solution is achieved by our bi-fidelity model compared to the low-fidelity model.
The low-fidelity solution may lack accuracy, but it still captures the essential characteristics of the high-fidelity solution in the parameter domain. This justification supports the effectiveness of the proposed bi-fidelity method.

Table \ref{tab:2dnlhigh} lists the computation time for the offline, and online stages of our method as well as the corresponding low-, high-fidelity models. Noticeably, the low-fidelity simulations, which are nearly three orders of magnitude faster than the high-fidelity runs, contribute to the majority of the online computation time. Consequently, our proposed algorithm achieves a speedup of more than $600$ times compared to the high-fidelity simulations.

\begin{table}[ht!]
\centering
\begin{tabular}{|c||c|c|c|c|c|}
\hline
$N_{rb}$&$T_{rb}^{(off)}$ & $T_{rb}^{(on)}$ &$T_{l}$ & $T_{h}$ & $T_{h}/T_{rb}^{(on)}$   \\
\hline
2 & 29.8642s & 4.0133e-4s &2.7674e-4s &3.0350e-1s& 756 \\  
4 & 30.0624s & 4.2712e-4s & 2.7674e-4s&3.0350e-1s& 711 \\  
6 &  30.9888s  & 4.4796e-4s  & 2.7674e-4s&3.0350e-1s& 678\\ 
8 & 31.1946s  & 4.8371e-4s &2.7674e-4s& 3.0350e-1s& 627 \\  
\hline
\end{tabular}
\caption{Example \ref{sec:2DNLHC}: Computation time for offline stage $T_{rb}^{(off)}$, online stage $T_{rb}^{(on)}$, low-fidelity run $T_{l}$ and high-fidelity run $T_{h}$ with respect to the number of reduced bases; Computation times $T_{rb}^{(on)}$, $T_{l}$ and $T_{h}$ are averaged over $512$ parameters.}
\label{tab:2dnlhigh}
\end{table}


\subsection{2d nonlinear system of multiscale equations}
\label{sec:2DCPNLHC}
\begin{figure}[!hbt]
	\centering
  	  \begin{subfigure}{0.34\textwidth}
  \includegraphics[width=\textwidth]{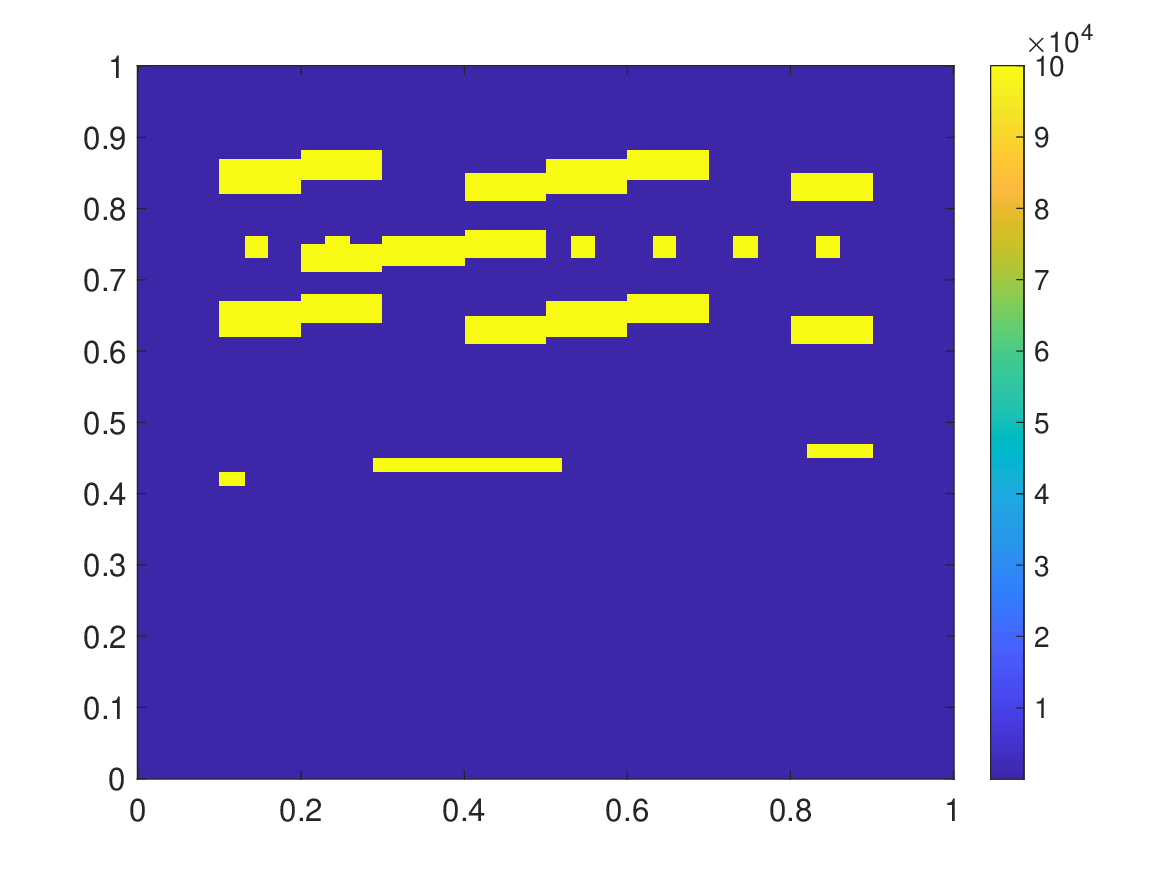}
  \caption{$\kappa_1(x,\mu)$, $\mu_1\mu_2 = 0$}
  \label{fig:k11}
 \end{subfigure}
  \hspace{-.17in}
		\begin{subfigure}{0.34\textwidth}
  \includegraphics[width=\textwidth]{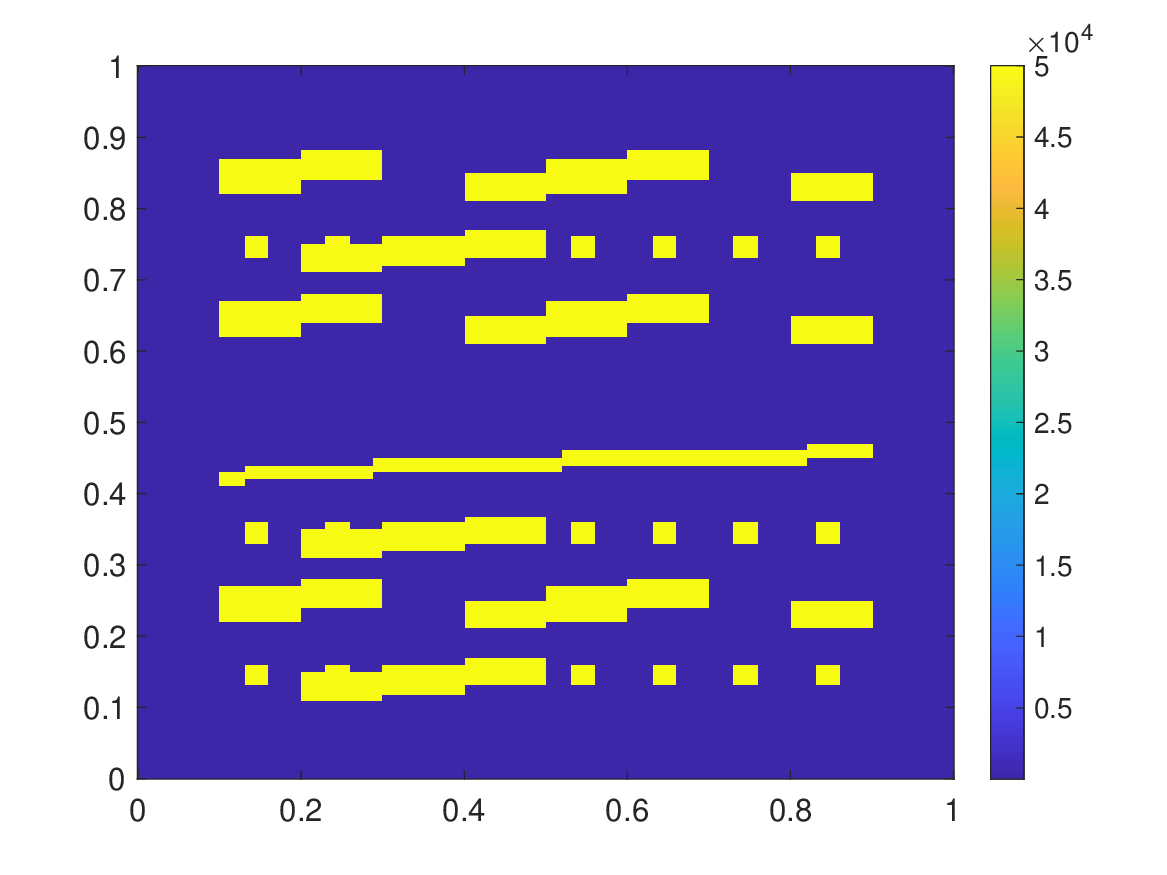}
  \caption{$\kappa_1(x,\mu)$, $\mu_1\mu_2 = 1/2$}
  \label{fig:k1}
\end{subfigure}
  \hspace{-.17in}
	\begin{subfigure}{0.34\textwidth}
  \includegraphics[width=\textwidth]{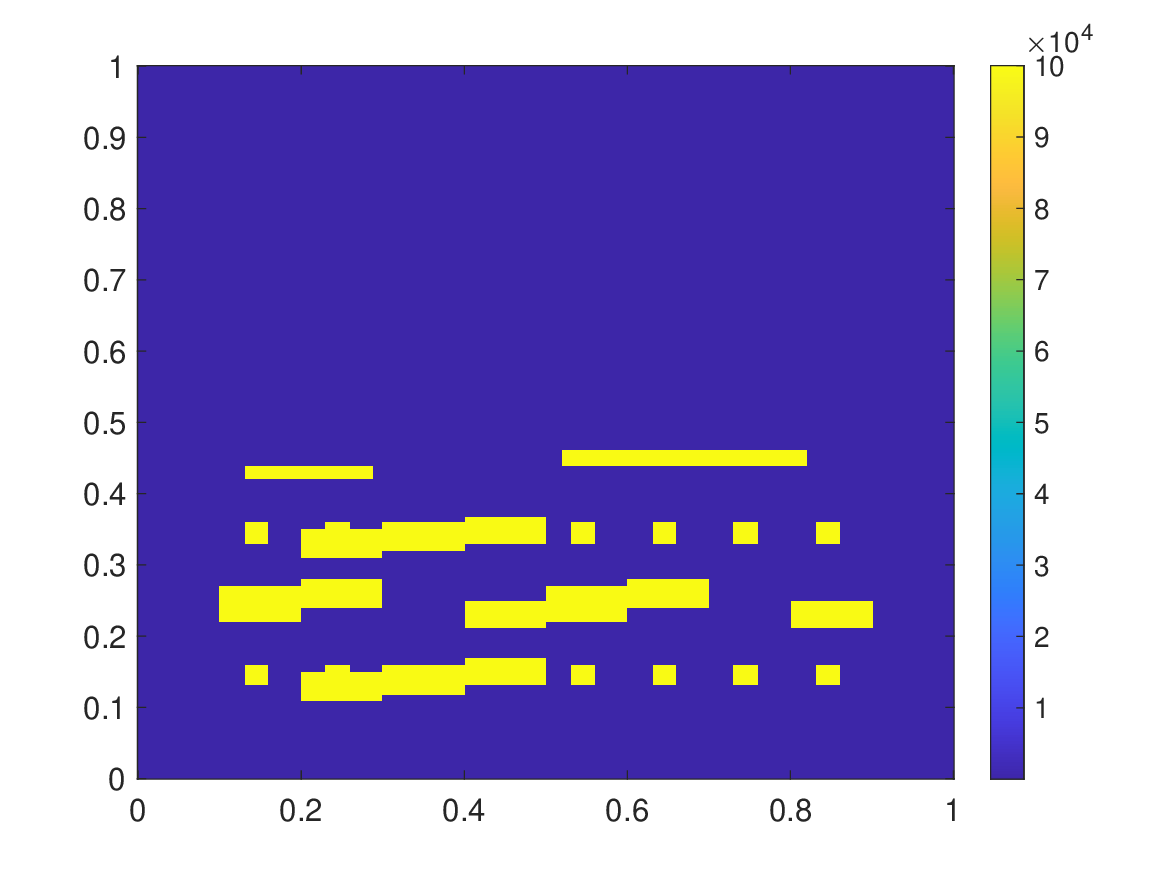}
  \caption{$\kappa_1(x,\mu)$, $\mu_1\mu_2 = 1$}
  \label{fig:k12}
\end{subfigure}
  \hspace{-.17in}
  \begin{subfigure}{0.34\textwidth}
  \includegraphics[width=\textwidth]{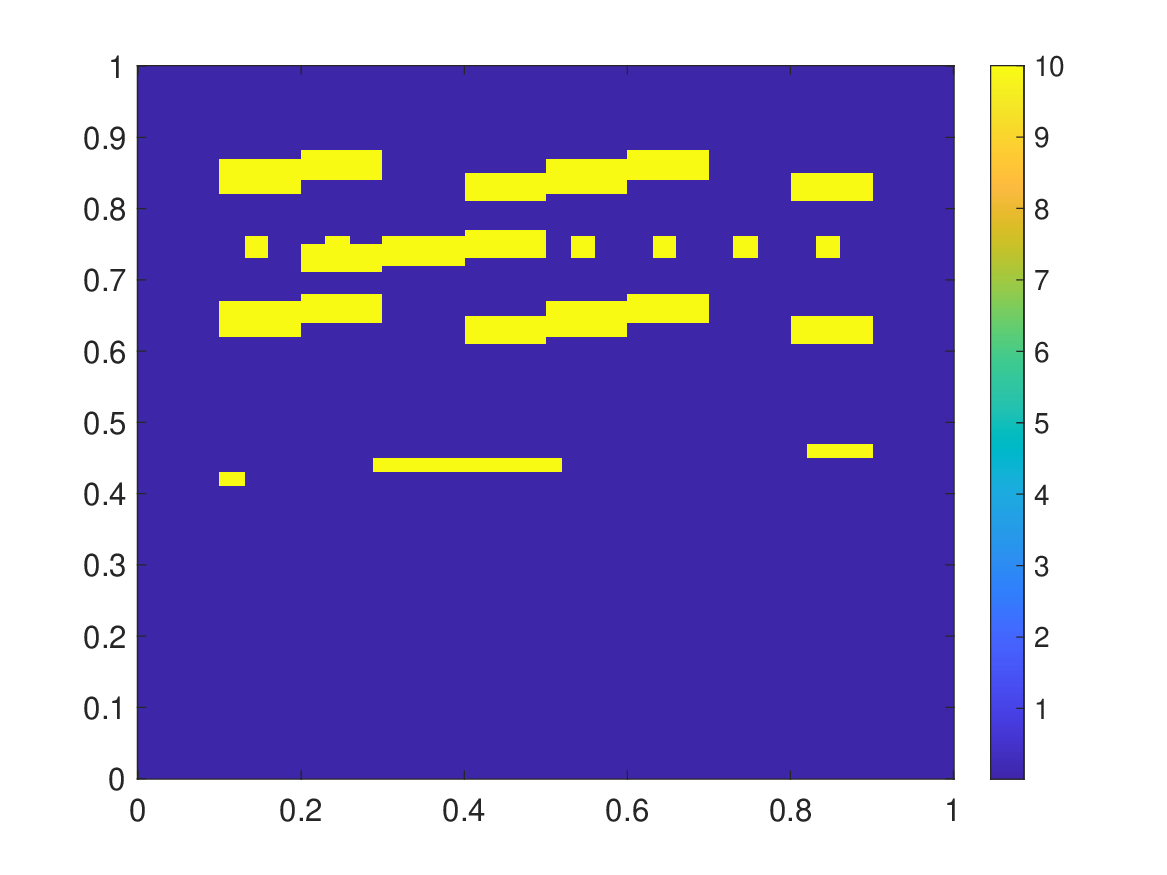}
  \caption{$\kappa_2(x,\mu)$, $\mu_1\mu_2^2 = 0$}
  \label{fig:k21}
 \end{subfigure}
  \hspace{-.17in}
  \begin{subfigure}{0.34\textwidth}
  \includegraphics[width=\textwidth]{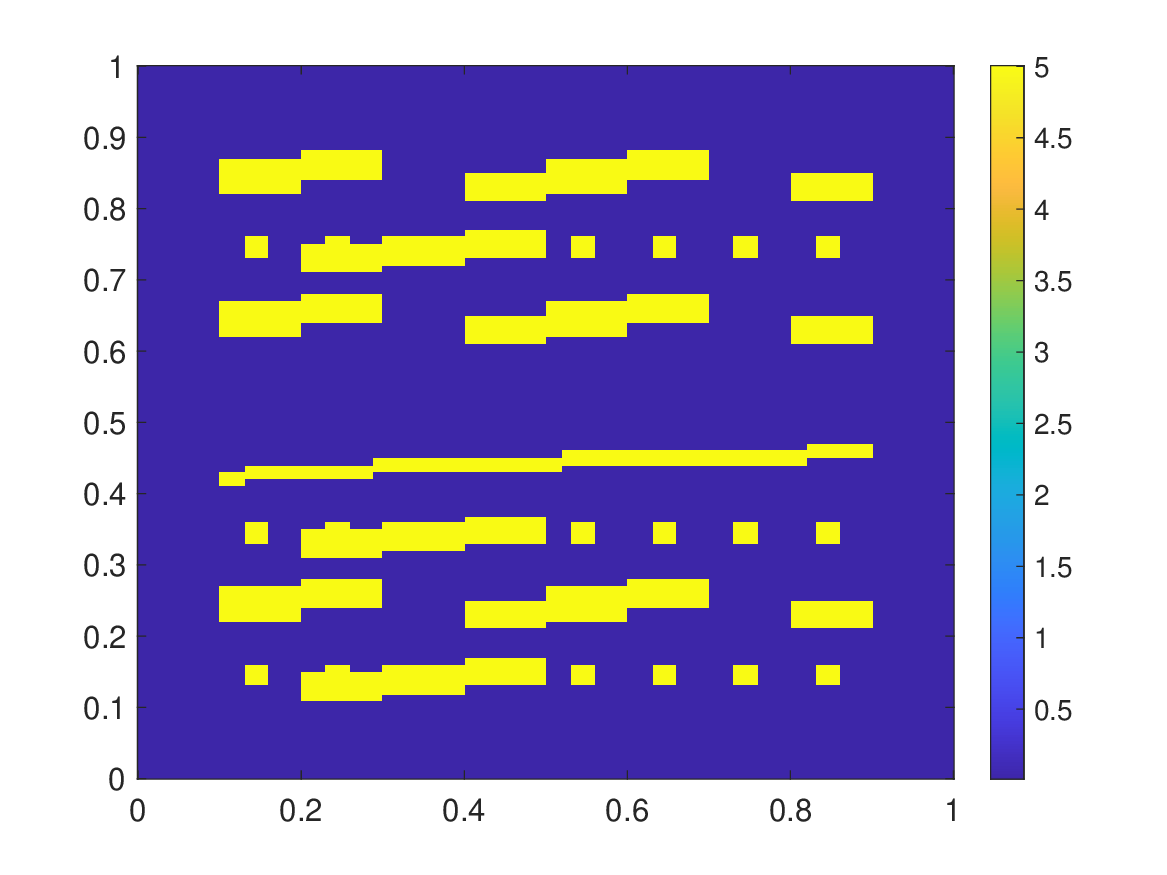}
  \caption{$\kappa_2(x,\mu)$, $\mu_1\mu_2^2 = 1/2$}
  \label{fig:k2}
 \end{subfigure}
  \hspace{-.17in}
  \begin{subfigure}{0.34\textwidth}
  \includegraphics[width=\textwidth]{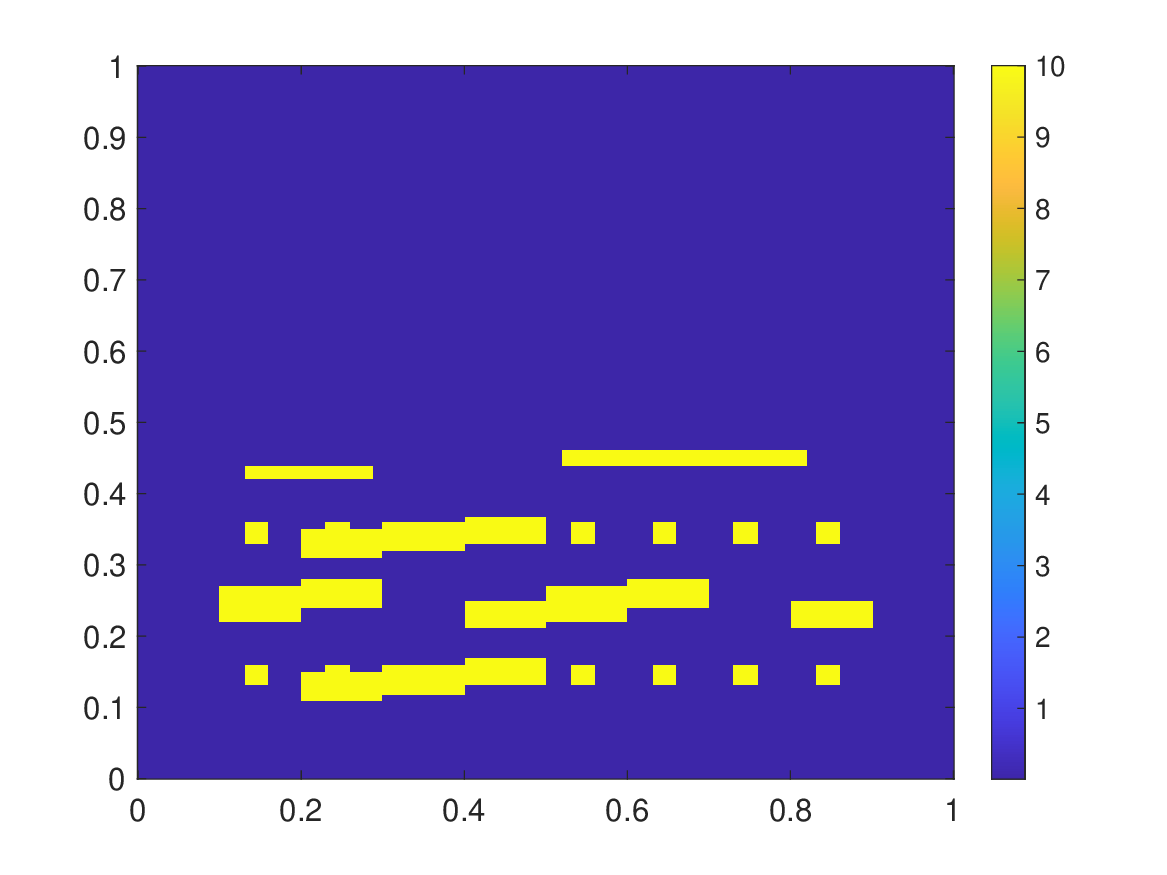}
  \caption{$\kappa_2(x,\mu)$, $\mu_1\mu_2^2 = 1$}
  \label{fig:k22}
 \end{subfigure}
 \vspace*{-4mm}
 \caption{Example \ref{sec:2DCPNLHC}: The coefficients $\kappa_1(x,\mu)$, $\kappa_2(x,\mu)$ with different values of $\mu$; The value in the background (blue region) is $1$ for all cases.}
\label{figure:highcontrast_kappas}
 \vspace*{-5mm}
\end{figure}

Motivated by an example in \cite{mai2022constraint}, we next consider the following parameterized coupled steady-state system in the domain $\Omega = [0,1]^2$. 
\beq\label{eq:nex}
\bsp
- \div \left(\frac{\kappa_1(x,\mu)}{1+|p_1(x,\mu)|} \nabla p_1(x,\mu) \right) + \frac{10^5}{1+|p_1(x,\mu)+p_2(x,\mu)|}(p_1(x,\mu)-p_2(x,\mu)) &= 1 \,,\\
- \div \left(\frac{\kappa_2(x,\mu)}{1+|p_2(x,\mu)|} \nabla p_2(x,\mu) \right) + \frac{10^5}{1+ |p_1(x,\mu)+p_2(x,\mu)|}(p_2(x,\mu)-p_1(x,\mu)) &= 1 \,,
\end{split}
\eeq
where it has zero Dirichlet boundary condition and the parameters $\mu = (\mu_1, \mu_2)$ are chosen from $[0,1]^2$. We have the high-contrast permeability fields $\kappa_1(x,\mu)$ and $\kappa_2(x,\mu)$ defined as follows:
\beq
\bsp
\kappa_1(x,\mu) = \mu_1\mu_2\kappa_{1,1}(x)+(1-\mu_1\mu_2)\kappa_{1,2}(x),\\
\kappa_2(x,\mu) = \mu_1\mu_2^2\kappa_{2,1}(x)+(1- \mu_1\mu_2^2)\kappa_{2,2}(x),
\end{split}
\eeq
where $\kappa_{i,j}(x)$ has the high-contrast features that represent inclusions, channels. This type of coefficient is introduced in \cite{efendiev2013generalized}. Figure \ref{figure:highcontrast_kappas} illustrates the coefficients $\kappa_1$, $\kappa_2$ for different parameters.
In our bi-fidelity approach, both high- and low-fidelity models are given by the same FEM solver with different grid sizes, $128\times 128$ and $4\times 4$ respectively. The FEM solver uses the Picard's iteration to handle the nonlinearity.
We consider $400$ different parameters randomly sampled from the parameter domain to construct the dense point set $\Gamma$.
We construct the reduced operator using high-fidelity operators at $n_L=12$ selected parameters, and the right-hand side is independent of the parameter, i.e., $n_f = 1$. 

\begin{figure}[!htb]
	\centering
 \begin{subfigure}{0.45\textwidth}
\includegraphics[width=\textwidth]{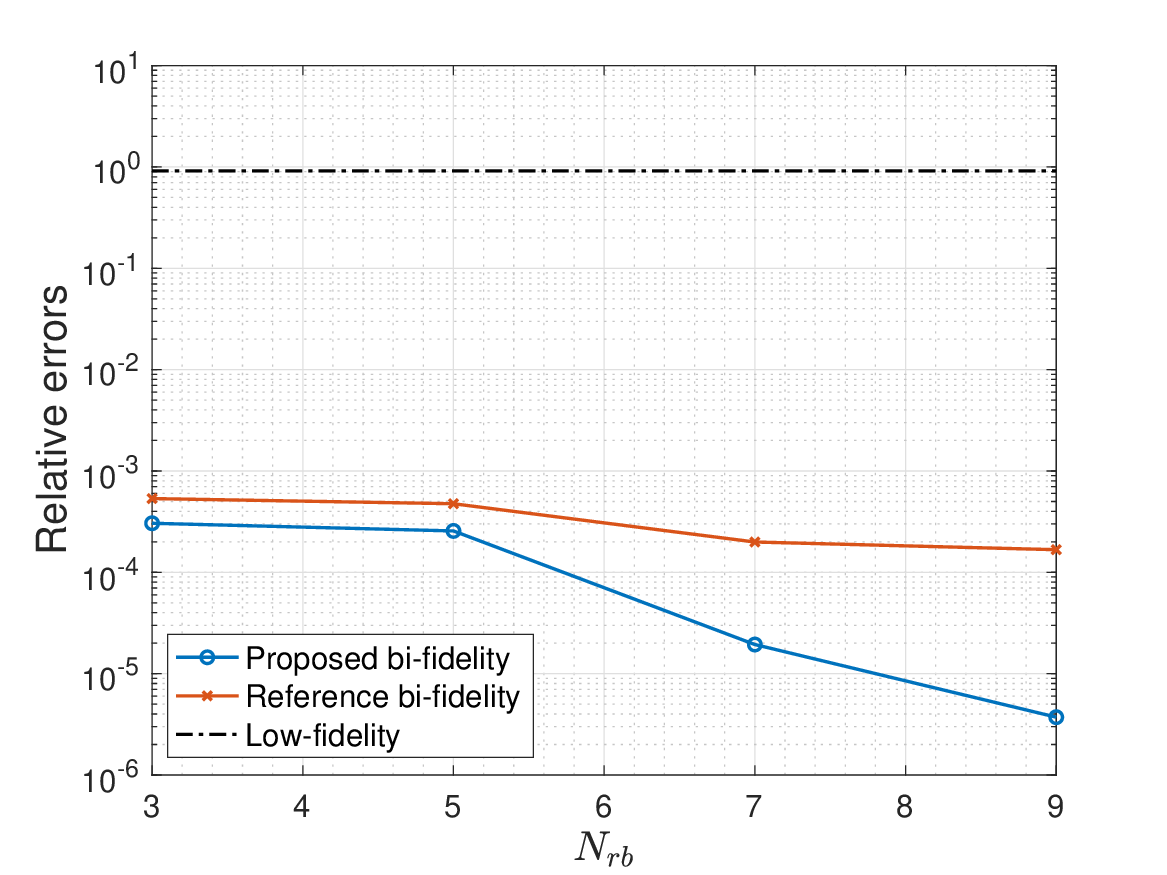}
\caption{}
    \label{fig:Errors_2DCPNLHC}
    \end{subfigure}
    \hspace{.1in}
 \begin{subfigure}{0.45\textwidth}
  \includegraphics[width=\textwidth]{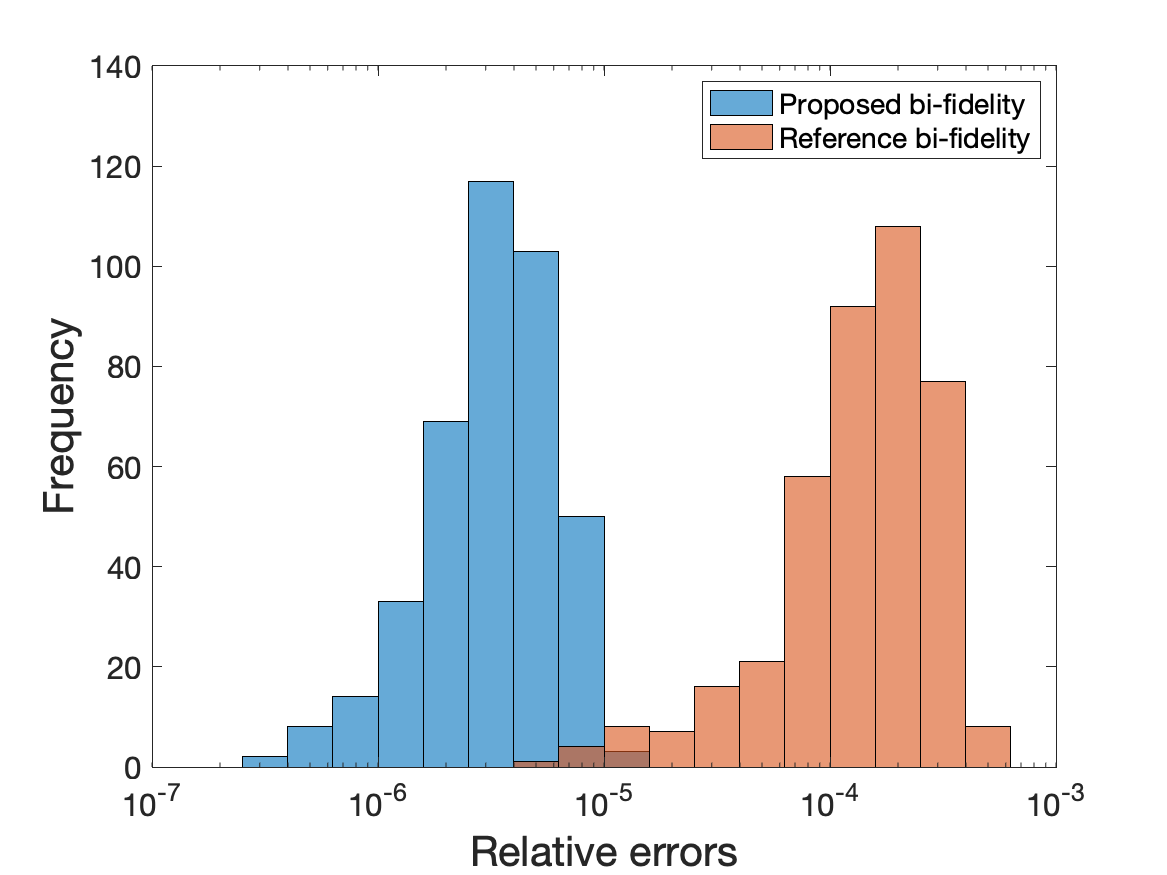}
  \caption{}
  \label{fig:Histo_2DCPNLHC}
\end{subfigure}
\vspace{-.2in}
  \caption{Example \ref{sec:2DCPNLHC}:  a. The relative errors of the reduced solution by our new bi-fidelity method and the reference bi-fidelity method with different number of reduced basis; b. Histogram of the relative errors for both new and the reference bi-fidelity methods when $N_{rb} = 9$.}
  \label{figure:Errors_2D_CPnonlinear}
 \end{figure}

Figure \ref{figure:Errors_2D_CPnonlinear} shows the averaged relative errors of the reduced solutions by the proposed and the reference bi-fidelity algorithms based on $N_{rb} = 3, 5, 7, 9$ reduced bases,  corresponding to the total number of $9$ high-fidelity runs are required for all cases. The relative error converges clearly as the number of basis increases and it reaches the level of $10^{-6}$ when $N_{rb} = 9$ reduced bases are used. 
We obtain the relative errors computed and averaged over $400$ randomly sampled parameters, where the parameter samples are independent of $\Gamma$.
We observe a significant improvement in accuracy achieved by the bi-fidelity algorithm compared to the corresponding low-fidelity model and the reference framework. This highlights the benefits of respecting the reduced order equation (\ref{eq:reduced_short}) during the online stage. This can be further supported by the error histogram of both bi-fidelity methods with $N_{rb} = 9$ in Figure \ref{figure:Errors_2D_CPnonlinear}.

\begin{figure}[!htb]
	\centering
	\begin{subfigure}{0.35\textwidth}
  \includegraphics[width=\textwidth]{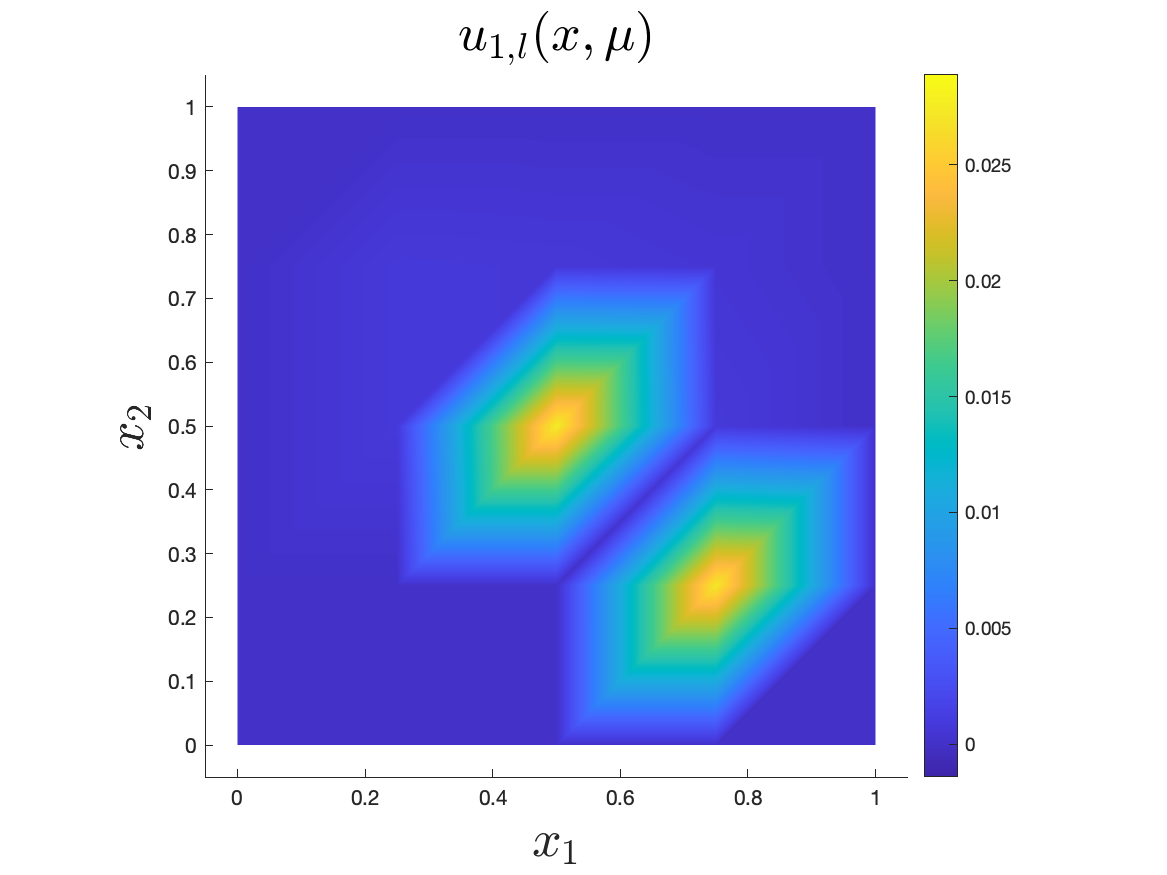}
  \caption{$u_{1,l}$}
  \label{fig:Sol_plots_ul_nrb9_2DCPNLHC}
\end{subfigure}
 \hspace{-.3in}
	\begin{subfigure}{0.35\textwidth}
  \includegraphics[width=\textwidth]{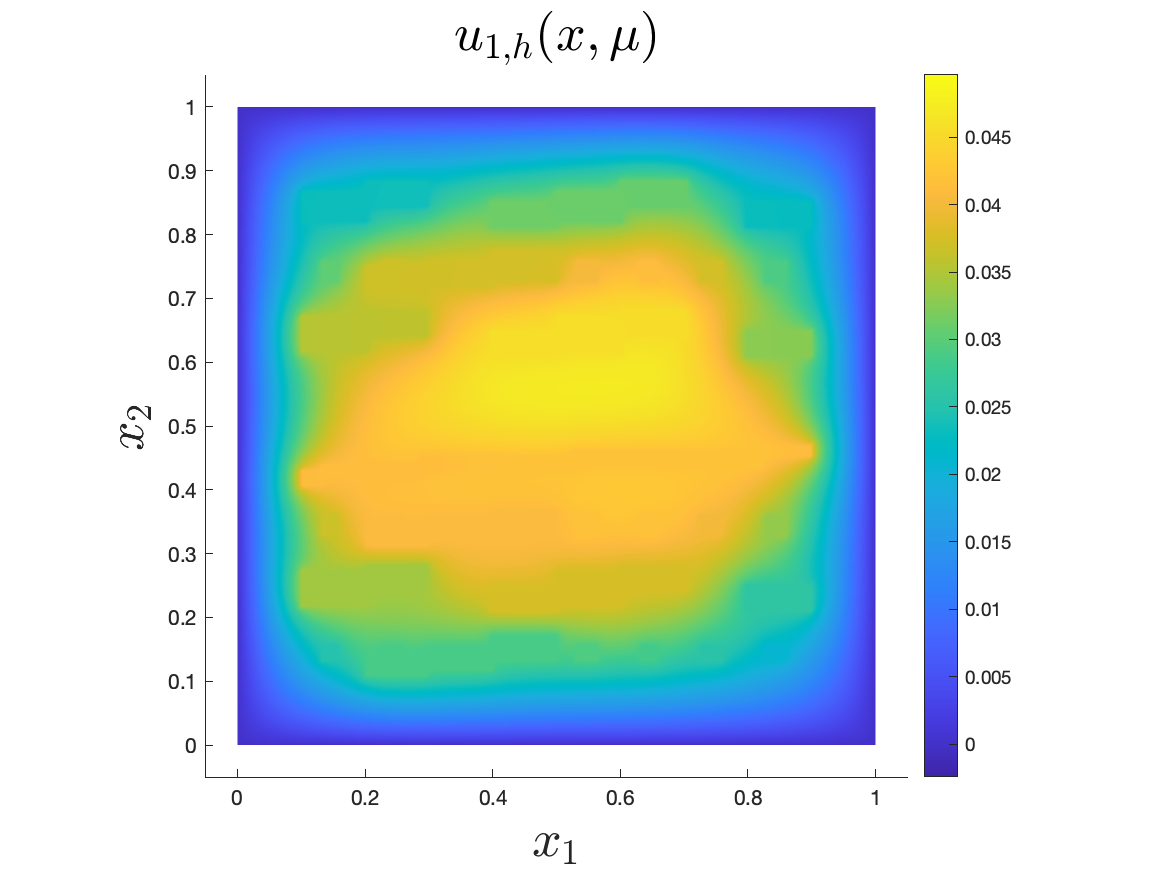}
  \caption{$u_{1,h}$}
  \label{fig:Sol_plots_uh_nrb9_2DCPNLHC}
\end{subfigure}
\hspace{-.3in}
  \begin{subfigure}{0.35\textwidth}
  \includegraphics[width=\textwidth]{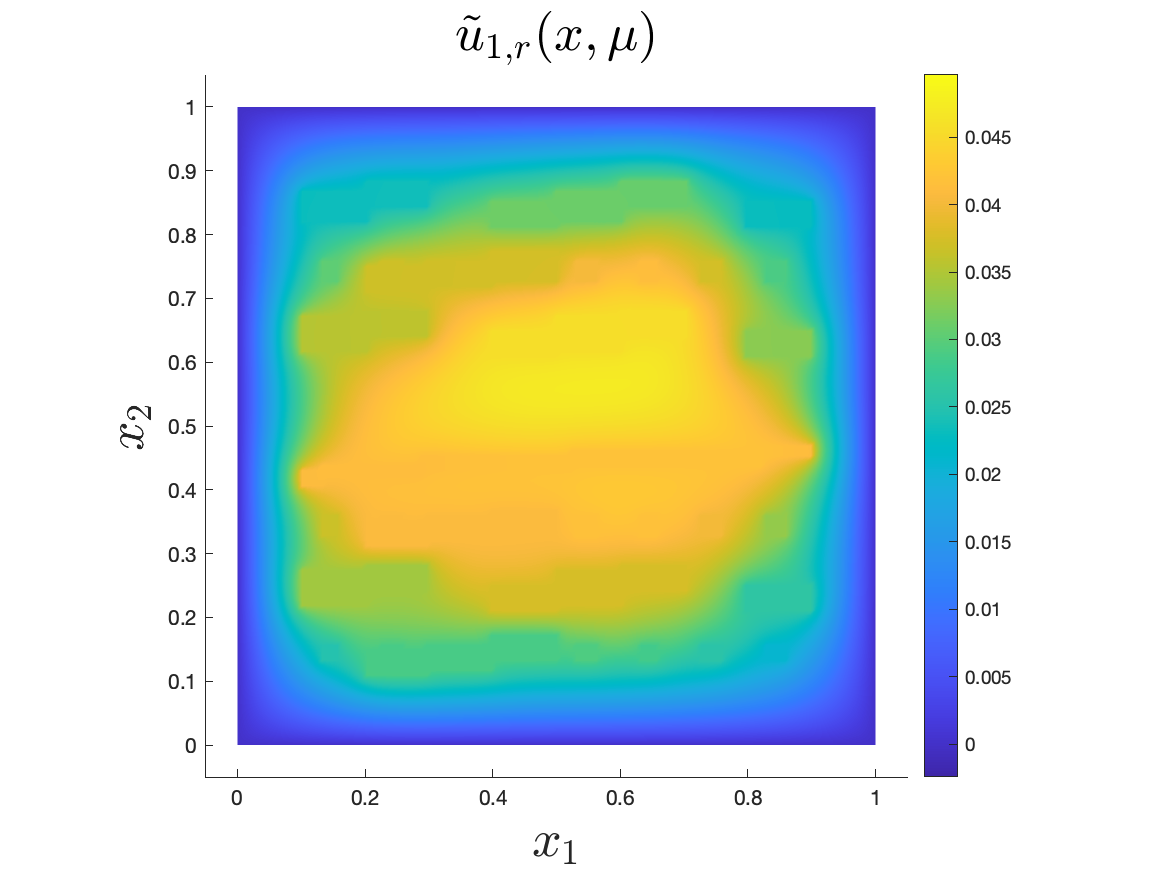}
  \caption{$\tilde{u}_{1,r}$}
  \label{fig:Sol_plots_ur_nrb9_2DCPNLHC}
 \end{subfigure}
  \vspace*{-4mm}
 \caption{Example \ref{sec:2DCPNLHC}:Solutions obtained by coarse model $u_l(x,\mu)$, full model $u_h(x,\mu)$ and reduced model $\tilde{u}_r(x,\mu)$ at $(\mu_1,\mu_2) \approx (0.326, 0.108)$ with the number of reduced basis  $N_{rb}=9$, $n_L=12, n_f=1$.}
\label{figure:solutions_2CPDNLHigh}
  \vspace*{-3mm}
\end{figure}

Figure \ref{figure:solutions_2CPDNLHigh} plots the low-fidelity and high-fidelity solutions and the proposed reduced solution at $(\mu_1,\mu_2) \approx (0.326, 0.108)$. The proposed bi-fidelity method delivered a much better solution compared to the corresponding low-fidelity solution. The result is surprising yet reasonable, suggesting
that even though the low-fidelity model may be inaccurate in the physical space, it still can
capture the behaviors and characteristics of the high-fidelity solution in the parameter space, which
is informative and crucial for the effectiveness of  our bifidelity method.

Table \ref{tab:2dnonlinearhighcontrast} presents the computation time for the offline, online stage of the proposed algorithm and the high- and low-fidelity simulations. 
The computation time of the low-fidelity model is dominant in the online stage and is more than $600$ times faster than the high-fidelity runs, which gives rise to a speedup range of $420-460$, demonstrating the efficiency of our proposed method.


\begin{table}[ht!]
\centering
\begin{tabular}{|c||c|c|c|c|c|}
\hline
$N_{rb}$&$T_{rb}^{(off)}$ & $T_{rb}^{(on)}$ &$T_{l}$ & $T_{h}$ & $T_{h}/T_{rb}^{(on)}$  \\
\hline
3 & 24.0406s &2.1559e-3s &1.4522e-3s &9.8342e-1s& 456\\ 
5 & 25.5374s &2.2328e-3s & 1.4522e-3s&9.8342e-1s&  440\\ 
7 & 27.7250s  &2.3001e-3s & 1.4522e-3s&9.8342e-1s&  428\\ 
9 & 28.6409s &2.3312e-3s & 1.4522e-3s &9.8342e-1s& 422\\ 
\hline
\end{tabular}
\caption{Example \ref{sec:2DCPNLHC}: Computation time for offline stage $T_{rb}^{(off)}$, online stage $T_{rb}^{(on)}$, low-fidelity run $T_{l}$ and high-fidelity run $T_{h}$ with respect to the number of reduced bases; Computation times $T_{rb}^{(on)}$, $T_{l}$ and $T_{h}$ are averaged over $400$ parameters.}
\label{tab:2dnonlinearhighcontrast}
\end{table}

%% file: Summary.tex
\section{Conclusions}

In this paper, we introduce a non-intrusive reduced basis technique for parameterized time-independent differential equations. Our method leverages both the efficiency of the low-fidelity model and the accuracy of the high-fidelity model to construct the reduced solution efficiently. Unlike many other non-intrusive reduced order models (ROMs), our approach enforces the reduced equation in the high-fidelity reduced approximation during the online stage. This enforcement enhances the accuracy in approximating the reduced solution for a new given parameter. Furthermore, the non-intrusive nature of our proposed method ensures straightforward implementation in practical applications.
We demonstrate the versatility and efficiency of our proposed method through fast convergence and significant speedup for various linear and non-linear nontrivial PDEs as long as the low-fidelity model can mimic the parameter dependence of high-fidelity solution/operator in the parameter space. While this study primarily focuses on the applicability and computational aspects of the technique, we acknowledge the need to develop a theoretical foundation for error estimation in future work.

%% file: Main_ROM_Coarse.bbl
\begin{thebibliography}{10}

\bibitem{barrault2004empirical}
Maxime Barrault, Yvon Maday, Ngoc~Cuong Nguyen, and Anthony~T Patera.
\newblock An ‘empirical interpolation’method: application to efficient
  reduced-basis discretization of partial differential equations.
\newblock {\em Comptes Rendus Mathematique}, 339(9):667--672, 2004.

\bibitem{bonomi2017matrix}
Diana Bonomi, Andrea Manzoni, and Alfio Quarteroni.
\newblock A matrix deim technique for model reduction of nonlinear parametrized
  problems in cardiac mechanics.
\newblock {\em Computer Methods in Applied Mechanics and Engineering},
  324:300--326, 2017.

\bibitem{chakir2009two}
Rachida Chakir and Yvon Maday.
\newblock A two-grid finite-element/reduced basis scheme for the approximation
  of the solution of parameter dependent pde.
\newblock In {\em 9e Colloque national en calcul des structures}, 2009.

\bibitem{chaturantabut2010nonlinear}
Saifon Chaturantabut and Danny~C Sorensen.
\newblock Nonlinear model reduction via discrete empirical interpolation.
\newblock {\em SIAM Journal on Scientific Computing}, 32(5):2737--2764, 2010.

\bibitem{chen2016model}
Peng Chen and Christoph Schwab.
\newblock Model order reduction methods in computational uncertainty
  quantification.
\newblock {\em Handbook of uncertainty quantification}, pages 1--53, 2016.

\bibitem{chen2021physics}
Wenqian Chen, Qian Wang, Jan~S Hesthaven, and Chuhua Zhang.
\newblock Physics-informed machine learning for reduced-order modeling of
  nonlinear problems.
\newblock {\em Journal of computational physics}, 446:110666, 2021.

\bibitem{chen2019l1}
Yanlai Chen, Sigal Gottlieb, Lijie Ji, Yvon Maday, and Zhenli Xu.
\newblock L1-roc and r2-roc: L1-and r2-based reduced over-collocation methods
  for parametrized nonlinear partial differential equations.
\newblock {\em arXiv preprint arXiv:1906.07349}, 2019.

\bibitem{efendiev2013generalized}
Yalchin Efendiev, Juan Galvis, and Thomas~Y Hou.
\newblock Generalized multiscale finite element methods (gmsfem).
\newblock {\em Journal of computational physics}, 251:116--135, 2013.

\bibitem{efendiev2014generalized}
Yalchin Efendiev, Juan Galvis, Guanglian Li, and Michael Presho.
\newblock Generalized multiscale finite element methods. nonlinear elliptic
  equations.
\newblock {\em Communications in Computational Physics}, 15(3):733--755, 2014.

\bibitem{efendiev2012systematic}
Yalchin Efendiev, Juan Galvis, and Florian Thomines.
\newblock A systematic coarse-scale model reduction technique for
  parameter-dependent flows in highly heterogeneous media and its applications.
\newblock {\em Multiscale Modeling \& Simulation}, 10(4):1317--1343, 2012.

\bibitem{etter2023coarse}
Philip~A Etter, Yuwei Fan, and Lexing Ying.
\newblock Coarse-proxy reduced basis methods for integral equations.
\newblock {\em Journal of Computational Physics}, 475:111835, 2023.

\bibitem{forrester2007multi}
Alexander~IJ Forrester, Andr{\'a}s S{\'o}bester, and Andy~J Keane.
\newblock Multi-fidelity optimization via surrogate modelling.
\newblock {\em Proceedings of the royal society a: mathematical, physical and
  engineering sciences}, 463(2088):3251--3269, 2007.

\bibitem{guo2018reduced}
Mengwu Guo and Jan~S Hesthaven.
\newblock Reduced order modeling for nonlinear structural analysis using
  gaussian process regression.
\newblock {\em Computer methods in applied mechanics and engineering},
  341:807--826, 2018.

\bibitem{hesthaven2016certified}
Jan~S Hesthaven, Gianluigi Rozza, Benjamin Stamm, et~al.
\newblock {\em Certified reduced basis methods for parametrized partial
  differential equations}, volume 590.
\newblock Springer, 2016.

\bibitem{hesthaven2018non}
Jan~S Hesthaven and Stefano Ubbiali.
\newblock Non-intrusive reduced order modeling of nonlinear problems using
  neural networks.
\newblock {\em Journal of Computational Physics}, 363:55--78, 2018.

\bibitem{kast2020non}
Mariella Kast, Mengwu Guo, and Jan~S Hesthaven.
\newblock A non-intrusive multifidelity method for the reduced order modeling
  of nonlinear problems.
\newblock {\em Computer Methods in Applied Mechanics and Engineering},
  364:112947, 2020.

\bibitem{langtangen2016solving}
Hans~Petter Langtangen.
\newblock Solving nonlinear ode and pde problems.
\newblock {\em Center for Biomedical Computing, Simula Research Laboratory and
  Department of Informatics, University of Oslo}, 2016.

\bibitem{larson2013finite}
Mats~G Larson and Fredrik Bengzon.
\newblock {\em The finite element method: theory, implementation, and
  applications}, volume~10.
\newblock Springer Science \& Business Media, 2013.

\bibitem{lu2021bifidelity}
Chuan Lu and Xueyu Zhu.
\newblock Bifidelity data-assisted neural networks in nonintrusive
  reduced-order modeling.
\newblock {\em Journal of Scientific Computing}, 87(1):1--30, 2021.

\bibitem{mai2022constraint}
Tina Mai, Siu~Wun Cheung, and Jun Sur~Richard Park.
\newblock Constraint energy minimizing generalized multiscale finite element
  method for multi-continuum richards equations.
\newblock {\em arXiv preprint arXiv:2205.11294}, 2022.

\bibitem{mai2023constraint}
Tina Mai, Siu~Wun Cheung, and Jun Sur~Richard Park.
\newblock Constraint energy minimizing generalized multiscale finite element
  method for multi-continuum richards equations.
\newblock {\em Journal of Computational Physics}, page 111915, 2023.

\bibitem{narayan2014stochastic}
Akil Narayan, Claude Gittelson, and Dongbin Xiu.
\newblock A stochastic collocation algorithm with multifidelity models.
\newblock {\em SIAM Journal on Scientific Computing}, 36(2):A495--A521, 2014.

\bibitem{negri2015efficient}
Federico Negri, Andrea Manzoni, and David Amsallem.
\newblock Efficient model reduction of parametrized systems by matrix discrete
  empirical interpolation.
\newblock {\em Journal of Computational Physics}, 303:431--454, 2015.

\bibitem{ou2020low}
Na~Ou, Guang Lin, and Lijian Jiang.
\newblock A low-rank approximated multiscale method for pdes with random
  coefficients.
\newblock {\em Multiscale Modeling \& Simulation}, 18(4):1595--1620, 2020.

\bibitem{quarteroni2015reduced}
Alfio Quarteroni, Andrea Manzoni, and Federico Negri.
\newblock {\em Reduced basis methods for partial differential equations: an
  introduction}, volume~92.
\newblock Springer, 2015.

\bibitem{rozza2008reduced}
Gianluigi Rozza, Dinh Bao~Phuong Huynh, and Anthony~T Patera.
\newblock Reduced basis approximation and a posteriori error estimation for
  affinely parametrized elliptic coercive partial differential equations:
  application to transport and continuum mechanics.
\newblock {\em Archives of Computational Methods in Engineering}, 15(3):229,
  2008.

\bibitem{wang2019non}
Qian Wang, Jan~S Hesthaven, and Deep Ray.
\newblock Non-intrusive reduced order modeling of unsteady flows using
  artificial neural networks with application to a combustion problem.
\newblock {\em Journal of computational physics}, 384:289--307, 2019.

\bibitem{xiao2017parameterized}
D~Xiao, F~Fang, CC~Pain, and IM~Navon.
\newblock A parameterized non-intrusive reduced order model and error analysis
  for general time-dependent nonlinear partial differential equations and its
  applications.
\newblock {\em Computer Methods in Applied Mechanics and Engineering},
  317:868--889, 2017.

\bibitem{xiao2015non}
Dunhui Xiao, Fangxin Fang, Christopher Pain, and Guangwei Hu.
\newblock Non-intrusive reduced-order modelling of the navier--stokes equations
  based on rbf interpolation.
\newblock {\em International Journal for Numerical Methods in Fluids},
  79(11):580--595, 2015.

\bibitem{zhu2017multi}
Xueyu Zhu, Erin~M Linebarger, and Dongbin Xiu.
\newblock Multi-fidelity stochastic collocation method for computation of
  statistical moments.
\newblock {\em Journal of Computational Physics}, 341:386--396, 2017.

\bibitem{zhu2014computational}
Xueyu Zhu, Akil Narayan, and Dongbin Xiu.
\newblock Computational aspects of stochastic collocation with multifidelity
  models.
\newblock {\em SIAM/ASA Journal on Uncertainty Quantification}, 2(1):444--463,
  2014.

\end{thebibliography}
